%% file: leib.tex
\input amstexl
\input lamstex

\font\blackital=cmmib10  \skewchar\blackital='177
\font\sblackital=cmmib7  \skewchar\sblackital='177
\font\ssblackital=cmmib5  \skewchar\ssblackital='177
\font\sanss=cmss10
\font\ssanss=cmss8 scaled 900
\font\sssanss=cmss8 scaled 600

\font\teneurm=eurm10

\font\blackboard=msbm10
\font\sblackboard=msbm7
\font\ssblackboard=msbm5
\font\caligr=eusm10
\font\scaligr=eusm7
\font\sscaligr=eusm5

\font\fraktur=eufm10
\font\sfraktur=eufm7
\font\ssfraktur=eufm5

\newfam\bbfam
\textfont\bbfam=\blackboard
\scriptfont\bbfam=\sblackboard
\scriptscriptfont\bbfam=\ssblackboard
\def\bb#1{{\fam\bbfam\relax#1}}

\newfam\ssfam
\textfont\ssfam=\sanss
\scriptfont\ssfam=\ssanss
\scriptscriptfont\ssfam=\sssanss
\def\ss#1{{\fam\ssfam\relax#1}}

\newfam\clfam
\textfont\clfam=\caligr
\scriptfont\clfam=\scaligr
\scriptscriptfont\clfam=\sscaligr
\def\cl#1{{\fam\clfam\relax#1}}

\newfam\frfam
\textfont\frfam=\fraktur
\scriptfont\frfam=\sfraktur
\scriptscriptfont\frfam=\ssfraktur
\def\fr#1{{\fam\frfam\relax#1}}
\newfam\gpfam
\def\eu#1{\hbox{$\fam\gpfam\relax#1\textfont\gpfam=\teneurm$}}

\catcode`\"=12
\mathchardef\za="710B  
\mathchardef\zb="710C  
\mathchardef\zg="710D  
\mathchardef\zd="710E  
\mathchardef\zve="710F 
\mathchardef\zz="7110  
\mathchardef\zh="7111  
\mathchardef\zvy="7112 
\mathchardef\zi="7113  
\mathchardef\zk="7114  
\mathchardef\zl="7115  
\mathchardef\zm="7116  
\mathchardef\zn="7117  
\mathchardef\zx="7118  
\mathchardef\zp="7119  
\mathchardef\zr="711A  
\mathchardef\zs="711B  
\mathchardef\zt="711C  
\mathchardef\zu="711D  
\mathchardef\zvf="711E 
\mathchardef\zq="711F  
\mathchardef\zc="7120  
\mathchardef\zw="7121  
\mathchardef\ze="7122  
\mathchardef\zy="7123  
\mathchardef\zvp="7124  
\mathchardef\zvr="7125 
\mathchardef\zvs="7126 
\mathchardef\zf="7127  
\mathchardef\zG="7000  
\mathchardef\zD="7001  
\mathchardef\zY="7002  
\mathchardef\zL="7003  
\mathchardef\zX="7004  
\mathchardef\zP="7005  
\mathchardef\zS="7006  
\mathchardef\zU="7007  
\mathchardef\zF="7008  
\mathchardef\zC="7009  
\mathchardef\zW="700A  
\catcode`\"=\active

\loadmsam
\def\eza{\eu\za}

\def\eze{\eu\ze}

\def\ezk{\eu\zk}
\def\ezt{\eu\zt}

\def\ezp{\eu\zp}
\def\sT{{\ss T}}
\def\cC{{\cl C}}
\def\cR{{\cl R}}
\def\cT{{\cl T}}
\def\ts{{\fr T}}

\TagsOnRight
\def\*{{\textstyle *}}
\def\s*{{\scriptstyle *}}
\font\ssans=cmss8
\font\smss=cmss8 scaled 800
\def\sss#1{\hbox{\ssanss #1}}
\def\R{{\bb R}}

\def\Z{{\bb Z}}

\def\ssT{\sss T}
\def\smT{\hbox{\smss T}}

\def\xd{\operatorname{d}\!}
\def\der{\operatorname{Der}}
\def\ix{\operatorname{i}}
\def\vta{\operatorname{v}_\zt}
\def\vtb{\operatorname{v}_\zp}

\def\sdt{\text{\sevenrm d}_{\smT}}
\def\dt{\xd_{\sss T}}
\def\vt{\text{\tenrm v}_{\sss T}}

\def\idt{\operatorname{id}}

\def\ix{\operatorname{i}}

\def\xD{\operatorname{ D}}
\def\ll{\text{\it \$}}

\newsymbol\leqs 1336
\newsymbol\geqs 133E
\newsymbol\rarrows 1313
\def\proof{\demo{Proof}}
\def\endproof{\hfill \vrule height4pt width6pt depth2pt \enddemo}

\input paper.st\relax
\hsize=33pc
\hoffset=30pt
\vsize=54pc
\voffset=6pt

\RefWarnings
    \title
                Algebroids -- general differential calculi on vector bundles
                \endtitle
        \thanks{Supported by KBN, grant No 2 PO3A 074 10}
    \author
                J. Grabowski$^1$ and
                P. Urba\'nski$^2$
        \endauthor
        \affil
     $^1$ Institute of Mathematics, University of Warsaw\\
        Banacha 2, 02-097 Warsaw, Poland\\
     {\sanss email:\ jagrab\@mimuw.edu.pl}\\
        \\
        $^2$    Division of Mathematical Methods in Physics,
                University of Warsaw \\
                Ho\D{z}a 74, 00-682 Warsaw, Poland\\
        {\sanss email:\ urbanski\@fuw.edu.pl}
        \endaffil
     \address{
     Division of Mathematical Methods in Physics, University of
     Warsaw, Ho\D{z}a 74, 00-682 Warsaw, Poland }

     \address{
     Institute of Mathematics, University of Warsaw, Banacha 2,
     02-097 Warsaw, Poland}



      \thanks{Supported by KBN, grant No 2 PO3A 074 10}
   \abstract
A notion of an algebroid - a generalization of a Lie
algebroid structure on a vector bundle is introduced.
        We show that many objects  of the differential calculus on a manifold
$M$ associated with the canonical Lie algebroid structure on $\hbox{\ssans T}
M$ can be obtained in the  framework of a general algebroid. Also a
compatibility condition which leads,  in general, to a concept of a bialgebroid.

   \endabstract

        \document
    \maketitle

\hl1"0"{Introduction}
The classical Cartan differential calculus on a manifold $M$, including the
exterior derivative $\xd\,$, the Lie derivative $\ll$, etc., can be viewed as
being associated with the canonical Lie algebroid structure on $\sT M$
represented by the Lie bracket of vector fields. Lie algebroids have been
introduced repeatedly into differential geometry since the early 1950s, and
also into physics and algebra, under a wide variety of names. They have been
also recognized as infinitesimal objects for Lie groupoids ([\Ref{Pr}]). We
refer to [\Ref{Ma}] for basic definitions, examples, and an extensive list
of publications in these directions.

        Being related to many areas of geometry, like connection theory,
cohomology theory, invariants of foliations and pseudogroups, symplectic and
Poisson geometry, etc., Lie algebroids became recently an object of extensive
studies.

        What we propose in this paper is to find out what are, in fact, the
structures responsible for the presence of a version of the Cartan
differential calculus on a vector bundle and how are they related. This
leads to the notion of a general algebroid.

  It is well known that there exists a one-one correspondence between Lie
algebroid
structures on a vector bundle $\zt\colon E\rightarrow M$ and linear
Poisson structures on the dual vector bundle $\zp\colon E^\*\rightarrow
M$. This correspondence can be extended to much wider class of binary
operations (brackets) on sections of $\zt$ on one side, and linear
contravariant 2-tensor fields on $E^\*$ on the other side. It is not
necessary for these  operations to be  skew-symmetric or to satisfy  the
Jacobi identity. The vector bundle $\zt$ together with  a
bracket operation, or the equivalent contravariant 2-tensor field, will be
called an algebroid. This terminology is justified by the fact that
contravariant 2-tensor fields define   certain binary operations on the space
$\cC^\infty(E^\*)$. The algebroids constructed in this way include all
finite-dimensional algebras over real numbers (e.g. associative, Jordan,
etc.) as particular examples. The base manifold $M$ is in these cases  a
single point.

     Searching for structures which give us differential calculi on vector
bundles, we look at objects of analytical mechanics as related to the
Lie algebroid structure of the tangent bundle.

  The tangent bundle $\sT M$ is the canonical Lie algebroid associated with
the canonical Poisson tensor (symplectic form) on $\sT^\*M$. Other
canonical objects associated with $\sT M$ are: the canonical isomorphism
    $$ \eza_M\colon \sT\sT^\* M \longrightarrow  \sT^\* \sT M
                              $$
of double vector bundles, discovered by Tulczyjew [\Ref{Tu}], dual to the
well-known flip
    $$ \ezk_M \colon \sT\sT M\longrightarrow \sT \sT M,
                              $$ and the tangent lift $\dt$ of tensor
fields on $M$ to tensor fields on $\sT M$ (cf. [\Ref{YI}], [\Ref{PT}],
[\Ref{GU1}]).

        The algebroid structure of the tangent bundle is not used in
analytical mechanics directly via its Lie bracket of vector fields.
In the Lagrange formulation of the infinitesimal dynamics one uses
the mentioned isomorphism $\eza_M$. The canonical symplectic
structure is the basic for
the Hamiltonian mechanics, but to obtain Euler-Lagrange equation we use
the complete tangent lift $\dt$. Following Weinstein ([\Ref{We}])
and  Libermann  ([\Ref{Li}])  also  Lie  algebroids   other   than
that of the tangent bundle can be used in the
variational formulation of the dynamics. It is important to know if one can
characterize an algebroid structure directly in terms of objects like
$\dt$, $\eza_M$, etc.

        The search for a basis of the Cartan differential calculus on a vector
bundle is not the only reason for study of general algebroids.
     We are motivated also by the fact that the interest in
Lie-like, but  not skew-symmetric, or not satisfying the Jacobi
identity, structures has been growing in last years.

     For example, structures  more  general  than  Lie  bialgebras
appear as semi-classical limits of quasi-Hopf algebras (Drinfeld in
[\Ref{Dr}]). The general theory of such objects, known  under  the
name of Jacobian quasi-bialgebras, was developed by Kosmann--Schwarzbach
([\Ref{KS1}]). They form the infinitesimal  part  of  quasi-Poisson
Lie groups. The dual objects are  Poisson  quasi-groups,  so  that
their infinitesimal parts are quasi-Lie algebras, i.e. `not  quite
Lie algebras'. We are sure all this has an algebroid counterpart  in
the spirit of our work.
Some examples are sketched in the last section.

     Also in the theory of webs one can find  objects  similar  to
Lie groupoids but not associative in general. We think that  their
infinitesimal parts are algebroids in our sense.

   Let us also mention the concept of Loday algebras, i.e. Leibniz
algebras  in  the  sense  of  Loday   [\Ref{Lo}]   (cf. also
[\Ref{KS2}])  which  are  `non-skew-symmetric  Lie  algebras'.

     Finally, let us  recall  the  use  of  Nijenhuis  tensors  in
defining deformed Lie algebroid  structures  on  tangent  bundles.
Considering properties of $(1,1)$-tensors  less  restrictive  than
vanishing of the Nijenhuis torsion, we get algebroids in the sense
we propose.

       All this provides another motivation for our work.

     The paper is organized as follows.

  In Sections~1,~2,~3 we show that objects similar to $\eza_M$,
$\ezk_M$, and $\dt$ can be associated with any algebroid structure. Each
of these objects separately provides a complete description of the
algebroid structure.
  Other constructions, known for Lie algebroids are extended to
general algebroids.  We define the tangent and the cotangent lifts of an
algebroid in Section~4 and the Lie derivative in Section~5. We find an
operation which acts as the exterior derivative in a limited way.
Conditions for an algebroid to be a Lie algebroid are formulated in terms
of  these objects (Theorems 8 and 10).

  The discussion of algebroids is extended to bialgebroids in Section~6.
Alternative definitions  of bialgebroids are considered.    A link  to the
original definition of Mackenzie and Xu  ([\Ref{MX1}]) is established by
Theorem~13.

        In Sections~7,~8 we discuss the algebroid constructions in the
important
case of an algebroid defined by a linear connection on the tangent bundle
$\sT M$.  The Lie derivative for this algebroid coincides with the
covariant derivative. The algebroid of a linear connection is not
skew-symmetric and, consequently, is not a Lie algebroid.

 The canonical example of an algebroid is the Lie algebroid $A(G)$ of a
Lie groupoid $G$. The Lie algebroid structure on $A(G)$ is obtained from
the canonical Lie algebroid structure of the tangent bundle $\sT G$. Thus,
one can expect that objects of the Lie algebroid $A(G)$ can be obtained
directly from the corresponding objects of the canonical algebroid $\sT
G$. We show in Section 9 how  the algebroid lift can be obtained from the
tangent lift $\dt$ and we discuss also the relation between the
lift $\dt$ and the lifting of multiplicative vector fields,
described  in  [\Ref{MX2}].  All  this  can  be  done  (at   least
partially) for pre- (or quasi- ; the terminology is  not
fixed yet) Lie groupoids $G$, which are structures weaker than Lie
groupoids (we do not assume associativity). Having  in  mind  some
natural examples of this kind related to quasi-Poisson Lie groups,
webs ([\Ref{Du}]), quasi-Nijenhuis tensors, etc., we are  planning  to  discuss
these problems in a further publication.

        \hl1"0.1"{Notation}
        Let $M$ be a smooth manifold. We denote by $\ezt_M \colon \sT M
\rightarrow M$ the tangent vector bundle and by $\ezp_M \colon \sT^\*
M\rightarrow M$ the cotangent vector bundle.

        Let $\zt\colon E \rightarrow M$ be a vector bundle and let $\zp
\colon E^\* \rightarrow M$ be the dual bundle. We use the following
notation for tensor bundles:
                $$\align
        \zt^{\otimes k}&\colon E\otimes_ M \cdots \otimes _M E  = \otimes
^k_M(E) \longrightarrow M , \\
                \zt^{\wedge  k} &\colon E\wedge _ M \cdots \wedge  _M E=
\wedge ^{k}_M (E) \longrightarrow M ,
            \endalign$$
         the module of sections over $\cC^\infty(M)$:
                $$\align
        \otimes ^k(\zt)&= \zG(\otimes ^k_M(E)),  \\
        \zF^k(\zt)&= \zG(\wedge ^k_M(E)),
                                                              \endalign$$
and the corresponding tensor algebras
                $$\align
        \otimes (\zt)&= \oplus _{k\in \Z} \otimes(\zt) ^k(E),  \\
                \zF(\zt)&=  \oplus_{k\in \Z} \zF^k(\zt), \\
          \ts(\zt)  &= \oplus_{k\in \Z} \otimes ^k(\zt \oplus \zp)
          = \oplus_{k\in \Z}\ts^k(\zt),  
                                                               \endalign$$
where
    $$ \zt\oplus \zp \colon E\oplus E^\* \longrightarrow M
                              $$
is the Whitney sum. In particular,
    $$\otimes ^0(\zt) = \zF^0(\zt) = \ts^0(\zt) =\cC^\infty (M)
                              $$
  and
    $$ \otimes ^k(\zt)= \zF^k(\zt) = \ts^k(\zt) = \{0\} \ \ \text{for
$k<0$}.
                              $$

  By $\langle \cdot ,\cdot \rangle $, we denote the canonical pairing
between $E$ and $E^\*$ as well as pairings between  the corresponding  tensor
bundles, e.g.,
    $$ \langle \cdot ,\cdot \rangle \colon \otimes ^k_M (E) \times_M
\otimes ^k_M (E^\*) \longrightarrow \R,
                              $$
  and pairings of sections, e.~g.,
    $$ \langle \cdot ,\cdot \rangle \colon \otimes ^k (\zt) \times
\otimes ^k (\zp) \longrightarrow \cC^\infty (M).
                              $$
  Let  $K$ be a section of a tensor bundle $\otimes ^k_M(E)$, $K \in
\otimes ^k(\zt)$. We denote by $\zi (K)$ the corresponding linear function
on the dual bundle
    $$\align
\zi(K)  &\colon \otimes^ k_M(E^\*) \rightarrow  \R \\
  &  \colon a  \mapsto \langle K(m), a\rangle, \ \ \  m = \zp^{\otimes k}(a).
                    \tag \label{F0.1}\endalign$$
  For a section $X$ of $\zt$ ($X\in \otimes ^1(\zt) $), we have the usual
operators of insertion
    $$\align
\ix^l_X  &\colon \otimes ^{k+1}(\zp)\rightarrow \otimes ^k(\zp)  \\
  &\colon \zm_1\otimes \cdots \otimes \zm_{k+1} \mapsto \langle X,
\zm_l\rangle  \zm_1\otimes \overset{ \overset l \to \vee} \to \cdots
\otimes \zm_{k+1},
                    \tag \label{F0.2}\endalign$$
  where $\zm_j \in \otimes ^1(\zp)$ and $ \overset l \to \vee$ stands for
the omission. We denote by $\ix_X$ the derivation of the tensor algebra
$\otimes (\zp)$:
    $$ \ix_X(\zm_1\otimes \cdots \otimes \zm_{k+1})= \sum_{l=1} ^{k+1}
(-1)^{l+1}\ix_X^l (\zm_1\otimes \cdots \otimes \zm_{k+1}).
                              $$
  Similarly, for a section  $K$ of $E^\*\otimes_M E$, we have the
derivation $\ix_K$ of the full tensor algebra $\ts(\zt)$ determined by
    $$ \ix_{\zm\otimes X}(Y) = -\langle \zm, Y\rangle X, \ \
\ix_{\zm\otimes X}(\zn) = \langle \zn,X \rangle \zm, \ \ \ix_{\zm\otimes
X}(f) =0
                              $$
  for $X,Y\in \otimes ^1(\zt)$,  $\zm,\zn \in \otimes ^1(\zp)$, $f\in
\cC^\infty(M)$.

  Let $\zL\in \otimes ^2(\zt)$. We denote by $\widetilde{\zL}$ the
mapping
    $$ \widetilde{\zL}\colon E^\* \rightarrow E, \ \
\widetilde{\zL}\circ \zm = \ix^1_\zm \zL.
                              $$

  \hl1"0.2"{Local coordinates}
  Let $(x^a), \ a=1,\dots,n$, be a coordinate system in $M$. We introduce
the induced coordinate systems
    $$\align
    (x^a, {\dot x}^b) &\ \text{in} \ \sT M,\\
    (x^a, p_b) &\ \text{in} \ \sT^\* M.
                    \endalign      $$
  Let $(e_1,\dots,e_m)$  be a basis of local sections of $\zt\colon
E\rightarrow M$ and let $(e^{1}_*,\dots, e^{m}_*)$ be the dual basis of
local sections of $\zp\colon E^\*\rightarrow M$. We have the induced
coordinate systems:
    $$\align
    (x^a, y^i),\quad & y^i=\zi(e^{i}_*), \quad \text{in} \ E,\\
    (x^a, \zx_i), \quad &\zx_i = \zi(e_i),\quad \text{in} \ E^\* ,\\
    (x^a, y^{i_1\dots i_r}) & \quad \text{in}\quad \otimes^r_M E ,\\
    (x^a, \zx_{i_1\dots i_r}) & \quad \text{in}\quad \otimes^r_M E^\*,
                    \endalign      $$
  and 
    $$\align
    (x^a, y^i,{\dot x}^b, {\dot y}^j ) &  \quad \text{in} \ \sT E,\\
    (x^a, \zx_i, {\dot x}^b, {\dot \zx}_j) & \quad \text{in} \ \sT E^\* ,\\
    (x^a, y^{i_1\dots i_r}, {\dot x}^b, {\dot y}^{j_1\dots j_r}) &
\quad \text{in}\ \sT\otimes^r_M E ,\\
    (x^a, y^i, p_b, \zp_j) & \quad \text{in}\ \sT^\*E, \\
    (x^a, \zx_i, p_b, \zf^j) & \quad \text{in}\ \sT^\* E^\* ,\\
    (x^a, \zx_{i_1\dots i_r}) & \quad \text{in}\ \otimes^r_M E^\*,\\
    (x^a,y^i, \zp_{\zb_1\dots \zb_r}) & \quad \text{in} \ \otimes ^r_E
\sT^\* E,
                    \endalign      $$
  where $\zb_i \in \{1',\dots,n'\}\cup\{1,\dots,m\}$.

   We have the canonical symplectic forms:
    $$ \zw_{E^*} = \xd p_a\wedge \xd x^a +\xd \zf_i \wedge \xd \zx_i
                              $$
on $\sT^\*E^\*$ and
    $$ \zw_{E} = \xd p_a\wedge \xd x^a +\xd \zp_i \wedge \xd y^i
                              $$
  on $\sT^\*E$, and the corresponding Poisson tensors
    $$ \zL_{E^*} = \partial _{p_a}\wedge \partial _{x^a} + \partial
_{\zf^j} \wedge \partial _{\zx_j}\ \ \text{and}\ \
  \zL_{E} = \partial _{p_a}\wedge \partial _{x^a} + \partial
_{\zp_j} \wedge \partial _{y^j}.
                              $$

  There is also a canonical isomorphism (cf. [\Ref{KU}])
    $$\cR_\zt \colon \sT^\*E^\* \longrightarrow \sT^\* E
                    \tag\label{F0.3}$$
  being an anti-symplectomorphism and also an isomorphism of double
vector bundles

                $$ \cgaps{.8;.8;1.2;.8;.8} \rgaps{.8;.4;.8}
        \CD  & \sT^\*E^\*  @()\L{\sT^\*\zp} @(1,-1) @() \L{\ezp_{E^*}}
@(-1,-2) @() \L{\cR_\zt} @(3,0) & & & \sT^\* E @()\dL{3} \L{\sT^\*\zt} @(-1,-2)
@() \L{\ezp_E} @(1,-1) & \\
        & & E @()\dl{-2}\l{ \zt} @(-1,-2) @()\dL{-11}
\L{\idt} @(3,0) & & & E @()\L{ \zt} @(-1,-2) \\
        E^\*  @()\L{\zp} @(1,-1) @()\dL{9}\L{\text{id} } @(3,0) & & & E^\*
@()\L{ \zp} @(1,-1)  & &  \\
        & M @()\L{\text{id}} @(3,0) & & & M &
                             \endCD \ \ \ .         \tag \label{F0.4} $$

  In local coordinates, $\cR_\zt$ is given by
    $$ (x^a,y^i, p_b, \zp_j)\circ \cR_\zt = (x^a, \zf^i, -p_b,\zx_j).
                              $$


        \hl1{Leibniz structures and algebroids\label{S1}}
                \claim \c{d}{Definition}{}{\rm           \label{C1.1}
        A {\it Leibniz structure} is a pair $(M,\zL)$, where $M$ is a manifold
and $\zL$ is a contravariant 2-tensor field. A Leibniz structure defines
the {\it Leibniz bracket} $\{\cdot ,\cdot \}_\zL$ on $\cC^\infty(M)$ by
                $$ \{f,g\}_\zL = \langle \zL, \xd f\otimes \xd g \rangle .
                                                                        $$
                }\endclaim
        The bracket $\{\, ,\,\}$ is a bilinear operation satisfying the Leibniz
rules
                $$\align
        \{f,gh\}_\zL&=\{f,g\}_\zL h + g\{f,h\}_\zL  \\
        \{fh,g\}_\zL&= \{f,g\}_\zL h + f\{h,g\}_\zL.
                                                \tag \label{F1.1}\endalign$$
        A Leibniz structure is called {\it skew-symmetric} if the tensor $\zL$
or, equivalently, the bracket $\{\,,\,\}_\zL$, is skew-symmetric.

        {\bf Remark.} A Leibniz bracket, which is skew-symmetric and satisfies
the Jacobi identity is called a {\it Poisson bracket} and the
corresponding tensor -- a {\it Poisson structure}. It is well known that
Lie algebroid structures on a vector bundle $E$ correspond to linear
Poisson structures on $E^\*$.

        A Leibniz structure $\zL$ on $E^\*$ is called {\it linear} if the
corresponding mapping $\widetilde{\zL} \colon \sT^\* E^\* \rightarrow \sT
E^\*$ is a morphism of double vector bundles.

        Let $\zt \colon  E\rightarrow M$ be a vector bundle. The commutative
diagram
                $$ \CD
                \sT^\* E^\* @() \L{\cR_\zt} @(0,-1) @() \L{\widetilde{\zL}}
@(1,0) & \sT E^\* \\
                \sT^\* E @() \L{\ze} @(1,1)
                        \endCD,                            \tag \label{F1.2}$$
          describes a one-to-one correspondence between linear Leibniz
structures $\zL$ on $E^\*$ and homomorphisms of double vector bundles (cf.
[\Ref{KU}])
                $$ \cgaps{.8;.8;1.2;.8;.8} \rgaps{.8;.4;.8}
        \CD  & \sT^\*E @()\L{\sT^\*\zt} @(-1,-2) @() \L{\ezp_E} @(1,-1)
@() \L{\ze} @(3,0) & & & \sT E^\* @()\dL{3} \L{\ezt_{E^*}} @(-1,-2) @()
\L{\sT\zp} @(1,-1) & \\
        & & E @()\dl{-2}\l{ \zt} @(-1,-2) @()\dL{-9}
\L{\ze_r} @(3,0) & & & \sT M @()\L{ \ezt_M} @(-1,-2) \\
        E^\*  @()\L{\zp} @(1,-1) @()\dL{7}\L{\text{id} } @(3,0) & & & E^\*
@()\L{ \zp} @(1,-1)  & &  \\
        & M @()\L{\text{id}} @(3,0) & & & M &
                              \endCD \ \ \ .         \tag \label{F1.3}$$
        The core of a double vector bundle is the intersection of the kernels
of the projections. It is obvious that the core of $\sT^\* E$ ($\sT E^\*$) can
be identified with $\sT^\* M$ ($E^\*$).
        With these identifications the induced by $\ze$ morphism of cores is a
morphism
                $$              \ze_c \colon \sT^\* M \rightarrow E^\*.
                                                                           $$


        In local coordinates, every  $\ze$ as in \Ref{F1.3} is of the form
                $$ (x^a, \zx_i, {\dot x}^b, {\dot \zx}_j) \circ \ze = (x^a,
\zp_i, \zr^b_k(x)y^k, c^k_{ij}(x) y^i\zp_k + \zs^a_j(x) p_a)
                                                            \tag \label{F1.4}$$
and corresponds to the linear Leibniz tensor
                $$ \zL_\ze = c^k_{ij}(x)\zx_k \partial _{\zx_i}\otimes \partial
_{\zx_j} + \zr^b_i(x) \partial _{\zx_i} \otimes \partial _{x^b} -
\zs^a_j(x) \partial _{x^a} \otimes \partial _{\zx_j}.
                                                            \tag \label{F1.5}$$
        We have also
                $$      \align
                (x^a, \dot x^b)\circ \ze_r &= (x^a, \zr^b_k(x)y^k), \\
                (x^a,  \zx_i)\circ \ze_c &= (x^a, \zs^b_i(x)p_b).
                                                                \endalign    $$
        {\bf Remark.} In [\Ref{GU3}] by {\it pseudo--Lie algebroids} (resp. {\it
pre--Lie algebroids}), we called the linear (resp. linear and
skew-symmetric) Leibniz tensors on $E^\*$. Throughout this note the
corresponding algebroid structures (represented by pairs $(E,\ze)$ or,
equivalently, by $(E^\*,\zL_\ze)$) will be called  simply {\it algebroids}
({\it skew algebroids}).  Note also that the notion of a skew algebroid
was introduced by Kosmann-Schwarzbach and Magri in [\Ref{KS-M}] under the
name of pre--Lie algebroid. The relation to the canonical definition of Lie
algebroid is given by the following theorem.

                \claim \c{t}{Theorem}{[\Ref{GU3}]}                 \label{C1.3}
        An algebroid structure $(E,\ze)$ can be equivalently defined as a
bilinear bracket $[\cdot ,\cdot]_\ze $ on sections of $\zt\colon
E\rightarrow M$ together with vector bundle morphisms $a^\ze_l,\, a^\ze_r
\colon E\rightarrow \sT M$ (left and right anchors) such that
                $$ [fX,gY]_\ze = f(a^\ze_l\circ X)(g)Y -g(a^\ze_r \circ Y)(f) X
+fg [X,Y]_\ze
                                                             \tag \label{F1.6}$$
         for $f,g \in \cC^\infty (M)$, $X,Y\in \otimes ^1(\zt)$.

        The bracket and anchors are related to the Leibniz tensor $\zL_\ze$ by
the formulae            $$\align
        \zi([X,Y]_\ze)&= \{\zi(X), \zi(Y)\}_{\zL_\ze},  \\
        \zp^\*(a^\ze_l\circ X(f))       &= \{\zi(X), \zp^\*f\}_{\zL_\ze}, \\
        \zp^\*(a^\ze_r\circ X(f))       &= \{\zp^\* f, \zi(X)\}_{\zL_\ze}.
                                                   \tag \label{F1.7}\endalign$$
        We have also $a^\ze_l=\ze_r$ and $a^\ze_r = (\ze_c)^\*$.
        The algebroid $(E,\ze)$ is a Lie algebroid if and only if the tensor
$\zL_\ze$ is a Poisson tensor.
        \endclaim

        The canonical example of a mapping $\ze$ in the case of $E=\sT M$ is
given by $\ze = \eze_M = \eza^{-1}_M$ -- the inverse to the Tulczyjew
isomorphism, which can be defined as the dual to the isomorphism of double
vector bundles
                $$ \cgaps{.8;.8;1.2;.8;.8} \rgaps{.8;.4;.8}
        \CD  & \sT \sT M @()\L{\sT\ezt_M} @(-1,-2) @() \L{\ezt_{\ssT M}}
@(1,-1) @() \L{\ezk_M} @(3,0) & & & \sT \sT M @()\dL{3} \L{\ezt_{\ssT M }}
@(-1,-2) @() \L{\sT\ezt_M} @(1,-1) & \\
        & & \sT M @()\dl{-2}\l{ \ezt_M} @(-1,-2) @()\dL{-10}
\L{\text{id}} @(3,0) & & & \sT M @()\L{ \ezt_M} @(-1,-2) \\
        \sT M  @()\L{\ezt_M} @(1,-1) @()\dL{7}\L{\text{id} } @(3,0) & & &
\sT M @()\L{ \ezt_M} @(1,-1)  & &  \\
        & M @()\L{\text{id}} @(3,0) & & & M &
                               \endCD \ \ \ .         \tag \label{F1.8}$$
        In general, the Leibniz structure map $\ze$ is not an isomorphism and,
consequently, its dual $\zk^{-1} = \ze^{\s*_r}$  with respect to the right
projection
                $$ \cgaps{.8;.8;1.2;.8;.8} \rgaps{.8;.4;.8}
        \CD  & \sT E @()\L{\sT\zt} @(-1,-2) @() \L{\ezt_E} @(1,-1) @()
\L{\zk} @(3,0) & & & \sT E @()\dL{3} \L{\ezt_E} @(-1,-2) @() \L{\sT\zt} @(1,-1)
& \\
        & & E @()\dl{-2}\l{ \zt} @(-1,-2) @()\dL{-7}
\L{\zk_r} @(3,0) & & & \sT M @()\L{ \ezt_M} @(-1,-2) \\
        \sT M @()\L{\ezt_M} @(1,-1) @()\dL{7}\L{\zk_l } @(3,0) & & & E
@()\L{ \zt} @(1,-1)  & &  \\
        & M @()\L{\text{id}} @(3,0) & & & M &
                           \endCD \ \ \ .         \tag \label{F1.9}$$
is a relation and not a mapping. In this diagram $\zk_r = \ze_r$.

        Similarly as in [\Ref{GU1}], we can extend $\ze$ to mappings
                $$ \ze^{\otimes r} \colon \otimes ^r_E \sT^\* E \longrightarrow
\sT \otimes ^r_M E^\*,
                                                                             $$
        $r\geqs 0$, as follows.
        We put, for $r=0$,
                $$\align
        \ze^{\otimes 0}&\colon E\times \R \rightarrow  \sT(M\times \R) = \sT
M\times \R\times \sT_0\R = \sT M\times \R \times \R\\
        &\colon (e,t) \mapsto (\ze_r(e), t, 0)
                                                  \tag \label{F1.10}\endalign$$
        and, for $r=1$, $\ze^{\otimes 1} = \ze$. For $r>1$, we apply the tangent
functor  to the tensor product map
                $$ \otimes ^r \colon E^\* \times _M \cdots \times _M E^\*
\longrightarrow \otimes ^r_ME^\*
                                                                            $$
        and we define $\ze^{\otimes r}$ by the commutativity of the diagram
                $$ \CD
        \sT E^\* \times _{\ssT M} \dots \times _{\ssT M} \sT E^\* @()
\L{\sT\otimes ^r} @(1,0) & \sT \otimes ^r_M E^\* \\
        \sT^\* E\times _E \cdots \times _E \sT^\* E @() \L{\times ^r\ze} @(0,1)
@() \L{\otimes ^r} @(1,0) & \otimes ^r_E \sT^\* E @() \l{\ze^{\otimes r}}
@(0,1)
                                                        \endCD  \ .
                                                           \tag \label{F1.11}$$
        In local coordinates, we have
                $$\multline
        (x^a, \zx_{i_1\dots i_r}, {\dot x}^b, \dot \zx_{j_1\dots j_r}) \circ
\ze^{\otimes r} \\ = \left(x^a, \zp_{i_1\dots i_r}, \zr^b_k y^k , \sum_l
\left( c^k_{ij_l}y^i \zp_{j_1\dots j_{l-1}kj_{l+1}\dots j_r} + \zs^{a'}_{j_l}
\zp_{j_1\dots j_{l-1}a'j_{l+1}\dots j_r} \right)\right).
        \endmultline
                                                           \tag \label{F1.12}$$
        Of course, $\ze^{\otimes r}$ may be reduced to skew-symmetric
(symmetric) tensors and, as in [\Ref{GU1}], we have also the commutative diagram
                $$\cgaps{1.9}\CD
        (\otimes ^i_E \sT^\*E)\times _E (\otimes ^{r-i}_E \sT^\*E) @()
\L{\ze^{\otimes i}\times \ze^{\otimes (r-i)}} @(1,0) @() \L{\otimes _E}
@(0,-1)
        & (\sT\otimes ^i_M E^\*) \times _{\ssT M} (\sT \otimes ^{r-i}_M E^\*)
        @() \L{\sT \otimes _M} @(0,-1)  \\
        \otimes ^r_E \sT^\* E @() \L{\ze^{\otimes r}} @(1,0) & \sT \otimes ^r_M
E^\*
                                          \endCD \ .      \tag \label{F1.13}$$

        \hl1{Leibniz relations}
        In the following definition we adapt standard concepts of the theory
of Poisson manifolds and Lie algebroids.
                \claim \c{d}{Definition}{}{\rm A submanifold $N$ of a Leibniz
manifold $(M,\zL)$ is called                                     \label{C1.4}
        {\it coisotropic } if $\widetilde{\zL} ((\sT N)^\circ ) \subset \sT N$,
where $(\sT N)^\circ \subset \sT^\*_NM$ is the annihilator of $\sT
N\subset \sT M$. A relation $L\subset M_1\times M_2$ between  Leibniz
manifolds $(M_1,\zL_1)$ and $(M_2,\zL_2)$ is a {\it Leibniz relation} if
$L$ is a coisotropic  submanifold of Leibniz manifold $(M_1\times M_2,
(-\zL_1) \times \zL_2)$.
                }\endclaim

                \claim \c{t}{Theorem}{}                          \label{C1.5}
        Let $F\colon M_1\rightarrow M_2$ be a differentiable mapping between
Leibniz manifolds $(M_1,\zL_1)$ and $(M_2,\zL_2)$. Then, $\zL_1$ and
$\zL_2$ are $F$-related if and only if the graph of $F$ is a Leibniz
relation.
                \endclaim
                \proof
        The proof is completely parallel to that in the Poisson case.
                \endproof

                \claim \c{d}{Definition}{}{\rm                   \label{C1.6}
        Let $\zt_i\colon E_i \rightarrow  M_i$ be a vector bundle, $i=1,2$.
Let $\ze_i$ be an algebroid structure on $\zt_i$ and let $\zL_i$ be the
corresponding Leibniz structure on $E^\*_i$. A vector bundle morphism
$\zc \colon E_1\rightarrow E_2$ is      an {\it algebroid morphism } if the
dual relation $\zc^\* \colon E_2^\* \rightarrow E_1^\*$ is a Leibniz
relation.
                }\endclaim

        \hl1{The algebroid lift $\dt^\ze$}
        For a tensor field $K\in \otimes ^k(\zt)$, we can define the vertical
lift $\text{\tenrm v}_\zt(K) \in \otimes ^k(\ezt_E)$ (cf.
[\Ref{GU2}],[\Ref{YI}]). In local coordinates,
                $$ \text{\tenrm v}_\zt(f^{i_1\dotsi_k} e_{i_1}\otimes \cdots
\otimes e_{i_k}) = f^{i_1\dotsi_k} \partial _{y^{i_1}} \otimes \cdots
\otimes \partial _{y^{i_k}}.
                                                      \tag \label{F1.14}$$
        A particular case of the vertical lift is the lift $\vt(K)$ of a
contravariant tensor field $K$ on $M$ into a contravariant tensor field on
$\sT M$.
                \claim \c{l}{Lemma}{}                         \label{C1.7}
        Let $(E,\ze)$ be an algebroid and let $K\in \otimes ^k(\zt), \ k\geqs
0$. Then
                $$  \zi(\vta(K) ) = \vt(\zi(K))\circ \ze^{\otimes k}.
                                                            \tag \label{F1.15}$$
                \endclaim
                \proof
        In local coordinates, both sides of \Ref{F1.15} are equal
                $$ f^{i_1\dots i_k}(x)\zp_{i_1\dots i_k}
                                                                            $$
for $k>0$. For a function $f$ on $M$, we get
                $$ \zi(\vta f)(e,t)= t f(\zt(e)) = tf(\zt_M(\ze_r(e))) =
\vt(\zi(f)) \circ \ze^0(e,t).
                       $$
                \endproof
        It is well  known (see [\Ref{YI}], [\Ref{GU1}]) that in the case of
$E=\sT M$ we have also the tangent lift $\dt \colon \otimes (\ezt_M)
\rightarrow \otimes (\ezt_{\ssT M})$ which is a $\vt$-derivation. It
turns out that the presence of such a lift for a vector bundle is
equivalent to the presence of an algebroid structure.

                \claim \c{t}{Theorem}{}                      \label{C1.8}
Let $(E,\ze)$ be an algebroid. For $K\in \otimes ^k(\zt)$, $k\geqs 0$, the
equality
                $$ \zi(\dt^\ze(K)) = \dt (\zi(K))\circ \ze^{\otimes k}
                                                            \tag \label{F1.17}$$
         defines the  tensor field $\dt^\ze(K) \in \otimes ^k(\ezt_E)$  which
is linear and the mapping
                $$ \dt^\ze \colon \otimes (\zt) \longrightarrow \otimes (\ezt_E)
                                                                             $$
is a $\vta$-derivation of degree $0$. In local coordinates,
                $$ \dt^\ze(f^i(x)e_i) = f_i(x)\zs^a_i(x)\partial _{x^a} +
\left( y^i \zr^a_i(x) \frac{\partial f^k}{\partial x^a}(x) +
c^k_{ij}(x)y^if^j(x) \right) \partial _{y^k}.
                                                           \tag \label{F1.18}$$
        Conversely, if $D\colon  \otimes (\zt)\rightarrow \otimes (\ezt_E)$ is
a $\vta$-derivation of degree $0$ such that $D(K)$ is linear for each $K
\in \otimes ^1(\zt)$, then there is an algebroid structure $\ze$ on
$\zt\colon E\rightarrow M$ such that $D = \dt^\ze$.
                \endclaim
                \proof
        Since $\dt (\zi(K))\circ \ze^{\otimes k}$ is a linear function on
$\otimes ^k_E\sT^{\textstyle *} E$, it defines a unique tensor field
$\dt^\ze(K) \in \otimes ^k(\ezt_E)$. In local coordinates,
                $$\multline
        \dt(\zi(f^{i_1\dots i_k}e_{i_1}\otimes \dots \otimes e_{i_k}))
        = \dt(f^{i_1\dots i_k}\zx_{i_1\dots i_k} ) = \frac{\partial
f^{i_1\dots i_k}}{\partial x^a}{\dot x}^a \zx_{i_1\dots i_k} + f^{i_1\dots
i_k} \dot \zx_{i_1\dots i_k} .
                                \endmultline              \tag \label{F1.19}$$
         Hence,
                $$\multline
        \dt(\zi(f^{i_1\dots i_k}e_{i_1}\otimes \dots \otimes e_{i_k}))\circ
\ze^{\otimes k} =\\
        =  \frac{\partial f^{i_1\dots i_k}}{\partial x^a}\zr^a_iy^i
\zp_{i_1\dots i_k} + f^{i_1\dots i_k} \sum_l \left( c^r_{ij_l}y^i
\zp_{j_1\dots j_{l-1}rj_{l+1}\dots j_r} + \zs^{a'}_{j_l} \zp_{j_1\dots
j_{l-1}a'j_{l+1}\dots j_r} \right)  .
                                \endmultline             \tag \label{F1.20}$$
        and
                $$\multline
        \dt^\ze(f^{i_1\dots i_k}e_{i_1}\otimes \dots \otimes e_{i_k})=\\
=\frac{\partial f^{i_1\dots i_k}}{\partial x^a}\zr^a_iy^i \partial
_{y^{i_1}} \otimes \cdots \otimes \partial _{y^{i_k}} + f^{i_1\dots i_k}
 \sum_l \partial _{y^{i_1}}\otimes \cdots\otimes \left( c^j_{ii_l}y^i
\partial _{y^i} + \zs^a_{i_l}\partial _{x^a}\right) \otimes \cdots\otimes
\partial_{y^{i_k}}.
                                \endmultline            \tag \label{F1.21}$$
        Now, it is easy to see that $\dt^\ze(K)$ is linear and that $\dt^\ze$
is a $\vta$-derivation of order $0$.

        We can give another, more intrinsic, proof of the fact that $\dt^\ze$
is a $\vta$-derivation.

Let $K\in \otimes ^k(\zt)$, $L \in\otimes ^l(\zt)$ and let us consider
$\widetilde{K} = \zi(K)\circ \otimes ^k$, $\widetilde{L} = \zi(L)\circ
\otimes ^l$, and $\widetilde{K\otimes L} = \widetilde{K}\cdot
\widetilde{L}$ as functions on $(\times ^{k+l}_M E^{\textstyle *}) \times
_M (\times _M^l)$. Since
                $$ \dt(\widetilde{K}\widetilde{L}) =
\dt\widetilde{K}\vt\widetilde{L} + \vt\widetilde{K}\dt\widetilde{L}
                                                                             $$
(cf. [\Ref{GU1}]), we get
                $$\dt(\widetilde{K}\widetilde{L}) \circ (\times^{k+l}\ze ) =
\left(\dt\widetilde{K} \circ (\times ^k\ze)\right)\left(\vt\widetilde{L}
\circ (\times ^l\ze)\right) + \left(\vt\widetilde{K}\circ (\times ^k\ze)
\right) \left( \dt\widetilde{L}\circ (\times ^l\ze) \right)
                                                                             $$
        and, in view of Lemma \Ref{C1.7},
                $$ \dt^\ze (K\otimes L) = \dt^\ze(K) \otimes \vta (L) + \vta(K)
\otimes \dt^\ze (L)  .
                                                                              $$
        To prove the converse, we use the method similar to the method used in
the proof that derivations of $\cC^\infty(M)$ are given by vector fields.

        Let $D\colon \otimes \zt\rightarrow \otimes (\ezt_E)$ be a
$\vta$-derivation of degree zero. It follows that $D\colon \cC^\infty(M)
\rightarrow \cC^\infty(E)$  is also a $\vta$-derivation.
Consequently,
                $$ D(1) = D(1\cdot 1)=D(1)\vta(1) +\vta(1) D(1) = 2D(1)
                                                            \tag \label{F1.22}$$
        and $D(1)=0$.

        It is well known that for every $m_0\in M$ we can find local coordinates
$(x^a)$ in a neighbourhood $U$ of $m_0$, $x^a(m_0)=0$, given by globally
defined functions, such that for each $f\in \cC^\infty(M)$ we have
                $$ f(m) = f(m_0) + (x^af_a)(m)
                                                                              $$
for $m\in U$ and some $f_a\in \cC^\infty(M)$. It is clear that $f_a(m_0) =
\dfrac{\partial f}{\partial x^a}(m_0) $.

        Let $e\in E$, $\zt(e) = m_0$. We have
                $$\align
        (Df)(e)&= D(f(m_0))(e) +D(x^a)(e)\vta(f_a)(e) +\vta(x^a)(e) D(f_a)(e) \\
                        &= f_a(m_0)D(x^a)(e) = \dfrac{\partial f}{\partial
x^a}(m_0) D(x^a) (e) .
                                                  \tag \label{F1.23}\endalign$$
        We can apply this arguments to every point $m\in U $ and coordinates
$x^a - x^a(m)$ and, since $D(x^a - x^a(m)) = D(x^a)$, we obtain the local
formula
                $$ D(f) = \dfrac{\partial f}{\partial x^a} D(x^a).
                                                            \tag \label{F1.24}$$

        Let $D(K)$ be linear for $K\in \otimes ^k(\zt)$. Then  the functions
$D(x^a)$ are linear on fibers of $E$ and, consequently, we can write
$D(x^a) = \zr^a_j y^j$ for some functions $\zr^a_j\in \cC^\infty(U)$, and
                $$ D(f) = \dfrac{\partial f}{\partial x^a} \zr^a_jy^j .
                                                                              $$
        Let $(e_i)$ be a basis of local sections of $E$. We have, due to
linearity of $D(e_i)$,
                $$ D(e_i) = \zs^a_i\partial _{x^a} + c^k_{ji}y^j \partial _{y^k},
                                                            \tag \label{F1.25}$$
        where $\zs^a_i, c^k_{ji} \in \cC^\infty(U)$.
        Finally, the equality
                $$\align
        D(f^ie_i)&= D(f^i)\vta(e_i) + \vta(f^i)D(e_i)  \\
                        &= \dfrac{\partial f^i}{\partial x^a}\zr^a_jy^j \partial
_{y^i} + f^i c^k_{ji}y^j\partial _{y^k} +f^i\zs^a_i\partial _{x^a}
                                                \tag \label{F1.26}\endalign$$
        shows that, when restricted to functions and sections of $E$, the
derivation  $D$ equals $\dt^\ze$ for $\ze$ as in \Ref{F1.4} when acting on
functions and sections of $E$. But the $\vta$-derivation $D$ of order $0$
is uniquely determined by its values on functions and on sections of $E$.
                \endproof

                \claim \c{t}{Theorem}{}                           \label{C1.9}
        For $X,Y \in \otimes ^1(\zt)$,
                $$ [\vta(X), \dt^\ze(Y)] = \vta([X,Y]_\ze).
                                                           \tag \label{F1.27}$$
                \endclaim
                \proof
        Easy calculations in local coordinates.
                \endproof
                \claim \c{c}{Corollary}{}                       \label{C1.10}
        For $X\in \otimes ^1(\zt)$ and $Y\in \otimes ^k(\zt)$, we have
                $$ [\vta(X), \dt^\ze(Y)] = \vta([X,Y]_\ze)
                                                                             $$
        and
                $$ [\vta(Y), \dt^\ze(X)] = \vta([Y,X]_\ze),
                                                                              $$
        where
                $$ [X,Y_1\otimes \dots\otimes Y_k] = \sum_iY_1\otimes
\dots\otimes [X,Y_i] \otimes \dots \otimes Y_k,
                                                                              $$
                $$ [Y_1\otimes \dots\otimes Y_k,X] = \sum_iY_1\otimes
\dots\otimes [Y_i,X] \otimes \dots \otimes Y_k,
                                                                             $$
        and the formulae for $[\,,\,]_\ze$ are similar.
                \endclaim
                \proof
        The proof follows directly from  Theorem \Ref{C1.9} and from the fact
that $\dt^\ze$ is a $\vta$-derivation of order $0$.
                \endproof

        The algebroid lift $\dt^\ze$ may be used to an alternative definition of
the bracket $[\,,\,]_\ze$ (cf. [\Ref{GU1}] for the case of the canonical
Lie algebroid in the tangent bundle).

        Let $X$ be a section of $E$. At points of $X(M)\subset E$, we have the
decomposition of $\sT E$ into the horizontal (tangent to $X(M)$) and the
vertical part. For a vector field $V$ on $E$, we denote by $X^{\textstyle
*} V$ the unique section of $E$ such that the vertical lift
$\vta(X^{\textstyle *} V)$ is, on $X(M)$, the vertical part of $V$.

                \claim \c{t}{Theorem}{}                          \label{C1.11}
        Let $(E,\ze)$ be an algebroid and let $X,Y\in \otimes ^1(\zt)$.
Then
                $$ X^{\textstyle *} (\dt^\ze Y) = [X,Y]_\ze.
                                                          \tag \label{F1.28}$$
                \endclaim
                \proof
        We have to prove that the vector field  $\dt^\ze(Y) - \vta([X,Y]_\ze)$
is tangent to $X(M) \subset E$. Let $X =f^ie_i$ and $Y=g^ie_i$, in local
coordinates, and let $\ze$ be as in \Ref{F1.4}. We have, according to
\Ref{F1.6}, \Ref{F1.14}, and \Ref{F1.18},
                $$\multline
                \dt^\ze (g^je_j) - \vta([f^ie_i,g^je_j]_\ze) =\\
= c^k_{ij}y^ig^j \partial _{y_k} + \zr^a_i y^i \frac{\partial
g^k}{\partial x^a} \partial _{y^k} + \zs^a_ig^i \partial _{x^a} -
c^k_{ij}f^ig^j \partial _{y_k} - \zr^a_i f^i \frac{\partial
g^k}{\partial x^a} \partial _{y^k} + \zs^a_ig^i \frac{\partial
f^k}{\partial x^a}  \partial _{y^k}\\
        = \zs^a_ig^i \left(\partial _{x^a} + \frac{\partial f^k}{\partial x^a}
\partial _{y^k} \right),
                                \endmultline              \tag \label{F1.29}$$
        since $y^i = f_i$ on $X(M)$. It is clear that the vector fields
$\partial _{x^a} + \dfrac{\partial f^k}{\partial x^a} \partial _{y^k}$ are
tangent to $X(M)$. \phantom{aaaaaa}
                \endproof

        Let $\zf\colon \sT^\*E\rightarrow \sT E^\*$ be a morphism of vector
bundles over $E^\*$. The vector bundle dual to $\ezt_{E^*}\colon \sT E^\*
\rightarrow E^\*$ can be identified (via $\cR_\zt$) with  the vector bundle
$\sT^\*\zt\colon \sT^\* E\rightarrow E^\*$ and vice versa.
        It follows that the dual to $\zf$ can be identified with a morphism
                $$ \zf^+ \colon  \sT^\* E \rightarrow  \sT E^\* .
                                                                       \tag $$

                \claim \c{t}{Theorem}{}                          \label{C1.12}
        Let $\ze$ be an algebroid structure as in \Ref{F1.4}. Then the morphism
                $$  \ze^+ \colon \sT^{\textstyle *} E \longrightarrow
\sT E^{\textstyle *} ,
                                                                     $$
        dual to $\ze$ with respect to the vector bundle structures over
$E^{\textstyle *} $ is again an algebroid structure on $E$. The
corresponding Leibniz tensor $\zL_{\ze^+}$ is the transposition of the
tensor $\zL_\ze$. In local coordinates,
                $$ (x^a, \zx_i, {\dot x}^b, {\dot \zx}_j) \circ \ze^+
= (x^a, \zp_i, \zs^b_k(x)y^k, c^k_{ji}(x) y^i\zp_k + \zr^a_j(x) p_a).
                                                      \tag \label{F1.30}$$

                \endclaim
                \proof
        The proof is an immediate consequence of the general properties of a
dual to a double vector bundle. We refer to ([\Ref{KU}], [\Ref{Ur}]) for
the theory of duality in the category of double vector bundles.
                \endproof
        It follows from the general theory of double vector bundles that
$(\ze^+)_r = (\ze_c)^\*$ and $(\ze^+)_c = (\ze_r)^\*$. Thus, we have the
diagram
                $$ \cgaps{.8;.8;1.2;.8;.8} \rgaps{.8;.4;.8}
    \CD  & \sT^\*E @()\L{\sT^\*\zt} @(-1,-2) @() \L{\ezp_E} @(1,-1) @()
\L{\ze^+} @(3,0) & & & \sT E^\* @()\dL{3} \L{\ezt_{E^*}} @(-1,-2) @() \L{\sT\zp}
@(1,-1) & \\
        & & E @()\dl{-2}\l{ \zt} @(-1,-2) @()\dL{-11}
\L{(\ze_c)^\*} @(3,0) & & & \sT M @()\L{ \ezt_M} @(-1,-2) \\
        E^\*  @()\L{\zp} @(1,-1) @()\dL{7}\L{\text{id} } @(3,0) & & & E^\*
@()\L{ \zp} @(1,-1)  & &  \\
        & M @()\L{\text{id}} @(3,0) & & & M &
                              \endCD \ \ \ .         \tag \label{F1.31}$$

We call an algebroid structure $(E,\ze)$ {\it skew--symmetric} if $\ze =
-\ze^+$, where the multiplication by $-1$ is with respect to the vector
bundle structure over $E^\*$. It is clearly equivalent to the fact that
$\zL_\ze$ is skew-symmetric. Skew-symmetric algebroids are sometimes
called {\it pre-Lie algebroids} ([\Ref{GU3}], [\Ref{KS-M}]).

        \hl1{The tangent and cotangent lifts of an algebroid}
        The theory of the tangent and cotangent lifts of  algebroids is
perfectly analogous to that of Lie algebroids (see [\Ref{Co}],
[\Ref{GU2}]). If $\ze$ is an algebroid structure on $\zt\colon
E\rightarrow M$ and $\zL_\ze$ is the corresponding Leibniz tensor on
$E^{\textstyle *} $, the tangent lift $\dt (\zL_\ze)$ is linear with
respect to both vector bundle structures on $\sT E^{\textstyle *} $ and,
consequently, it defines algebroid structures on dual bundles:
$\sT^{\textstyle *} E^{\textstyle *} \rightarrow E^{\textstyle *} $  and
$\sT E\rightarrow  \sT M$, called respectively the {\it cotangent} and {\it
tangent} lift of the algebroid $\ze$. Let us denote the corresponding
double vector bundle morphisms by $\cT^{\textstyle *} \ze$ and $\cT \ze$
respectively. It is easy to prove the following theorem.

                \claim \c{t}{Theorem}{}                           \label{C1.13}
        The tangent $\cT\ze$ and the cotangent $\cT^{\textstyle *} \ze$ lifts of
an algebroid structure $\ze$ on $\zt\colon E\rightarrow M$ are defined by
the following commutative diagram

                $$  \cgaps{1;1.7} \CD
                \sT\sT^\* E  @()\L{\sT\ze} @(0,-1) @()\L{\eza_E} @(1,0) & \sT^\*
\sT E @()\L{\cT\ze} @(0,-1)  @()\L{ \cR_{\ezp_{E^*}} \circ \cR_{\ssT\zp} }
@(1,0) & \sT^\* \sT ^\* E^\* @()\L{\cT^\*\ze} @(0,-1)\\
        \sT\sT E^\*  @()\L{\ezk_{E^*}} @(1,0)    & \sT\sT E^\*  @()\L{\idt}
@(1,0) &  \sT\sT E^\*
                             \endCD \ .                   \tag \label{F1.32}$$


                \endclaim
                \proof
        We have two commutative diagrams
                $$ \CD
                \sT^\* E^\* @() \L{\cR_\zt} @(0,-1) @() \L{\widetilde{\zL_\ze}}
@(1,0) & \sT E^\* \\
                \sT^\* E @() \L{\ze} @(1,1)
                        \endCD  \quad \text{and} \qquad
                \CD
                \sT^\*\sT E^\*  @() \L{\cR_{\ssT \zt}} @(0,-1) @()
\L{\eza_{E^*}} @(1,0) & \sT \sT^\* E^\*  @() \L{\sT \cR_\zt} @(0,-1)\\
                \sT^\* \sT E  @() \L{\eza_E} @(1,0)  & \sT \sT^\* E
                        \endCD  \ .                         \tag \label{F1.33}$$

It follows that the diagram
                $$ \cgaps{1.7;1;1}
                \CD
        \sT^\* \sT E^\* @() \L{\idt} @(1,0) @() \L{\cR_{\ezp_{E^*}}} @(0,-1) &
\sT^\*\sT E^\*  @() \L{\cR_{\ssT \zt}} @(0,-1) @() \L{\eza_{E^*}} @(1,0) &
\sT \sT^\* E^\* @() \L{\sT \cR_\zt} @(0,-1) @() \L{\sT\widetilde{\zL_\ze}}
@(1,0) & \sT \sT E^\* \\
        \sT^\*\sT^\* E^\* @() \L{\cR_{\ssT \zt}\circ \cR_{\ezt_{E^*}}} @(1,0)
&       \sT^\* \sT E  @() \L{\eza_E} @(1,0)  & \sT \sT^\* E @() \L{\sT\ze}
@(1,1) &
                        \endCD
                                                            \tag\label{F1.34}$$
         is also commutative.
        All the mappings in this diagram respect three vector bundle structures
and, since $\widetilde{\dt\zL_\ze} = \ezk_{E^*} \circ
\sT\widetilde{\zL_\ze} \circ \eza_{E^*}$ (cf. [\Ref{GU1}]),
        we have  $\cT\ze  $, $\cT^\*\ze $ as in the diagram \Ref{F1.32}, and
the corresponding double vector bundle morphisms:
                $$ \cgaps{.8;.8;1.2;.8;.8} \rgaps{.8;.4;.8}
        \CD  & \sT^\* \sT E @()\L{\sT^\*\sT\zt} @(-1,-2) @() \L{\ezp_{\ssT
E}} @(1,-1) @() \L{\cT\ze} @(3,0) & & & \sT\sT E^\* @()\dL{3} \L{\ezt_{\ssT
E^*}} @(-1,-2) @() \L{\sT\sT\zp} @(1,-1) & \\
        & & \sT E @()\dl{-2}\l{\sT \zt} @(-1,-2) @()\dL{-17} \L{\ezk_M\circ
\sT\ze_r} @(3,0) & & & \sT\sT M @()\L{ \ezt_{\ssT M}} @(-1,-2) \\
      \sT  E^\*  @()\L{\sT\zp} @(1,-1) @()\dL{7}\L{\text{id} } @(3,0) & &
& \sT E^\* @()\L{\sT \zp} @(1,-1)  & &  \\
        &\sT M @()\L{\text{id}} @(3,0) & & &\sT M &
                              \endCD \ \ \          \tag \label{F1.35}$$
        and
                $$ \cgaps{.8;.8;1.2;.8;.8} \rgaps{.8;.4;.8}
        \CD  &\sT^\* \sT^\*E^\* @()\L{\sT^\*\zp_{E^*}} @(-1,-2) @()
\L{\ezp_{\ssT^*E}} @(1,-1) @() \L{\cT^\*\ze} @(3,0) & & & \sT\sT E^\* @()
\dL{3}\L{\ezt_{\ssT E^*}} @(-1,-2) @() \L{\sT\ezt_{E^*}} @(1,-1) & \\
    & &\sT^\* E^\* @()\dl{-2}\l{ \ezp_{E^*}} @(-1,-2) @()\dL{-11}
\L{\widetilde{\zL}_\ze} @(3,0) & & & \sT E^\* @()\L{ \ezt_{E^*}} @(-1,-2)
\\
       \sT E^\*  @()\L{\ezt_{E^*}} @(1,-1) @()\dL{7}\L{\text{id} } @(3,0)
& & & \sT E^\* @()\L{ \ezt_{E^*}} @(1,-1)  & &  \\
        & E^\* @()\L{\text{id}} @(3,0) & & & E^\* &
                            \endCD \ \ \ .         \tag \label{F1.36}$$
        We see from these diagrams, that $(\cT\ze)_r = \ezk_M\circ \sT\ze_r$ and
$(\cT^\* \ze)_r= \widetilde{\zL_\ze}$. We have also
                $$ \align
                (\cT\ze)_c &= \sT\ze_c \circ \eza_M^{-1}\\
                (\cT^\*\ze)_c &= \widetilde{\zL_\ze}
                                                 \tag \label{F1.a}\endalign$$
                \endproof

                \claim \c{t}{Theorem}{}                          \label{C1.14}
        Let $\ze$ be an algebroid structure on $\zt\colon E\rightarrow M$. The
following properties of $\ze$ are equivalent
        \list
        \item"(a)" $\ze$ is a Lie algebroid structure,
        \item"(b)" $\zL_\ze$ is a Poisson tensor,
        \item"(c)" $\zL_E$ and $\dt\zL_\ze$ are $\ze$--related,
        \item"(d)" The following diagram is commutative (on the domain of the
relation $\sT^\* \ze$)
                $$      \CD
        \sT\sT^\* E @() \L{\sT \ze} @(0,-1)& \sT^\*\sT^\* E @()
\L{\widetilde{\zL}_E} @(-1,0) & \\
        \sT\sT E^\* & \sT^\* \sT E^\* @() \L{\sT^\*\ze} @(0,1) @()
\L{\eza^{-1}_{E^*}} @(1,-1) & \\
        \sT\sT E^\* @() \L{\ezk_{E^*}} @(0,1)& \sT\sT^\* E @() \L{\sT\ze}
@(-1,0) & \sT\sT^\* E^\*  @() \L{\sT\cR_\zt} @(-1,0)
                                            \endCD  \ \ , \tag \label{F1.37}$$

        \item"(e)" $\dt^\ze([X,Y]_\ze) = [\dt^\ze(X), \dt^\ze(Y)]$ for all
$X,Y\in \otimes ^1(\zt)$.

        \endlist
                \endclaim
                \proof
                (a)$\Longleftrightarrow$(b) by the definition of a Lie
algebroid.

                (b)$\Longleftrightarrow$(c)$\Longleftrightarrow$(d) is a
version of Theorem~4.4 and Theorem~4.5 from  [\Ref{GU1}]: a Leibniz tensor
$\zL$ is a Poisson tensor on a manifold $N$ if and only if $\zL_N$ and
$-\dt\zL$ are $\widetilde{\zL}$--related.

                Finally, (c)$\Longleftrightarrow$(e), since
                $$ \zi_{\ssT^*E}([\dt^\ze(X), \dt^\ze(Y)]) = \{ \zi\dt^\ze (X),
\zi\dt^\ze (Y) \}_{\zL_E} =\{ \dt(\zi(X))\circ \ze, \dt(\zi(Y)) \circ
\ze\}_{\zL_E}
                                                                        $$
which equals
                $$\multline
        \{\dt(\zi(X)), \dt(\zi(Y))\}_{\sdt\zL_\ze} \circ \ze = \dt\left( \{
\zi(X),\zi(Y)\}_{\zL_\ze} \right) \circ \ze \\
        \dt\left(\zi([X,Y]_\ze)\right) \circ \ze = \zi\left(
\dt^\ze([X,Y]_\ze)\right)
                                \endmultline                            $$
if and only if $\zL_E$ and $\dt\zL_\ze$ are $\ze$-related.
                \endproof

        \hl1{The Lie and exterior derivatives}

        let $(E,\ze)$ be an algebroid, $X\in \otimes ^1(\zt)$, $\zm \in
\otimes ^1 (\zp)$. Since $\dt^\ze(X)$ is a linear vector field on $E$,
$\dt^\ze(\zi(\zm))$ is a linear function on $E$. The corresponding section
of $E^\*$ we call the {\it Lie derivative of $\zm$ along $X$} and denote
$\ll^\ze_X(\zm)$, i.e.,
                $$ \zi(\ll^\ze_X(\zm)) = \dt^\ze(X)(\zi(\zm)).
                                                   \tag"(A)" \label{F1.38}$$
         In local coordinates,
                $$ \ll^\ze_{f^ie_i}(\zm_je_*^j) = f^i\zs_i^b \frac{\partial
\zm_k}{\partial x^b} e_*^k + \left(f^i c_{ji}^k + \frac{\partial
f^k}{\partial x^a} \zr^a_j \right) \zm_ke_*^j .
                                                   \tag "(B)"\label{F1.39}$$
        It is easy to see that
                $$ \ll^\ze_X(f\zm) = \ll^\ze_X(\zm) + a^\ze_r(X)(f)\cdot \zm
                                                    \tag"(C)" \label{F1.40}$$
        and
                $$ \ll^\ze_{fX}(\zm) = f\ll^\ze_X(\zm) + \langle X, \zm\rangle
\xd^\ze_l(f) ,
                                                   \tag"(D)" \label{F1.41}$$
        where {\it the left derivative} $\xd^\ze_l(f) \in \otimes ^1(\zp)$ is
defined by
                $$ \langle \xd^\ze_l(f), Y\rangle  = \langle \xd f,
a^\ze_l(Y)\rangle , \ \ Y\in \otimes ^1(\zt).
                                                                            $$
        Similarly, we can define the {\it right derivative}
                $$ \langle \xd^\ze_r(f), Y\rangle  = \langle \xd f,
a^\ze_r(Y)\rangle , \ \ Y\in \otimes ^1(\zt).
                                                                           $$
        We have clearly
                $$ \xd^\ze_r(fg) = f\xd^\ze_r(g) + g \xd^\ze_r(f)
                                                      \tag"(E)" \label{F1.42}$$
        and the analog identity for the left derivative.

                \claim \c{t}{Theorem}{}                           \label{C1.15}
        For every $X \in \otimes ^1(\zt)$ there is a unique derivative
$\ll^\ze_X$ of the tensor algebra $\ts(\zt)$, called the Lie derivative
along $X$, such that
        \list
        \item"(a)" $\ll^\ze_X (f)= a^\ze_r(X)(f) = \langle X,\xd^\ze_rf
\rangle $\ for\ $f\in \cC^\infty(M)$,
        \item"(b)" $\zi(\ll^\ze_X(\zm)) = \dt^\ze(X)(\zi(\zm))$\ for\ $\zm
\in\otimes ^1(\zp)$,
        \item"(c)" $\ll^\ze_X(Y) = -[Y,X]_\ze$\ for\ $Y \in \otimes ^1(\zt) $.
                \endlist
        Moreover,
        \list
        \item"(d)" $\ll^\ze_{fX}(K) = f\ll^\ze_X(K) + \ix_{\xd^\ze_l(f)\otimes
X}(K) $\ for\ $f\in C^\infty(M)$\ and\ $K\in \ts(\zt)$.
        \endlist
                \endclaim
                \proof
        Any derivative of the tensor algebra $\ts(\zt)$ is uniquely determined
by its values on functions and sections of $E$ and $E^\*$. On the other
hand, it is easy to see that the identities
                $$\align
        \ll^\ze_X(f\cdot g)&= \ll^\ze_X(f)\cdot g + f\cdot \ll^\ze_X(g),  \\
        \ll^\ze_X(f\zm) &= \ll^\ze_X(f)\zm + f\ll^\ze_X(\zm),  \\
        \intertext{and}
                \ll^\ze_X(fY)   &= \ll^\ze_X(f)Y +f \ll^\ze_X(Y),
                                                 \tag \label{F1.43}\endalign$$
        for $f,g\in \cC^\infty(M)$, $\zm \in \otimes ^1(\zp)$, $Y\in \otimes
^1(\zt)$, which follow from \Ref{F1.40}, \Ref{F1.42}, and the properties of
$[\,,\,]_\ze$, are sufficient to extend $\ll^\ze_X$ to  a derivation of the
tensor algebra:
                $$ \ll^\ze_X(A_1\otimes \cdots\otimes A_k) = \sum_i A_1\otimes
\cdots\otimes \ll^\ze_X(A_i)\otimes \cdots \times A_k.
                                                                           $$
        It remains to prove (d). Since $\ix_{\xd^\ze_l(f)\otimes X}$ is also a
derivation, it i sufficient to check (d) on functions and sections of $E$
and $E^\*$. On functions $\ix_{\xd^\ze_l(f)\otimes X}$ acts trivially and
(d) follows from (a).  For  $\zm\i \otimes ^1(\zm)$, we have
$\ix_{\xd^\ze_l(f)\otimes X}\zm = \langle X,\zm \rangle \xd^\zl_l(f)$  and
(d) follows from \Ref{F1.41}. Finally,
                $$\ix_{\xd^\ze_l(f)\otimes X}(Y) = -a^\ze_l(Y)(f)\cdot X =
-[Y,fX]_\ze + f[Y,X]_\ze
                                                                            $$
        for $Y\in \otimes ^1(\zt)$ and (d) follows from (c).
                \endproof

        {\bf Remark.} \Ref{F1.38} was introduced in [\Ref{MX2}] as the
definition of $\dt^\ze(X)$ for $\ze$ being a Lie algebroid structure.

                \claim \c{t}{Theorem}{}                          \label{C1.16}
        Let $\otimes ^1(\zt)\ni X \rightarrow \ll_X \in \der(\otimes (\zp))$ be
a linear mapping assigning to each section of $E$ a derivative of the
tensor algebra $\otimes (\zp)$. Then $\ll_X = \ll_X^\ze$ for an algebroid
structure $\ze$ on $E$ if and only if $\ll_{fX} = f\ll_X +
\ix_{\xd_l(f)\otimes X}$ for a derivation $\xd_l \colon \cC^\infty(M)
\rightarrow \otimes ^1(\zp)$.
                \endclaim
                \proof
        The ''only if'' part follows from the previous theorem. To prove the
''if'' part, we can start with the identity
                $$ \zi(\ll_X\zm) = \xD(X)(\zi(\zm))
                                                                              $$
        to define a lift $\xD(X)$ of $X$ to a tangent vector field on $E$, and
to show that this lift is of the form $\dt^\ze$.

        We can also find directly the bracket $[\,,\,]_\ze$ as follows. First,
since
                $$ \left(\ix_X\ll_Y -\ll_Y \ix_X\right)(f\zm) = f \left(\ix_X\ll_Y
-\ll_Y \ix_X\right)(\zm),
                                                                              $$
        there is an element $[X,Y] \in \otimes ^1(\zt)$ such that $ \ix_X\ll_Y
-\ll_Y \ix_X = \zi_{[X,Y]}$. The bracket $[\,,\,]$ is bilinear and
                $$  \left(\ix_{fX}\ll_{gY} -\ll_{gY} \ix_{fX}\right) =
fg\left(\ix_X\ll_Y \ll_Y \ix_X\right) + f\ix_X\ix_{\xd_lg\otimes Y} -
g\ll_Y(f)\cdot \zi_X
                                                                              $$
which shows that $[\,,\,]$ is an algebroid bracket with the left anchor
defined by $a_l(X)(g) = \langle X,\xd_l g\rangle $ and the right anchor
defined by  $a_r(Y)(f) = \ll_Y(f) $.  The right anchoris really tensorial
with respect to $Y$, since
                $$ \ll_{gY}(f) = g\ll_Y(f) + \ix_{\xd_lg\otimes Y}f = g\ll_Y(f).
                                                                              $$
        Hence, $[\,,\,] = [\,,\,]_\ze$ for an algebroid structure $\ze$ and it
is easy to see that $\ll_X =\ll^\ze_X$.
                \endproof

        We have already defined the right and the left exterior derivatives
                $$ \xd_r^\ze,\ \xd_l^\ze \colon \cC^\infty(M) \rightarrow \otimes
^1(\zp)
                                                                               $$
        with the property
                $$ \ll^\ze_X(f) = \ix_X\xd^\ze_r f + \xd_l^\ze\ix_X f, \quad f\in
\cC^\infty(M).
                                                            \tag \label{F1.44}$$
        In general, for $\ze$ which is not skew-symmetric, it is not true that
$\ll _X^\ze = \ix_X\xd^\ze + \xd^\ze\ix_X$ for some exterior derivative
$\xd^\ze$ even on the Grassman algebra $\zF(\zp)$.
        However, we can always find
                $$ \xd_r^\ze,\ \xd_l^\ze \colon \otimes ^1(\zp) \rightarrow
\otimes ^2(\zp)
              $$
        such that
                $$ \ll_X^\ze \zm - \xd_l^\ze\ix_X\zm = \ix^1_X\xd^\ze_r\zm =
-\ix^2_X\xd^\ze_l\zm.
                                                                             $$
        We have, in local coordinates,
                $$\align
        \xd^\ze_r (f_ke_*^k)&= c_{ij}^k f_k e_*^j\otimes e_*^i + \zs^a_i
\frac{\partial f_k}{\partial x^a} e_*^i\otimes e_*^k  -  \zs^a_i
\frac{\partial f_k}{\partial x^a} e_*^k\otimes e_*^i , \\
        \xd^\ze_l (f_ke_*^k)&=- c_{ij}^k f_k e_*^i\otimes e_*^j  -\zs^a_i
\frac{\partial f_k}{\partial x^a} e_*^k\otimes e_*^i +  \zs^a_i
\frac{\partial f_k}{\partial x^a} e_*^i\otimes e_*^k .
                                                  \tag \label{F1.45}\endalign$$
        It follows also that
                $$\align
        \xd_r^\ze(f\zm)&=  f\xd_r^\ze(\zm) + \xd_r^\ze(f) \otimes \zm
-\zm\otimes \xd_l^\ze(f),\\
                \xd_l^\ze(f\zm)&=f\xd_l^\ze(\zm) +\xd_l^\ze(f) \otimes \zm
-\zm\otimes \xd_r^\ze(f),\\
                                                 \tag \label{F1.46}\endalign$$
        and $\xd_l^{\ze^+} = \xd_r^\ze$.

The following theorem is the general algebroid version of Theorem~15~d) of
[\Ref{GU2}].

                \claim \c{t}{Theorem}{}                          \label{C1.17}
        There is a unique derivation $\xd^\ze_l \colon \otimes ^1(\zp)
\rightarrow \otimes ^2(\zp)$ such that
                $$ \ll_X^\ze \zm = \xd_l^\ze \ix_X\zm - \ix^2_X\xd_l^\ze \zm, \
\ \zm \in \otimes ^1(\zp),\ X\in \otimes ^1(\zt).
                                                                              $$
        Moreover,
                $$ \vtb(\xd^\ze_l\zm) = -\ll_{\vtb(\zm)} \zL^\ze.
                                                                              $$
                \endclaim
                \proof
        We have
                $$\align
        \vtb(\xd^\ze_l(f_k e_*^k)) &= -c^k_{ij}f_k \partial _{\zx_i}\otimes
\partial _{\zx_j} -\zs^a_i \frac{\partial f_k}{\partial x^a} \partial
_{\zx_k}\otimes \partial _{\zx_i} + \zr^a_i\frac{\partial f_k}{\partial
x^a} \partial _{\zx_i}\otimes \partial _{\zx_k}\\
        &= - \ll_{f_k\partial _{\zx_k}}\left( c^k_{ij} \zx_k \partial
_{\zx_i}\otimes \partial _{\zx_j} + \zr^a_i \partial _{\zx_i}\otimes
\partial _{x^a} - \zs^a_i \partial _{x^a}\otimes \partial _{\zx_j}\right)\\
        &= - \ll_{\vtb(\zm)}\zL^\ze .
                        \endalign $$
                \endproof
        {\bf Remark.} It is well known that for skew-symmetric $\ze$ we have
$\xd^\ze_r = \xd^\ze_l= \xd^\ze$ and, consequently, it can be extended to
a derivation $\xd^\ze$ of the Grassman algebra $\zF(\zp)$ ([\Ref{KS-M}]).
Then, $\xd^\ze\circ \xd^\ze =0$ if and only if $\zL^\ze$ is Poisson, i.e.,
$\ze$ is a Lie algebroid structure.

        \hl1{Bialgebroids}
        Let $\ze $ and $\widetilde{\ze}$ be algebroid structures on bundles
$\zt\colon E\rightarrow M$ and $\zp\colon E^\*\rightarrow M$ respectively.
Following the ideas of Xu and Mackenzie ([\Ref{MX1}], Theorem 6.2), we
call the pair $(\ze,\widetilde{\ze})$ a bialgebroid if $\ze$, regarded as
a vector bundle morphism
                $$      \CD
                \sT^\*E @() \L{\ze} @(1,0) @() \L{\ezp_E} @(0,-1) & \sT E^\* @()
\L{\sT\zp} @(0,-1) \\
                E   @() \L{\ze_r} @(1,0) & \sT M
                        \endCD \        ,                    \tag \label{F2.1}$$
is a morphism of algebroids $\cT^\*\widetilde{\ze}$ and $\cT
\widetilde{\ze}$, the cotangent and the tangent lifts of the algebroid
$\widetilde{\ze}$.  This means exactly that the dual relation $\ze^\*$
between the corresponding dual bundles
                $$      \CD
                \sT E  @() \L{\ezt_E} @(0,-1) & \sT E @()
\L{\sT\zt} @(0,-1) @() \L{\ze^\*} @(-1,0)\\
                E   @() \L{\zt} @(1,0) & \sT M
                        \endCD \        ,                 \tag \label{F2.2}$$
        is a Leibniz relation between the Leibniz tensors
$\dt\zL_{\tilde{\ze}}$ on both sides, i.e., this relation is a
coisotropic submanifold of $\sT E\times \sT E$ with respect to the Leibniz
tensor $\dt\zL_{\tilde{\ze}}\times (-\dt\zL_{\tilde{\ze}})$.
        Let us assume that, in local coordinates, $\ze$ has the form \Ref{F1.4}
and
                $$ \zL_{\tilde{\ze}} = {\widetilde{c}}_k^{ij}y^k \partial
_{y^i}\otimes \partial _{y^j} + {\widetilde{\zr}}^{ia} \partial _{\zx_i}
\otimes \partial _{x^b} - \widetilde{\zs}^{aj} \partial _{x^a} \otimes
\partial _{y^j}.
                                                            \tag \label{F2.3}$$
        It is easy to calculate that

                $$\multline
        \dt(\zL_{\tilde{\ze}}) = \left( \frac{\partial
\widetilde{c}^{ij}_k}{\partial x^a} y^k\dot x^a + \widetilde{c}^{ij}_ky^k
\dot y^k\right) \partial_{\dot y^i} \otimes \partial_{\dot y^j} +
\frac{\partial \widetilde{\zr}^{ib}}{\partial x^a }\dot x^a \partial_{\dot
y^i} \otimes \partial_{\dot x^b} -\frac{\partial
\widetilde{\zs}^{aj}}{\partial x^b}\dot x^b \partial_{\dot x^a} \otimes
\partial_{\dot y^j} \\
        + \widetilde{c}^{ij}_ky^k (\partial_{\dot y^i} \otimes \partial_{ y^j}
+ \partial_{ y^i} \otimes \partial_{\dot y^j})    + \widetilde{\zr}^{ib}
(\partial_{\dot y^i} \otimes \partial_{x^b} + \partial_{ y^i} \otimes
\partial_{\dot x^b}) - \widetilde{\zs}^{aj}( \partial_{\dot x^a} \otimes
\partial_{\dot y^j} + \partial_{x^a} \otimes \partial_{\dot y^j})
                                \endmultline               \tag \label{F2.4}$$
        and that the relation $\ze^\*$ in $\sT E\times \sT E$ with coordinates
$(x^a, y^i, \dot x^b, \dot y^j, {\bar  x}^{a'}, {\bar y}^{i'}, {\dot {\bar
x}}^{b'}, \dot {\bar y}^{j'})$ is defined by equations
                $$\align
        F^a_0&= {\bar x}^a - x^a =0  \\
        F^a_1&= {\dot x}^a  - \zs^b_k(x){\bar y}^k = 0  \\
        F^j_2&= {\dot y}^j - c^j_{ik}(x)y^i{\bar y}^k - {\dot {\bar y}}^j\\
        F^a_3&= {\dot {\bar x}}^a - \zr^a_k(x)y^k=0 .
                            \tag \label{F2.5}\endalign$$
        The equations define a coisotropic submanifold  $N$ with respect to
$\zL = (-\dt\zL_{\tilde{\ze}}) \times \dt\zL_{\tilde{\ze}}$ if and
only if the Leibniz bracket of the functions \Ref{F2.5} vanish on $N$. The
hamiltonian vector fields corresponding to these functions are the
following (we make use of the fact that $x^a$ and ${\bar x}^a$ are equal
on $N$):

                $$\align
        X^a_0   &= \widetilde{\zs}^{aj}(\partial _{\dot y^j} + \partial_{\dot
{\bar y}^j} ) \quad ; \\
        X^b_1   &= \left( \widetilde{\zs}^{aj} \frac{\partial
\zs^b_k}{\partial x^a} {\bar y^k} - \frac{\partial
\widetilde{\zs}^{bj}}{\partial x^a} \zs^a_l {\bar y^l}\right)\partial
_{\dot y^j} - \widetilde{\zs}^{bj} \partial _{y^j}
                        + \zs^b_i \widetilde{c}^{ij}_k {\bar y^k} \partial
_{\dot{\bar y}^j}  + \zs^b_i \widetilde{\zr}^{ia} \partial _{\dot {\bar
x}^a} \quad ; \\
        X^b_3   &= \left( \widetilde{\zs}^{aj} \frac{\partial \zr^b_k}{\partial
x^a} { y^k} -  \zr^b_i \widetilde{c}^{ij}_k { y^k}\right)\partial _{\dot
y^j} + \widetilde{\zs}^{bj} \partial _{\bar y^j}
                        + \frac{\partial \widetilde{\zs}^{bj}}{\partial x^a}
\zr^a_l {` y^l} \partial _{\dot{\bar y}^j}  - \zr^b_i \widetilde{\zr}^{ia}
\partial _{\dot {\bar x}^a} \quad ; \\
        X^j_2   &= \left( \frac{\partial \widetilde{c}^{jl}_k}{\partial x^a}
\zs^a_i y^k\bar y^i +\widetilde{c}^{jl}_k c^k_{si} y^s{\bar y^i} +
\widetilde{c}^{jl}_k \dot {\bar y}^k  + \widetilde{\zs}^{al}\frac{\partial
c^j_{ik}}{ \partial x^a} y^i {\bar y^k} - c^j_{ik} \widetilde{c}^{il}_s
y^s {\bar y^k}  \right) \partial _{\dot y^l}\\
                & + \left( \frac{\partial \widetilde{c}^{jl}_k}{\partial x^a}
\zr^a_i y^i{\bar y^k} + \widetilde{c}^{jl}_k \dot {\bar y}^k + c^j_{ik}
\widetilde{c}^{kl}_s y^i {\bar y^s} \right) \partial _{\dot {\bar y}^l}
                + \left( \frac{\partial \widetilde{\zr}^{jb}}{\partial x^a}
\zs^a_k {\bar y^k} - \widetilde{\zr}^{ib} c^j_{ik}{\bar y ^k} \right)
\partial _{\dot x^b}\\
        &+ \left( \frac{\partial
\widetilde{\zr}^{jb}}{\partial x^a} \zr^a_k { y^k} + \widetilde{\zr}^{kb}
c^j_{ik}{\bar y ^i} \right) \partial _{\dot {\bar x}^b}
                 + \widetilde{c}^{ji}_k y^k \partial _{y^i} +
\widetilde{\zr}^{jb} \partial _{x^b} + \widetilde{c}^{ji}_k {\bar y^k}\partial
_{\bar y^i} + \widetilde{\zr}^{jb} \partial _{\bar x^b} \ .
                                                 \tag \label{F2.6}\endalign$$
        One can easily see that $F^a_0$ commute, with respect to
$\{\,,\,\}_\zL$, with all defining functions, so that we have three types
of functions left. Moreover, $X_1^b(F^a_1) =0$, $X_3^b(F^a_3) =0$ and we
get the following seven non-trivial equations defining bialgebroids,
corresponding respectively to  $X_1^b(F^a_3) =0$, $X_3^b(F^a_1) =0$,
$X_3^b(F^j_2) =0$, $X_2^j(F^b_1) =0$, $X_1^b(F^j_2) =0$, $X_2^j(F^b_3) =0$,
and $X_2^j(F^i_2) =0$.

                \claim \c{t}{Theorem}{}                          \label{C2.1}
        The tensors $\zL_\ze$ and $\zL_{\tilde{\ze}}$ as in \Ref{F1.5}
 and \Ref{F2.3} constitute a bialgebroid structure if and only if the
following equations are satisfied
        \list
        \item"(1)" $\zs^b_i \widetilde{\zr}^{ia} + \widetilde{\zs}^{bi}\zr^a_i
=0$\medskip
        \item"(2)" $\zr^b_i \widetilde{\zr}^{ia} + \widetilde{\zs}^{bi}\zs^a_i
=0$\medskip
        \item"(3)" $\dfrac{\partial \zr^b_k}{\partial x^a} \widetilde{\zs}^{aj}
- \dfrac{\partial \widetilde{\zs}^{bj}}{\partial x^a}\zr^a_k -
\widetilde{c}^{ij}_k\zr^b_i - c_{ki}^j \widetilde{\zs}^{bi} = 0$\smallskip
        \item"(4)"$ -\dfrac{\partial \zs^b_j}{\partial x^a} \widetilde{\zr}^{ka}
+ \dfrac{\partial \widetilde{\zr}^{kb}}{\partial x^a}\zs^a_j -
\widetilde{c}^{ki}_j \zs^b_i - c_{ij}^k \widetilde{\zr}^{ka} = 0$\smallskip
        \item"(5)" $\dfrac{\partial \zs^b_k}{\partial x^a} \widetilde{\zs}^{aj}
- \dfrac{\partial \widetilde{\zs}^{bj}}{\partial x^a}\zs^a_k -
\widetilde{\zs}^{ij}_k\zs^b_i + c_{ik}^j \widetilde{\zs}^{bi} = 0$
\smallskip
        \item"(6)" $ -\dfrac{\partial \zr^b_k}{\partial x^a}
\widetilde{\zr}^{ja} + \dfrac{\partial \widetilde{\zr}^{jb}}{\partial
x^a}\zr^a_k - \widetilde{\zs}^{ij}_k\zr^b_i - c_{ki}^j
\widetilde{\zr}^{ib} = 0$ \smallskip
        \item"(7)" $\dfrac{\partial \widetilde{c}^{jl}_i}{\partial x^a} \zs^a_k
+\dfrac{\partial {c}_{ik}^j}{\partial x^a} \widetilde{\zs}^{al}
-\dfrac{\partial \widetilde{c}^{jl}_k}{\partial x^a}\zr^a_i
-\dfrac{\partial \widetilde{c}_{ik}^l}{\partial x^a} \widetilde{\zr}^{ja}
+ \widetilde{c}^{jl}_sc^s_{ik} - \widetilde{c}^{sl}_ic^j_{sk} -
\widetilde{c}^{sl}_kc^j_{is} - \widetilde{c}^{js}_ic^l_{sk} -
\widetilde{c}^{js}_kc^l_{is} =0$
        \endlist
                \endclaim

                \claim \c{c}{Corollary}{}                         \label{C2.2}
        The pair $(\ze, \widetilde{\ze})$ constitutes a bialgebroid if and only
if $(\widetilde{\ze}, \ze)$ constitutes a bialgebroid.
                \endclaim
                \proof
        The family of equations (1) -- (7) does not change when we interchange
''tilde'' with ''no tilde''.
                \endproof

        The original definition of a Lie bialgebroid by Mackenzie and Xu
([\Ref{MX1}]) was given in terms exterior derivatives and Lie brackets. In
the case of a general algebroid we have a slight substitute of the
exterior derivative only, but it is enough to get the full analogy.

                \claim \c{t}{Theorem}{}                           \label{C2.3}
        A pair $(\ze, \widetilde{\ze})$ of algebroid structures on $E$ and
$E^\*$ respectively, constitutes  a bialgebroid if and only if
                $$ \xd^{\tilde{\ze}}_l[X,Y]_\ze = [ \xd^{\tilde{\ze}}_lX,Y]_\ze +
[X, \xd^{\tilde{\ze}}_lY]_\ze
                                                             \tag \label{F2.7}$$
        for all $X,Y \in \otimes ^1(\zt)$, where the brackets (''the Schouten
brackets'') are defined by the formulae
                $$\align
        [X\otimes Y, Z]_\ze&= [X,Z]_\ze\otimes Y + X\otimes [Y,Z]_\ze   \\
        [X, Y\otimes Z]_\ze     &=[Z,X]_\ze\otimes Y + X\otimes [Z,Y]_\ze .
                                                    \tag \label{F2.8}\endalign$$
                \endclaim
                \proof
        We shall show, in local coordinates, that \Ref{F2.7} is equivalent to
equations (1)--(7). To reduce the problem to the case $X=e_i$, $Y=e_j$, we
find the relation of \Ref{F2.7} with the structure of
$\cC^\infty(M)$-module in $\otimes ^1(\zt)$. Replacing in \Ref{F2.7} $Y$
by $fY$ we get

                $$\align
        \xd^{\tilde{\ze}}_l\left([X,f]_\ze Y + f[X,Y]_\ze\right)&=
[\xd^{\tilde{\ze}}_lX, f]^2_\ze \otimes Y + Y\otimes [\xd^{\tilde{\ze}}_l
X,f]^1_\ze + f[\xd^{\tilde{\ze}}_l X, Y]_\ze  \\
                        & + [X, f\xd^{\tilde{\ze}}_l Y +
\xd^{\tilde{\ze}}_l(f)\otimes Y - Y\otimes \xd^{\tilde{\ze}}_l(f)]_\ze .
                                                    \tag \label{F2.9}\endalign$$
        (We use the following conventions: $[X,f]_\ze = a^\ze_l(X)(f)$, $[f,X] =
-a^\ze_r(X)(f)$, $[X\otimes Y,f]^2_\ze = [Y,f]_\ze X$, $[X\otimes
Y,f]^1_\ze = [X,f]_\ze Y$, etc.)
        We get furtherly
                $$\multline
        [X,f]_\ze \xd^{\tilde{\ze}}_l Y + \xd^{\tilde{\ze}}_l([X,f]_\ze)
\otimes Y - Y \otimes \xd^{\tilde{\ze}}_r[X,f]_\ze +
f\xd^{\tilde{\ze}}_l[X,Y]_\ze + \xd^{\tilde{\ze}}_l(f)\otimes [X,Y]_\ze
- [X,Y]_\ze \otimes \xd^{\tilde{\ze}}_r(f)\\
        = [X,f]_\ze \xd^{\tilde{\ze}}_lY + f[X,\xd^{\tilde{\ze}}_lY]_\ze
+ [X,\xd^{\tilde{\ze}}_l(f)]_\ze \otimes Y + \xd^{\tilde{\ze}}_l(f)
\otimes [X,Y]_\ze - [X,Y]_\ze \otimes \xd^{\tilde{\ze}}_r(f)\\
        - Y\otimes [X,\xd^{\tilde{\ze}}_r(f)]_\ze +
[\xd^{\tilde{\ze}}_lX,f]^2_\ze\otimes Y + Y\otimes
[\xd^{\tilde{\ze}}_lX,f]_\ze^1 + f[\xd^{\tilde{\ze}}_lX,Y]_\ze
                                \endmultline              \tag \label{F2.10}$$
        and, finally,
                $$\multline
        f\left(\xd^{\tilde{\ze}}_l[X,Y]_\ze - [\xd^{\tilde{\ze}}_lX,Y]_\ze -
[X,\xd^{\tilde{\ze}}_lY]_\ze \right) \\
        = Y\otimes \left(\xd^{\tilde{\ze}}_r[X,f]_\ze + [\xd^ {\tilde{\ze}}_l
X,f]_\ze^1 - [X,\xd^{\tilde{\ze}}_r(f)]_\ze\right) - \left(
\xd^{\tilde{\ze}}_l[X,f]_\ze - [\xd^{\tilde{\ze}}_lX,f]^2_\ze -
[X,\xd^{\tilde{\ze}}_l(f)]_\ze\right) \otimes Y.
                                \endmultline                  \tag \label{F}$$
        This shows that \Ref{F2.7} with $Y=e_i$ is equivalent to equations
                        $$\gather
        \xd^{\tilde{\ze}}_r[X,f]_\ze + [\xd^ {\tilde{\ze}}_l
X,f]_\ze^1 - [X,\xd^{\tilde{\ze}}_r(f)]_\ze =0,    \tag"(a)" \label{F1}\\
        \xd^{\tilde{\ze}}_l[X,f]_\ze - [\xd^{\tilde{\ze}}_lX,f]^2_\ze -
[X,\xd^{\tilde{\ze}}_l(f)]_\ze =0,                      \tag"(b)" \label{F2}
                                        \endgather                          $$
        where $X\in \otimes ^1(\zt)$ and $f\in \cC^\infty(M)$.

        We reduce \Ref{F1} and \Ref{F2} once more, this time with respect to
$X$. Let us put $gX$ instead of $X$. We get from \Ref{F1}
                $$
        \xd^ {\tilde{\ze}}_r\left(g[X,f]_\ze\right) + [g\xd^ {\tilde{\ze}}_lX +
\xd^ {\tilde{\ze}}_lg\otimes X - X\otimes \xd^ {\tilde{\ze}}_rg,f]_\ze^1
-g[X, \xd^ {\tilde{\ze}}_r(f)]_\ze - [g,\xd^ {\tilde{\ze}}_r(f)]_\ze X =0.
                                                                \tag \label{F3}$$
         Hence,
                $$ g\left(\xd^ {\tilde{\ze}}_r[X,f]_\ze + [\xd^ {\tilde{\ze}}_lX,
f]^1_\ze - [X,\xd^ {\tilde{\ze}}_l(f)]_\ze\right) + \left( [\xd^
{\tilde{\ze}}_l g,f]_\ze - [g, \xd^ {\tilde{\ze}}_rf]_\ze \right)X =0
                                                                         $$
        which shows that \Ref{F1} is equivalent to
                $$\gather
        [\xd^ {\tilde{\ze}}_l(g),f]_\ze - [g,\xd^ {\tilde{\ze}}_r(f)]_\ze = 0,
\tag"(c)" \label{Fc}\\
        \xd^ {\tilde{\ze}}_r[e_i,f]_\ze + [\xd^ {\tilde{\ze}}_le_i,f]^1_\ze -
[e_i, \xd^ {\tilde{\ze}}_r(f)]_\ze =0. \tag"(d)" \label{Fd}
                        \endgather                                         $$
        Similarly, \Ref{F2} is equivalent to

                $$\gather
        [\xd^ {\tilde{\ze}}_r(g),f]_\ze - [g,\xd^ {\tilde{\ze}}_l]_\ze = 0,
\tag"(c')" \label{Fc1}\\
        \xd^ {\tilde{\ze}}_l[e_i,f]_\ze + [\xd^ {\tilde{\ze}}_le_i.f]^2_\ze -
[e_i, \xd^ {\tilde{\ze}}_l(f)]_\ze =0. \tag"(d')" \label{Fd1}
                        \endgather                                           $$
        Relaxing now the $\cC^\infty(M)$-module structure in \Ref{F2.7} with
respect to $X$, we get analogously
                        $$\gather
        \xd^{\tilde{\ze}}_r[f,Y]_\ze + [f,\xd^ {\tilde{\ze}}_l
Y]_\ze^1 + [\xd^{\tilde{\ze}}_r(f),Y]_\ze =0,    \tag"(e)" \label{Fe}\\
        \xd^{\tilde{\ze}}_l[f,Y]_\ze - [f,\xd^{\tilde{\ze}}_lY]^2_\ze -
[\xd^{\tilde{\ze}}_l(f),Y]_\ze =0,                      \tag"(f)" \label{Ff}
                                        \endgather                            $$
        which are equivalent to \Ref{Fc}, \Ref{Fc1}, and \Ref{Fe}, \Ref{Ff} with
$Y=e_i$. Finally, in local basis, \Ref{F2.7} is equivalent to the system
of equations
        \list
        \item"(1')" \quad  $[\xd^{\tilde{\ze}}_r(g),f]_\ze -
[g,\xd^{\tilde{\ze}}_l(f)]_\ze = 0$  \medskip
        \item"(2')" \quad   $[\xd^{\tilde{\ze}}_l(g),f]_\ze -
[g,\xd^{\tilde{\ze}}_r(f)]_\ze = 0$ \medskip
        \item"(3')" \quad   $\xd^{\tilde{\ze}}_r[e_i,f]_\ze +
[\xd^{\tilde{\ze}}_l e_i, f]_\ze^1 - [e_i,\xd^{\tilde{\ze}}_r(f)]_\ze =0$
\medskip
        \item"(4')" \quad   $\xd^{\tilde{\ze}}_l[f,e_i]_\ze -
[\xd^{\tilde{\ze}}_l (f),e_i]_\ze - [f,\xd^{\tilde{\ze}}_le_i]_\ze^2 =0$
\medskip
        \item"(5')" \quad   $\xd^{\tilde{\ze}}_r[f,e_i]_\ze +
[\xd^{\tilde{\ze}}_r (f),e_i]_\ze + [f,\xd^{\tilde{\ze}}_re_i]_\ze^1 =0$
\medskip
        \item"(6')" \quad   $\xd^{\tilde{\ze}}_l[e_i,f]_\ze -
[\xd^{\tilde{\ze}}_l e_i, f]_\ze^2 - [e_i,\xd^{\tilde{\ze}}_l(f)]_\ze =0$
\medskip
        \item"(7')" \quad   $\xd^{\tilde{\ze}}_l[e_i,e_j]_\ze +
[\xd^{\tilde{\ze}}_l e_i, e_j]_\ze - [e_i,\xd^{\tilde{\ze}}_le_j]_\ze =0$
\medskip
        \endlist
        Now, it is a direct check that these equations are equivalent to the
system of equations (1)--(7).

                \endproof
        A canonical example of a bialgebroid is given by the following theorem.
                \claim \c{t}{Theorem}{}                              \label{C5}
        Let $\ze$ be a Lie algebroid structure on $\zt\colon E\rightarrow M$ and
let $\zL \in \otimes ^2(\zt)$. Let $\widetilde{\ze}$ be the algebroid
structure on $E^\*$ which corresponds to the linear Leibniz structure
$\dt^\ze\zL$ on $E$. Then the pair $(\ze,\widetilde{\ze})$ is a
bialgebroid.
                \endclaim
                \proof
        From Theorem~\Ref{C1.17} we get
                $$ \vta(\xd^{\tilde{\ze}}_l X) = -[\vta(X),\dt^\ze\zL] =
-\vta([X,\zL]_\ze),
                                                                             $$
        where we used a formula from Corollary~\Ref{C1.10}. It follows that
$\xd^{\tilde{\ze}}_lX = -[X,\zL]_\ze$ and \Ref{F2.7} reads now
                $$
                        [[X,Y]_\ze,\zL] = [[X,\zL]_\ze,Y]_\ze +
[X,[Y,\zL]_\ze]_\ze
                                                                             $$
        which easily follows from the Jacobi identity for $[\,,\,]_\ze$.
                \endproof

        {\bf Remark.} A standard example of a situation as above is provided by
the Lie bialgebroid induced by a Poisson structure on a manifold (cf.
[\Ref{KS-M}], [\Ref{MX1}]).  Moreover, the above theorem shows that a
bialgebroid may be constituted by a Lie algebroid and an algebroid which
is not even skew-symmetric.

        {\bf Example.} An extreme case are bialgebroids over a single point.
It means exactly that we have algebra structures $[\,,\,]$ and
$[\,,\,]'$ on a finite dimensional vector space $E$ and on its dual
$E^\*$ respectively. (An algebra structure on $E$ is a bilinear operation
$[\,,\,] \colon E\times E\rightarrow E$.) They form a bialgebroid (or,
simply, bialgebra) if and only if
                $$ \xd'[X,Y] = [\xd'X,Y] + [X,\xd'Y]
                                                                            $$
for all $X,Y \in E$, where $\xd'\colon E\rightarrow E\otimes E$ is the
dual map to $[\,,\,]' \colon E^\* \otimes E^\* \rightarrow E^\*$.  In the
case of Lie algebra we recognize the definition of a Lie bialgebra.

        \hl1{The algebroid of a linear connection on $\sT M$}
                Important examples of  algebroids which are not skew-symmetric
are provided  by linear connections on a tangent bundle.
                Let a linear connection be given on $\sT M$. It can be
represented by the covariant derivative
                $$ (X,Y)\mapsto \nabla_XY,   \qquad X,Y\in\otimes ^1(\ezt_M)
                                                            \tag \label{F3.1}$$
        or, equivalently, by the horizontal projection
                        $$ P_h\colon \sT\sT M \rightarrow \sT\sT M
                                                                         $$
        or by the vertical projection
                $$      P_v\colon \sT\sT M \rightarrow \sT\sT M  .
                                                                           $$
                The linearity of the connection implies that $P_h, P_v$ are
double vector bundle morphisms.

        Let $(x^a,{\dot x}^b, {x'}{}^c, {\dot x'}{}^d)$ be a coordinate system on
$\sT\sT M$. We have
                $$ \gather (x^a,{\dot x}^b, {x'}{}^c, {\dot x'}{}^d) \circ P_v
= (x^a,{\dot x}{}^b, 0, {\dot x'}{}^d + \zG^d_{ba}(x)\dot x^a x'{}^b), \\
(x^a,{\dot x}^b, {x'}^c, {\dot x'}{}^d) \circ P_h = (x^a,{\dot x}^b, {x'}{}^c,
- \zG^d_{ba}(x)\dot x^a x'{}^b), \\
                \nabla_X Y = \frac{\partial Y^a}{\partial x^b}X^b
\frac{\partial }{\partial x^a} + \zG^a_{cb} Y^b X^c \frac{\partial }{\partial
x^a}.
                                              \tag \label{F3.2} \endgather  $$

                \claim \c{t}{Theorem}{}                           \label{C3.1}
                There is a unique algebroid structure $\ze$ on $\sT M$ such that
                $$   \ll^\ze_X Y = \nabla _XY.
                                                             \tag \label{F3.3}$$
                For this algebroid  $a^\ze_r = \operatorname{id}, \ a_l^\ze = 0$.

        Conversely, any algebroid structure $(\sT M,\ze)$ such that $a^\ze_r =
\operatorname{id}, \ a_l^\ze = 0$ is the algebroid of a linear connection
on $\sT M$
                \endclaim
                \proof
        We define a bracket $[\,,\,]_\ze$ of vector fields by the formula
                $$ [X,Y]_\ze = - \nabla_Y X.
                                                         \tag \label{F3.4}$$
        It is obvious that it satisfies the condition \Ref{F1.6}  with
$a_r^\ze = \operatorname{id}$ and $a^\ze_l = 0$:
                $$ \align
                [fX,gY]_\ze &= -\nabla_{gY}(fX) = -gf\nabla_Y X -gY(f) X \\
                     &= gf [X,Y]_\ze - gY(f) X \\
              &= gf [X,Y] - g a^\ze_r(Y)(f)X + f a^\ze _l (X)(g) Y.
                                                                  \endalign  $$

        To prove the converse it is enough to notice that $a^\ze_r =
\operatorname{id}, \ a_l^\ze = 0$ implies the following form of $\ze
\colon \sT^\* \sT M \rightarrow  \sT \sT^\* M$:
                $$ (x^a,p_b, {\dot x}^c, \dot p_d) \circ  \ze = (x_a, \zp_b,0,
f_d - \zG^a_{db}{\dot x}^b \zp_a),
                                                             \tag \label{F3.5}$$
         where $(x^a, {\dot x}^b, f_c,\zp_d)$ are coordinates in $\sT^\* \sT M$.

                \endproof

        Let $P_h$ be the horizontal projection of a linear connection. The
formula
                $$      P^t_h = \ezk_M\circ P_h\circ \ezk_M
                                                         \tag \label{F3.6}$$
        defines the horizontal lift of the {\it transposed linear connection}
on $\sT M$. The connection is called {\it symmetric} if it is equal to the
transposed connection, i.e., if $ P_h = P_h^t$.

        It is well known that for each linear connection on $\sT M$  with the
covariant derivative $\nabla $ there exists the {\it dual connection} on
$\sT^\* M$, with the covariant derivative $\nabla^+$, such that
                $$ X(\langle \zm, Y\rangle ) = \langle \nabla^+_X \zm, Y\rangle
+ \langle \zm, \nabla_X Y\rangle , \qquad \zm\in \otimes^1 (\ezp_M), \ X,Y \in
\otimes^1 (\ezt_M).
                                                                      $$

                \claim \c{t}{Theorem}{}                    \label{C3.3}
        Let $(\sT M, \ze)$ be the algebroid structure of a linear connection
with the horizontal projection $P_h$.
        The following relations are satisfied:
        \list
        \item"(a)" $\ze = (P^t)^+_v \circ \eze_M$, \ where  $\eze_M =
\eza_M^{-1}$ is the canonical algebroid on $\sT M$ and  $(P^t)^+_v$ is the
vertical projection of the dual to the transposed   connection,
        \medskip
        \item"(b)" $\widetilde{\zL_\ze} = (P^t)^+_v \circ \widetilde{\zL_M}$,
where $\zL_M$ is the canonical Poisson tensor on $\sT^\*M$,
        \medskip
        \item"(c)" $\dt^\ze = P_h \circ \dt$, i.e., $\dt^\ze$ is the horizontal
lift,
        \medskip
\item"(d)" $\ll_X^\ze \zm = \nabla_X^+\zm $ \ for $\zm\in \otimes
^1(\ezp_M)$,
        \medskip
        \item"(e)" $\xd_r^\ze f =\xd f$, $\xd^\ze_l f =0$ for $f\in
\cC^\infty(M)$,
        \medskip
        \item"(f)"  $\xd^\ze_r \zm = \nabla^+ \zm$ for  $\zm \in \otimes
^1(\ezp_M)$.
        \endlist
                \endclaim
                \proof
        \ (a) \ In local coordinates,
                        $$\align
                        (x^a, p_b, \dot x^c, \dot p_d)\circ (P^t)^+_v&=  (x^a,
p_b, 0, \dot p_d -  (\zG^t)^{a}_{bd} p_a\dot x^b ) \\
                        &=  (x^a, p_b, 0, \dot p_d -  \zG^{a}_{db} p_a\dot x^b
), \\
          (x^a, p_b, \dot x^c, \dot p_d)\circ (P^t)^+_v\circ \eze_M&= (x^a,
\zp_b, 0, f_d -  \zG^{a}_{db} \zp_a\dot x^b )     \tag \label{F3.7}\endalign$$
and the equality follows from \Ref{F3.5}.

        \qquad (b) \  We have from (a)
                $$ \widetilde{\zL_\ze} =  \ze\circ \cR_{\ezt_M} =
(P^t)^+_v\circ \eze_M \circ \cR_{\ezt_M} = (P^t)^+_v \circ \widetilde{\zL_M}.
\tag \label{F3.8}$$

        \qquad (c) \ It follows from the formula \Ref{F1.18} that
                $$ \align
                \dt^\ze(X) &= X^a \partial _{x^a} - \zG^a_{cb} \dot x^b X^c \partial
_{\dot x^a} = P_h (X^a\partial _{x^a})\\
                        &=  P_h(X^a\partial _{x^a} + \frac{\partial
X^a}{\partial x^b} \dot x^b \partial _{\dot x^a}) = P_h \circ  \dt X .
                                                   \tag \label{F3.9} \endalign$$
                \qquad (d) \  By the definition of the Lie derivative,
                $$ \zi (\ll^\ze_X(\zm)) = \dt^\ze X(\zi(\zm)) = (P_h\circ \dt X)
(\zi(\zm)).
                                                           \tag \label{F3.10}$$
                Since $\zi(\zm) = \zm_a\dot x^a$ and $P_h\circ \dt X = X^a
\partial _{x^a} - \zG^a_{cb} \dot x^b X^c \partial _{x^a}$, we get
                $$ (P_h\circ \dt X) (\zi(\zm)) = \frac{\partial \zm_b}{\partial
x^a} X^a \dot x^b - \zG^a_{cb} \dot x^b X^c\zm_a
                                                           \tag \label{F3.11}$$
          and, consequently,
        $$ \ll^\ze_X (\zm) =  \frac{\partial \zm_b}{\partial x^a} X^a
\xd x^b - \zG^a_{cb}  X^c\zm_a \xd^b = \nabla^+_X\zm .
                                                                              $$

        \qquad (e) \
                                        $$\align
        \langle \xd^\ze_r(f), X\rangle &= a^\ze_r(X) (f) = X(f) = \langle \xd f,
x\rangle , \\
                       \langle \xd^\ze_l(f), X\rangle&= a^\ze_l(X) (f) =0.
                                                \tag \label{F3.12}\endalign$$

                \qquad (f) \ It follows from the definition of the exterior
derivatives and from  (d) that              $$\align
        \ix^1_X \xd_r^\ze       &= \ll^\ze_X \zm - \xd^\ze_l \ix_X\zm =
\ll^\ze_X \zm  \\
                        &= \nabla^+_X\zm = \ix^1_X \nabla ^+\zm .
                                                              \endalign     $$
                                                \endproof

        \hl1{Metric connections}
        In this section we give an interpretation of the Levi-Civita connection
in terms of algebroids.
                Let $g$ be a contravariant metric tensor on $M$ and let
$\widetilde{g} \colon \sT^\* M \rightarrow \sT M$ be the corresponding
isomorphism of vector bundles. The tensor $\dt g$ on $\sT M$ is linear and
defines an algebroid $(\sT^\* M, \ze)$, i.e., $\zL_{\ze} = \dt g$. In local
coordinates,
                $$ \gather
                g= g^{ab} \partial _{x^a}\otimes \partial _{x^b} ,\\
             \dt g = \zL_{\ze} = \frac{\partial g^{ab}}{\partial x^c} \dot x^c
\partial _{\dot x^a}\otimes \partial _{\dot x^b} + g^{ab}\left( \partial
_{x^a}\otimes \partial _{\dot x^b} + \partial _{\dot x^a}\otimes \partial
_{x^b}\right),\\
        [\zm,\zn]_{\ze} = \frac{\partial g^{ab}}{\partial x^c} \zm_a\zn_b \xd x^c
+ g^{ab} (\frac{\partial \zm_c}{\partial x^a} \zn_b + \zm_a \frac{\partial
\zn_c}{\partial x^b}) \xd x^c .
                                             \tag \label{F3.13}\endgather $$
        It follows that $a_l^{\ze} = \widetilde{g}$ and $a_r^{\ze} = -
\widetilde{g}$.

        Since $\widetilde{g}$ is an isomorphism, $\widetilde{g}^{-1}_*
\zL_{\ze}$ is a linear Leibniz structure on $\sT^\* M$. It
defines an algebroid on $\sT M$ with the bracket
                $$ [X,Y]_{\ze(g)} = \widetilde{g}\left([\widetilde{g}^{-1}X,
\widetilde{g}^{-1} Y]\right).
                                                     \tag \label{F3.14}$$
The left anchor is $a^{\ze(g)}_l = \operatorname{id}$, the right anchor is
$a^{\ze(g)}_r = -\operatorname{id}$, and  coordinates of the algebroid
bracket are  given by the formulae
                $$   ([X,Y]_{\ze(g)})^a = \frac{\partial  X^a}{\partial x^b}
Y^b + \frac{\partial Y^a}{\partial x^b} X^b + g^{ab}\left(\frac{\partial
g_{bc}}{\partial x^d} + \frac{\partial g_{bd}}{\partial x^c} - \frac{\partial
g_{cd}}{\partial x^b}\right) X^c Y^d.
                                                         \tag \label{F3.b1}$$

        The formula
                $$   \zL_{\ze_g} = \frac{1}{2}( \zL_M -  \widetilde{g}^{-1}_*
\dt g)
                                                           \tag \label{F3.15}$$
         defines then an algebroid structure $(\sT, \ze_g)$ with the anchors
$a_r^{\ze_g} = \operatorname{id}$, $a_l^{\ze_g} =0 $, i.e., $(\sT M, \ze_g)$
is the algebroid of a linear connection on $\sT M$. The covariant
derivative of this connection is given by the formula
                $$ \align
                \nabla_XY &= -[Y,X]_{\ze_g}= \frac{1}{2} \left( [X,Y] +
[Y,X]_{\ze(g)}\right)\\
           &=\frac{\partial Y^a}{\partial x^b} X^b \partial _{x^a} +
\frac{1}{2} g^{ab}\left(\frac{\partial g_{bc}}{\partial x^d} + \frac{\partial
g_{bd}}{\partial x^c} - \frac{\partial g_{cd}}{\partial x^b}\right)X^c
Y^d\partial _{x^a} .
                                                \tag \label{F3.b2} \endalign$$
        We recognize $\nabla _X Y$ as a covariant derivative with respect to the
Levi-Civita connection of the metric $g$.

        \hl1{Lifting processes on differentiable grupoids}
        In [\Ref{MX2}] a procedure of lifting of multiplicative vector fields on
Lie grupoids is described. This procedure, what the authors admit
implicitely, has not to do with the whole grupoid structure, but rather
with the structure of a fibration only. In other words, it is just the
standard complete lift procedure, but applied to specific vector fields.
The general scheme is based on the following theorem.

                \claim \c{t}{Theorem}{}                            \label{C3.4}
        If $\za\colon P\rightarrow B$ is a fibration, $\zg\colon B\rightarrow P$
is its section, and $X$ is a projectable vector field on $P$, tangent to
$\zg(B)$, then the complete lift $\dt X$ is a vector field on $\sT P$,
tangent to the subbundle $E\subset \sT_{\zg(B)} P$ of $\za$-vertical
vectors.
                \endclaim
                \proof
        In a neighbourhood of a point $p\in \zy(B)$, we can choose coordinates
$(x^a,f^i)$ such that $(x^a)$ are coordinates near $\za(p)\in B$ and
$(f^i)$ are coordinates in fibers, vanishing on $\zg(B)$.

        Let us write
                $$ X(x,f) = g^a(x) \partial _{x^a} + h^i(x,f) \partial _{f^i},
                                                                              $$
        where $h^i(x,0) =0$. Thus
                $$ \dt X = g^a(x) \partial _{x^a} + \frac{\partial g^a}{\partial
x^b}(x) \dot x^b \partial _{\dot x^a} + h^i(x,f) \partial _{f^i} +
\frac{\partial h^i}{\partial x^a} (x,f)\dot x^a \partial _{\dot f^i} +
\frac{\partial h^i}{\partial f^j}(x,f) \dot f^j\partial _{\dot f^i}.
                                                                          $$
        The subbundle $E\subset \sT P$ is described by the conditions $\dot
x^a = 0,\ f^i =0$, so that $\dt X$ is clearly tangent to $E$ and
                $$ \dt X|_{E} = g^a(x)\partial _{x^a} + \frac{\partial
h^i}{\partial f^j}(x,0)\dot f^j \partial _{\dot f^i}.
                                                                     $$
                \endproof

What we use to lift multiplicative vector fields on a Lie grupoid $\za,\zb
\colon G\rarrows B $ to vector fields on the bundle $A(G)$ of the
corresponding Lie algebroid is only the $\za$--fibration structure over
the base $B$ and the fact that multiplicative vector fields have the
desired properties (are star-vectors in the terminology of [\Ref{MX2}]).

        On the other hand, one can use the above lift to obtain the lift
$\dt^\ze$ for the corresponding Lie algebroid. This lift for sections of
$A(G)$ was introduced in [\Ref{MX2}]  by the formula
                $$ \dt^\ze X (\zi(\zm)) = \zi(\ll^\ze_X\zm)
                                                    \tag \label{F3.16}$$
(cf. our Theorem~\Ref{C1.15}). Since we want to present this procedure in
the whole generality, we have to start with a more general object than a
Lie grupoid.

        A {\it pre-Lie grupoid} is, roughly speaking, an object we obtain by
relaxing the associativity condition in the definition of a Lie grupoid.
Having submersions $\za, \zb \colon G\rightarrow B$ onto the manifold $B$
of units and the  inclusion map $e\colon B\rightarrow G$, we define the
vector bundle $\zt\colon A(G)\rightarrow B$ as usual, to be the inverse
image of the vertical bundle $\ss V^\za G\rightarrow G$ across the
embedding $e \colon B\rightarrow G$. We shall consider, for simplicity,
$B$ to be just a submanifold of $G$, so that $A(G)= \ss V^\za_B G$.
Starting with a section $X\colon B\rightarrow A(G) \subset \sT G$, we
define the right and left prolongations of $X$ to a vector field on $G$
                $$ \overrightarrow{X}(g) =\sT(R_g)X(\zb g) \ \ \text{and} \ \
\overleftarrow{X}(g) =\sT(L_g)T(i)X(\za g),
                                                        \tag \label{F3.17}$$
        where $R_g, L_g$ are the right and left translations, and $i \colon
G\rightarrow G$ is the inverse mapping. Now, we define a bracket on
sections of $A(G)$ putting
                $$ [X,Y]_\ze(p) = [\overrightarrow{X}, \overrightarrow{Y}](p)
                                                         \tag \label{F3.18}$$
         for $p\in B$. Let us note that $\overrightarrow{X} $ and
$\overleftarrow{X} $ are no longer invariant vector fields, since the
pre-grupoid product is not assumed to be associative. Hence, the bracket
$[\overrightarrow{X}, \overrightarrow{Y}]$ is no longer the
right-prolongation of any element of $\otimes ^1(\zt)$. On the other hand,
the definition \Ref{F3.18} makes sense and we have the following.

                \claim \c{t}{Theorem}{}                          \label{C3.5}
        The bracket $[\,,\,]_\ze$ defines a skew-symmetric algebroid structure
on the bundle $\zt\colon A(G)\rightarrow B$ with the anchor
                $$ a\colon A(G)\longmapsto \sT B,\ \ \ a= \sT\zb|_{A(G)}.
                                                                             $$
                \endclaim
                \proof
        The bracket $[\,,\,]_\ze$ is obiously skew-symmetric and, for $f\in
\cC^\infty (B)$, we have
                $$\multline
        [X,fY]_\ze(p) = [\overrightarrow{X}, \overrightarrow{fY}](p) =
[\overrightarrow{X}, (f\circ \zb)\overrightarrow{Y}] (p) \\=
f(\zb(p))[\overrightarrow{X}, \overrightarrow{Y}] (p) +
\overrightarrow{X}(f\circ \zb)(p) \overrightarrow{Y}(p) = f[X,Y]_\ze(p) +
((\sT\zb(X))(f))(p)Y(p).
                                \endmultline                                 $$
                \endproof

                \claim \c{t}{Theorem}{}                            \label{C3.6}
The vector field $\overleftrightarrow{X} = \overrightarrow{X}
+\overleftarrow {X}$ projects to $a(X)$ under both projections, $\za$ and
$\zb$, and it is tangent to $B$--regarded as a submanifold of $G$.
                \endclaim
                \proof
        No difference with respect to the classical case. Let us  mention only,
that in the pre-algebroid case the vector field $\overleftrightarrow{X} $
is, in general,  no longer multiplicative. This shows that the
multiplicativity is not essential for the lifting procedures.
                \endproof
                \claim \c{c}{Corollary}{}                            \label{C3.7}
        The vector field $\dt \overleftrightarrow{X} $ is tangent to the
submanifold $A(G)$ of $\sT G$.
                \endclaim

        Let us denote $\dt \overleftrightarrow{X}|_{A(G)} $ by $\dt ^G X$ -- the
pre-Lie grupoid lift of a section $X$ of $A(G)$.

                \claim \c{t}{Theorem}{}                             \label{C3.8}
        For each section $X\in \otimes ^1{A(G)}$, we have
                $$ \dt^G X = \dt^\ze X.
                                                                             $$
                \endclaim
                \proof
        Since the set of functions $\{\zi(\zm)\colon \zm\in \otimes ^1(\zp)\}$,
where $\zp \colon A^\*(G)\rightarrow B$ is the bundle du;at to $\zt$, Is a
complete set of functions almost everywhere on $A(G)$, it is sufficient, in
view of \Ref{F3.16}, to prove  the equality
                $$ \dt^G X(\zi(\zm)) = \zi(\ll^\ze_X\zm) \ \ \ \text{for all}\ \
\zm\in \otimes ^1(\zp),
                                                                             $$
i.e.,
                $$ \dt^G X(\zi(\zm)) \circ Y = a(X) \langle \zm,Y\rangle  -
\langle \zm, [X,Y]_\ze\rangle
                                                         \tag \label{F3.19}$$
        for all $\zm\in \otimes ^1(\zp)$, $X,Y \in \otimes ^1(\zt)$.

        We can find (at least locally) a function $f$ on $G$ vanishing on $B$
and such that $\xd f$ represents $\zm$, i.e., $\langle \zm,Y \rangle =
\langle \xd f|_B ,Y\rangle $.
        We have,
                $$\multline
        \dt^GX(\zi(\zm))\circ Y = \langle \dt \overleftrightarrow{X} ,
\xd\zi(\zm) \rangle \circ Y = \langle \dt \overleftrightarrow{X} , \dt \xd
f \rangle \circ Y = \dt \langle \overleftrightarrow{X} , \xd f\rangle
\circ Y \\
        = \zi\left(\xd\,\langle \overleftrightarrow{X} , \xd f\rangle \right)
\circ Y = \langle \xd\,\langle \overleftrightarrow{X} , \xd f\rangle ,
Y\rangle  = \langle \overrightarrow{Y}, \xd\,\langle
\overleftrightarrow{X} ,\xd f\rangle |_B = \overrightarrow{Y}(\langle
\overleftrightarrow{X} ,\xd f\rangle )|_B\\
        =\left( \langle [\overrightarrow{Y}, \overleftrightarrow{X} ], \xd
f\rangle  + \overleftrightarrow{X} (\langle \overrightarrow{Y}, \xd
f\rangle )\right)|_B = \langle [\overrightarrow{Y}, \overrightarrow{X}],
\xd f\rangle |_B + a(X)\langle Y, \zm\rangle\\
        = a(X) \langle Y,\zm \rangle  -\langle \zm, [X,Y]_\ze \rangle  =
\langle Y,\ll_X^\ze{\zm}\rangle .
                                \endmultline                           $$
        Here we used the fact that $[\overrightarrow{Y}, \overleftarrow{X}]|_B
=0$. In general, the left and right prolongations do not commute as in the
case of Lie groupoids, but they commute on $B$, and it is sufficient for
our purposes.
                \endproof

        It is clear that relaxing the associativity assumption we get
skew-algebroids as introduced in Section~\Ref{S1}  which are no longer Lie
algebroids.

        {\bf Example.} Let $V$ be a vector space and let $D\colon V\times V
\rightarrow V$ be a skew-symmetric, bi-linear mapping. We define a pre-Lie
group structure on V (a pre-Lie grupoid over a single point) as follows.
Let  $X\star Y$ be the pre-group product given by
        $$X\star Y = X +Y +D(X,Y).
                                                                         $$
        It is clear that $0$ is the unit element and $X\mapsto -X$ is the
inverse mapping. $A(V)$ is canonically identified with $V$. We have the
following formulae for prolongations
                $$\align
        \overrightarrow{X}(Z) &= X + D(X,Z),   \\
        \overleftarrow{X}(Z)&=  -X + D(X,Z),
                                                                 \endalign$$
        and we obtain the bracket
                $$ [X,Y]_\ze = [\overrightarrow{X}, \overrightarrow{Y}](0) =
2D(X,Y).
                                                                            $$
        It is easy to see that this bracket satisfies the Jacobi identity if
and only if '$\star$' is associative. Also the vector field
$\overleftrightarrow{X}(Z) = 2D(X,Z)$ is multiplicative if and only if
'$\star$' is associative.
        \bigskip

        In the above considerations concerning pre-Lie groupoids we implicitly
assumed that the product $q_1q_2$ exists if and only if $\zb(q_1)=
\za(q_2)$. However, relaxing the associativity assumption, we should
probably change also this axiom as shows the following example.

        {\bf Example.} Suppose that in a Lie group $D$ we have chosen a Lie
subgroup $G$ and a `complementary' submanifold such that

        \list
\item $e\in M$ and $M^{-1} = M$,
        \item the composition $G\times M \ni (g,u) \mapsto gu \in D$ is a
diffeomorphism,
        \item for each pair $u,v\in M$ there is a unique $g\in G$, denoted by
$\zf(u,v)$, such that $ugv\in M$.
        \endlist
        We denote the element $u\zf(u,v)v$  by $u \star v $.

        For a given subgroup $G$ one can find such a triple $(D,M,G)$, at least
locally, putting $M = \exp \fr m$, where $\fr m$ is a complementary subspace to
the Lie algebra $\fr g$ of $G$ in the Lie algebra $\fr d$ of $D$.  The reader
easily recognizes some similarities with double Lie groups (e.g. if $M$ is a
subgroup) and quasi-Poisson Lie groups (cf. [\Ref{KS1}] and the concept of
quasi-triple in [\Ref{AKS}]).

        Since every element $d\in D$ has unique decompositions $d=gu =vh$,
where $g,h\in G$ and $u,v \in M$, we shall write $d = (g,u,v,h)$. We have
obvious projections $\za, \zb \colon D\rightarrow G$,
                $$ \za(g,u,v,h) = g, \ \ \ \zb(g,u,v,h) = h.
                                                                            $$
        Now, we define a partial product in $D$ defined for pairs
        $$d_1 =(g_1,u_1,v_1,h_1),\  d_2 =(g_2,u_2,v_2,h_2)$$
        such that
        $$ \zb(d_1)\zf(u_1,u_2) = \zf(v_1,v_2) \za(d_2)
                                                                            $$
        by
                $$ d_1\circ d_2 = (g_1, u_1\star u_2, v_1\star v_2, h_2).
                                                                            $$
        This definition is correct, because
                $$ \multline
        g_1(u_1\star u_2) = g_1 u_1 \zf(u_1,u_2)u_2 = v_1h_1\zf(u_1,u_2)u_2  \\
        = v_1\zf(v_1,v_2)q_2u_2 = v_1\zf(v_1,v_2)v_2h_2 = (v_1\star v_2)h_2.
                                                         \endmultline    $$
        Moreover, $\za(d_1\circ d_2) = \za(d_1)$ and $\zb(d_1\circ d_2) =
\zb(d_2)$. Of course, this partial operation is not associative unless the
$\star$-product is not associative. It is easy to see that elements $d=
(g,e,e,g)$, i.e., $d=g \in G$, form the set of units and that every element
$d= (g,u,v,h)$ has the inverse $i(d)= (g,u^{-1}, v^{-1},h)$.
        If $M$ is also a subgroup, then what we get is exactly the Lie groupoid of
the double group.

        The question what are the algebroids being the infinitesimal versions of
such structures we postpone to a separate paper.

\newcounter\bibno
\def\bb{\bib\key \bibno}

        \makebib
        \bb \label{AKS} \by  A.~Alekseev and Y.~Kosmann-Schwarzbach \pp
1243--1274 \paper Manin pairs and moment maps  \paperinfo preprint, Centre de
Math\'ematique, Ecole Polytechnique, 98-9, 1998
                \endbib

        \bb \label{Co} \by T.~Courant \pp 4527--4536 \paper Tangent Lie
algebroids \yr 1994 \vol 27 \jour J. Phys. A
                \endbib

        \bb \label{Dr} \by V.~G. Drinfeld \pp 1419--1457 \paper Quasi--Hopf
algebras \yr 1990 \vol 1 \jour Leningrad. Math. J.
                \endbib

        \bb \label{Du}\by J.-P. Dufour \pp 55--76 \paper Introduction aux
tissus \yr 1991 \jour S\'eminaire GETODIM
       \endbib

        \bb \label{GU1} \by J.~Grabowski and P.~Urba\'nski \pp 6743--6777
        \paper Tangent lifts of Poisson and related structures
        \yr  1995 \vol 28
        \jour J. Phys. A
                \endbib
        \bb \label{GU2} \by J.~Grabowski and P.~Urba\'nski \pp 447--486 \paper
Tangent and cotangent lifts and graded Lie algebras associated with Lie
algebroids \yr 1997 \vol 15 \jour Ann. Global Anal. Geom.
                \endbib
        \bb \label{GU3} \by J.~Grabowski and P.~Urba\'nski \pp 195--208
\paper Lie algebroids and Poisson-Nijenhuis structures \yr 1997 \vol 40
\jour Rep. Math. Phys.
                \endbib

        \bb \label{KU} \by K.~Konieczna and P.~Urba\'nski \paper Double vector
bundles and duality \paperinfo dg-ga/9710014
                \endbib

        \bb \label{KS1} \by Y.~Kosmann-Schwarzbach \pp 453--489 \paper
Jacobian quasi-bialgebras and quasi-Poisson Lie groups \yr 1992 \vol 132
\jour Contemporary Math.
                \endbib

        \bb \label{KS2} \by Y.~Kosmann-Schwarzbach \pp 1243--1274 \paper
From Poisson algebras to Gerstenhaber algebras \yr 1996 \vol 46 \jour
Ann.~Inst.~Fourier
                \endbib

        \bb \label{KS-M} \by Y.~Kosmann-Schwarzbach and F.~Magri \pp
35--81
        \paper Poisson-Nijenhuis structures
        \yr 1990  \vol 53
        \jour Ann. Inst. Henri Poin\-ca\-r\'e, Phys. Th\'eor.
                \endbib

        \bb \label{Li} \by P.~Libermann \pp 147--162 \paper Lie algebroids
and mechanics \yr 1996 \vol 32 \jour Archivum Mathematicum
                \endbib

        \bb \label{Lo} \by J.--L.~Loday \pp
269--293 \paper Une version non
commutative des alg\`ebres de Lie: les alg\`ebres de Leibniz \yr 1993 \vol
321 \jour Ensegn.~Math.
                        \endbib

        \bb \label{Ma} \by K.~C.~H.~Mackenzie \pp 97--147 \paper Lie
algebroids
and Lie pseudoalgebras \yr 1995 \vol 27 \jour Bull. London Math. Soc.
        \endbib
        \bb \label{MX1} \by K.~C.~H.~Mackenzie and P. Xu  \pp 415--452
        \paper Lie bialgebroids and Poisson groupoids
        \yr 1994  \vol 73
        \jour  Duke Math. J.
                \endbib
        \bb \label{MX2} \by K.~C.~H.~ Mackenzie and P. Xu \pp 59--85 \paper
Classical
lifting processes and multiplicative vector fields \yr 1998
\vol 49 \jour  Quarterly J. Math. Oxford (2)
        \endbib
\bb \label{PT} \by G.~Pidello and W.~Tulczyjew \pp 249--265 \paper
Derivations of differential forms on jet bundles \yr 1987 \vol 147 \jour
Ann. Mat. Pura Appl.
                \endbib
        \bb \label{Pr} \by J.~Pradines \pp 245--248 \paper Th\'eorie de Lie
pour
groupo\"\i des diff\'erentiables. Calcul diff\'erentiel dans la cat\'egorie
de groupo\"\i des ininit\'esimaux \yr 1967 \vol 264 \jour C.R. Acad. Sci.
Paris S\'er. A
        \endbib
        \bb \label{Tu} \by W.~Tulczyjew \pp 101--114 \paper Hamiltonian
systems,
Lagrangian systems, and the Legendre transformation \yr 1974 \vol 14 \jour
Symposia Math.
                \endbib
        \bb \label{Ur} \by P.~Urba\'nski \pp 405--421 \paper Double vector
bundles in classical mechanics \yr 1996 \vol 54 \jour Rend. Sem. Mat.
Univ. Pol. Torino
                \endbib

         \bb \label{We} \by A.~Weinstein \paper Lagrangian mechanics and
grupoids \paperinfo preprint, Univ. of California, Berkeley 1992
                \endbib

        \bb \label{YI} \by K.~Yano and S.~Ishihara
        \book Tangent and Cotangent Bundles
          \publ Marcel Dekker, Inc., New York \yr 1973
          \endbib

\enddocument

%% file: lamstex.tex


\catcode`\@=11
\ifx\amstexloaded@\relax\else
 \errmessage{AmS-TeX must be loaded before LamS-TeX}\fi
\ifx\laxread@\undefined\else\catcode`\@=\active \fi
\def\err@#1{\errmessage{LamS-TeX error: #1}}
\def^^L{\par}
\let\+\tabalign
\def\newcount{\alloc@0\count\countdef\insc@unt}
\def\newdimen{\alloc@1\dimen\dimendef\insc@unt}
\def\newskip{\alloc@2\skip\skipdef\insc@unt}
\def\newmuskip{\alloc@3\muskip\muskipdef\@cclvi}
\def\newbox{\alloc@4\box\chardef\insc@unt}
\let\newtoks\relax
\def\newhelp#1#2{\newtoks#1#1\expandafter{\csname#2\endcsname}}
\def\newtoks{\alloc@5\toks\toksdef\@cclvi}
\def\newread{\alloc@6\read\chardef\sixt@@n}
\def\newwrite{\alloc@7\write\chardef\sixt@@n}
\def\newfam{\alloc@8\fam\chardef\sixt@@n}
\def\newlanguage{\alloc@9\language\chardef\@cclvi}
\def\newinsert#1{\global\advance\insc@unt by\m@ne
  \ch@ck0\insc@unt\count
  \ch@ck1\insc@unt\dimen
  \ch@ck2\insc@unt\skip
  \ch@ck4\insc@unt\box
  \allocationnumber=\insc@unt
  \global\chardef#1=\allocationnumber
  \wlog{\string#1=\string\insert\the\allocationnumber}}
\def\newif#1{\count@\escapechar \escapechar\m@ne
  \expandafter\expandafter\expandafter
   \edef\@if#1{true}{\let\noexpand#1=\noexpand\iftrue}%
  \expandafter\expandafter\expandafter
   \edef\@if#1{false}{\let\noexpand#1=\noexpand\iffalse}%
  \@if#1{false}\escapechar\count@}

\def\Err@#1{\errhelp\defaulthelp@\err@{#1}}
{\catcode`\@=\active
 \edef\next{\gdef\noexpand@{\futurelet\noexpand\next
  \csname at\string@\endcsname}}
 \next
}
\def\at@{\ifcat\noexpand\next a\let\next@\at@@\else
 \ifcat\noexpand\next0\let\next@\at@@\else
 \ifcat\noexpand\next\relax\let\next@\at@@\else
 \let\next@\at@@@\fi\fi\fi\next@}
\def\at@@@{\errhelp\athelp@\err@{Invalid use of @}}
\def\at@@#1{\expandafter
 \ifx\csname\string#1@at\endcsname\relax\let\next@\at@@@\else
 \DN@{\csname\string#1@at\endcsname}\fi\next@}
\def\atdef@#1{\expandafter\def\csname\string#1@at\endcsname}
\newif\iftest@
\def\tagin@#1{\tagin@false
 \DN@##1\tag##2##3\next@{\test@true\ifx\tagin@##2\test@false\fi}%
 \next@#1\tag\tagin@\next@\tagin@false\iftest@\tagin@true\fi}
\let\lkerns@\relax
\def\nolinebreak{\RIfM@\mathmodeerr@\nolinebreak\else
 \ifhmode\saveskip@\lastskip\unskip
 \nobreak\ifdim\saveskip@>\z@\hskip\saveskip@\fi\lkerns@
 \else\vmodeerr@\nolinebreak\fi\fi}
\def\allowlinebreak{\RIfM@\mathmodeerr@\allowlinebreak\else
 \ifhmode\saveskip@\lastskip\unskip
 \allowbreak\ifdim\saveskip@>\z@\hskip\saveskip@\fi\lkerns@
 \else\vmodeerr@\allowlinebreak\fi\fi}
\def\linebreak{\RIfM@\mathmodeerr@\linebreak\else
 \ifhmode\unskip\unkern\break\lkerns@
 \else\vmodeerr@\linebreak\fi\fi}
\let\nkerns@\relax
\def\newline{\RIfM@\mathmodeerr@\newline\else
 \ifhmode\unskip\unkern\null\hfill\break\nkerns@
 \else\vmodeerr@\newline\fi\fi}%
\def\newbox@{\alloc@@4\box\chardef\insc@unt}
\def\newcount@{\alloc@@0\count\countdef\insc@unt}
\def\accentedsymbol#1#2{\expandafter\newbox@\csname\exstring@#1@box\endcsname
 \setbox\csname\exstring@#1@box\endcsname\hbox{$\m@th#2$}%
 \define#1{\copy\csname\exstring@#1@box\endcsname{}}}
\def\rightadd@#1\to#2{\toks@{\\#1}\toks@@\expandafter{#2}\xdef#2{\the\toks@@
 \the\toks@}\toks@{}\toks@@{}}
\def\fontlist@{\\\tenrm\\\sevenrm\\\fiverm\\\teni\\\seveni\\\fivei
 \\\tensy\\\sevensy\\\fivesy\\\tenex\\\tenbf\\\sevenbf\\\fivebf
 \\\tensl\\\tenit}
\def\font@#1=#2 {\rightadd@#1\to\fontlist@\font#1=#2 }
\def\ismember@#1#2{\global\let\Next@ F\let\next@= #2%
 {\def\\##1{\let\nextii@##1\ifx\nextii@\next@\global\let\Next@ T\fi}#1}%
 \test@false\ifx\Next@ T\test@true\fi\let\next@\relax}
\def\FNSS@#1{\let\FNSS@@#1\FN@\FNSS@@@}
\def\FNSS@@@{\ifx\next\space@\def\FNSS@@@@. {\FN@\FNSS@@@}\else
 \def\FNSS@@@@.{\FNSS@@}\fi\FNSS@@@@.}
\atdef@"{\unskip
 \DN@{\ifx\next`\DN@`{\FN@\nextii@}%
  \else\ifx\next\lq\DN@\lq{\FN@\nextii@}%
  \else\DN@####1{\FN@\nextiii@}\fi\fi
  \next@}%
 \DNii@{\ifx\next`\DN@`{\sldl@``}%
  \else\ifx\next\lq\DN@\lq{\sldl@``}%
  \else\DN@{\dlsl@`}\fi\fi\next@}%
 \def\nextiii@{\ifx\next'\DN@'{\srdr@''}%
  \else\ifx\next\rq\DN@\rq{\srdr@''}%
  \else\DN@{\drsr@'}\fi\fi\next@}%
 \FNSS@\next@}
\def\root{%
  \DN@{\ifx\next\uproot\let\next@\nextii@\else
   \ifx\next\leftroot\let\next@\nextiii@\else
   \let\next@\plainroot@\fi\fi\next@}%
  \DNii@\uproot##1{\uproot@##1\relax\FNSS@\nextiv@}%
  \def\nextiv@{\ifx\next\leftroot\let\next@\nextv@\else
   \let\next@\plainroot@\fi\next@}%
  \def\nextv@\leftroot##1{\leftroot@##1\relax\plainroot@}%
  \def\nextiii@\leftroot##1{\leftroot@##1\relax\FNSS@\nextvi@}%
  \def\nextvi@{\ifx\next\uproot\let\next@\nextvii@\else
   \let\next@\plainroot@\fi\next@}%
  \def\nextvii@\uproot##1{\uproot@##1\relax\plainroot@}%
  \bgroup\uproot@\z@\leftroot@\z@
 \FNSS@\next@}
\def\loop#1\repeat{\def\iterate{#1\relax\expandafter\iterate\fi}%
 \iterate\let\iterate\relax}
\def\gloop@#1\repeat{\gdef\iterate@{#1\relax\expandafter\iterate@\fi}%
 \iterate@\global\let\iterate@\relax}
\def\printoptions{\W@{Do you want S(yntax check),
  G(alleys) or P(ages)?^^JType S, G or P, follow by <return>: }\loop
 \read\m@ne to\ans@
 \edef\next@{\def\noexpand\Ans@{\ans@}}\uppercase\expandafter{\next@}%
 \ifx\Ans@\S@\test@true\syntax\else
 \ifx\Ans@\G@\test@true\galleys\else
 \ifx\Ans@\P@\test@true\else
 \test@false\fi\fi\fi
 \iftest@\else\W@{Type S, G or P, follow by <return>: }%
 \repeat}
\expandafter\let\csname A@;\endcsname;
\expandafter\let\csname A@:\endcsname:
\expandafter\let\csname A@?\endcsname?
\expandafter\let\csname A@!\endcsname!
\def\APdef#1{\def\next@{\expandafter\let\csname A@\string#1\endcsname#1}%
 \afterassignment\next@\def#1}
\let\fextra@\,
\def\tdots@{\unskip
 \DN@{$\m@th\mathinner{\ldotp\ldotp\ldotp}\,
   \ifx\next,\,$\else\ifx\next.\,$\else
   \ifx\next;\,$\else
   \expandafter\ifx\csname A@\string;\endcsname\next\fextra@$\else
   \ifx\next:\,$\else
   \expandafter\ifx\csname A@\string:\endcsname\next\fextra@$\else
   \ifx\next?\,$\else
   \expandafter\ifx\csname A@\string?\endcsname\next\fextra@$\else
   \ifx\next!\,$\else
   \expandafter\ifx\csname A@\string!\endcsname\next\fextra@$\else
   $ \fi\fi\fi\fi\fi\fi\fi\fi\fi\fi}%
 \ \FN@\next@}
\def\extrap@#1{%
 \ifx\next,\DN@{#1\,}\else
 \ifx\next;\DN@{#1\,}\else
 \expandafter\ifx\csname A@\string;\endcsname\next\DN@{#1\fextra@}\else
 \ifx\next.\DN@{#1\,}\else\extra@
 \ifextra@\DN@{#1\,}\else
 \let\next@#1\fi\fi\fi\fi\fi\next@}
\def\dotsc{\DN@{\ifx\next;\plainldots@\,\else
 \expandafter\ifx\csname A@\string;\endcsname\next\plainldots@\fextra@\else
 \ifx\next.\plainldots@\,\else\extra@\plainldots@
 \ifextra@\,\fi\fi\fi\fi}%
 \FN@\next@}
\def\keybin@{\keybin@true
 \ifx\next+\else\ifx\next=\else\ifx\next<\else\ifx\next>\else\ifx\next-\else
 \ifx\next*\else\ifx\next:\else
 \expandafter\ifx\csname A@\string;\endcsname\next\else
 \keybin@false\fi\fi\fi\fi\fi\fi\fi\fi}
\def\boldkey#1{\ifcat\noexpand#1A%
  \ifcmmibloaded@{\fam\cmmibfam#1}\else
   \Err@{First bold symbol font not loaded}\fi
 \else
 \let\next=#1%
 \ifx#1!\mathchar"5\bffam@21 \else
 \expandafter\ifx\csname A@\string!\endcsname\next\mathchar"5\bffam@21 \else
 \ifx#1(\mathchar"4\bffam@28 \else\ifx#1)\mathchar"5\bffam@29 \else
 \ifx#1+\mathchar"2\bffam@2B \else\ifx#1:\mathchar"3\bffam@3A \else
 \expandafter\ifx\csname A@\string:\endcsname\next\mathchar"3\bffam@3A \else
 \ifx#1;\mathchar"6\bffam@3B \else
 \expandafter\ifx\csname A@\string;\endcsname\next\mathchar"6\bffam@3B \else
 \ifx#1=\mathchar"3\bffam@3D \else
 \ifx#1?\mathchar"5\bffam@3F \else
 \expandafter\ifx\csname A@\string?\endcsname\next\mathchar"5\bffam@3F \else
 \ifx#1[\mathchar"4\bffam@5B \else
 \ifx#1]\mathchar"5\bffam@5D \else
 \ifx#1,\mathchari@63B \else
 \ifx#1-\mathcharii@200 \else
 \ifx#1.\mathchari@03A \else
 \ifx#1/\mathchari@03D \else
 \ifx#1<\mathchari@33C \else
 \ifx#1>\mathchari@33E \else
 \ifx#1*\mathcharii@203 \else
 \ifx#1|\mathcharii@06A \else
 \ifx#10\bold0\else\ifx#11\bold1\else\ifx#12\bold2\else\ifx#13\bold3\else
 \ifx#14\bold4\else\ifx#15\bold5\else\ifx#16\bold6\else\ifx#17\bold7\else
 \ifx#18\bold8\else\ifx#19\bold9\else
  \Err@{\noexpand\boldkey can't be used with #1}%
 \fi\fi\fi\fi\fi\fi\fi\fi\fi\fi\fi\fi\fi\fi\fi
 \fi\fi\fi\fi\fi\fi\fi\fi\fi\fi\fi\fi\fi\fi\fi\fi\fi\fi}
\def\arabic#1{#1}
\def\alph#1{\count@#1\relax\advance\count@96 \ifnum\count@>122
 \Err@{\noexpand\alph invalid for numbers > 26}\else\char\count@\fi}
\def\Alph#1{\count@#1\relax\advance\count@64 \ifnum\count@>90
 \Err@{\noexpand\Alph invalid for numbers > 26}\else\char\count@\fi}

\def\Roman#1{\uppercase\expandafter{\romannumeral#1}}
\def\fnsymbol#1{\count@#1\relax
 \count@@\count@
 \advance\count@\m@ne\divide\count@7
 \count@@@\count@\advance\count@@@\@ne
 \multiply\count@7 \advance\count@@-\count@
 \count@\count@@@
 {\loop
  \ifcase\count@@\or*\or\dag\or\ddag\or\P\or\S\or\text{$\|$}\or\#\fi
  \advance\count@\m@ne\ifnum\count@>\z@\repeat}}
\def\cardnine@#1{\ifcase#1\or one\or two\or three\or four\or five\or
 six\or seven\or eight\or nine\fi}
\let\alloc@\alloc@@
\newcount\ten@
\ten@10
\def\cardinal#1{\count@#1\relax
 \ifnum\count@>99 \number\count@
 \else
  \ifnum\count@=\z@ zero%
  \else
   \ifnum\count@<\ten@\cardnine@\count@
   \else
    \ifnum\count@<20
     \advance\count@-\ten@
     \ifcase\count@ ten\or eleven\or twelve\or thirteen\or fourteen\or
      fifteen\or sixteen\or seventeen\or eighteen\or nineteen\fi
    \else
     \count@@\count@\count@@@\count@@
     \divide\count@\ten@\multiply\count@\ten@
     \advance\count@@@-\count@\divide\count@\ten@
     \ifcase\count@\or\or twenty\or thirty\or forty\or fifty\or sixty\or
      seventy\or eighty\or ninety\fi
     \ifnum\count@@@=\z@\else-\cardnine@\count@@@\fi
    \fi
   \fi
  \fi
 \fi}
\def\ordnine@#1{\ifcase#1\or first\or second\or third\or fourth\or fifth\or
 sixth\or seventh\or eighth\or ninth\fi}
\newcount\count@@@@
\def\ordsuffix@{\count@@@@\count@
 \divide\count@\ten@
 \count@@@\count@\count@@\count@
 \divide\count@@\ten@\multiply\count@@\ten@
 \advance\count@@@-\count@@
 \ifnum\count@@@=\@ne th%
 \else
  \count@@@\count@@@@
  \count@@\count@@@@
  \divide\count@@\ten@\multiply\count@@\ten@
  \advance\count@@@-\count@@
  \ifcase\count@@@ th\or st\or nd\or rd\else th\fi
 \fi}
\def\nordinal#1{\count@#1\relax\number\count@\ordsuffix@}
\def\spordinal#1{\count@#1\relax\number\count@$^{\text{\ordsuffix@}}$}
\def\ordinal#1{\count@#1\relax
 \ifnum\count@>99 \number\count@\ordsuffix@
 \else
   \ifnum\count@=\z@ zeroth%
  \else
    \ifnum\count@<\ten@\ordnine@\count@
    \else
     \ifnum\count@<20 \advance\count@-\ten@
      \ifcase\count@ tenth\or eleventh\or twelfth\or thirteenth\or
       fourteenth\or fifteenth\or sixteenth\or seventeenth\or eighteenth\or
       nineteenth\fi
     \else
      \count@@\count@
      \divide\count@\ten@\multiply\count@\ten@
      \count@@@\count@@\advance\count@@@-\count@
      \divide\count@\ten@
      \ifcase\count@\or\or twent\or thirt\or fort\or fift\or sixt\or sevent\or
       eight\or ninet\fi
      \ifnum\count@@@=\z@ ieth\else y-\ordnine@\count@@@\fi
     \fi
    \fi
  \fi
 \fi}
\font@\tensmc=cmcsc10
\textonlyfont@\smc\tensmc
\newtoks\noexpandtoks@
\noexpandtoks@{\let\arabic\relax\let\alph\relax\let\Alph\relax
 \let\Roman\relax\let\fnsymbol\relax\let\rm\relax
 \let\it\relax\let\bf\relax\let\sl\relax\let\smc\relax
 \let\/\relax\let\null\relax}
\def\noexpands@{\the\noexpandtoks@}
\def\Nonexpanding#1{\global\noexpandtoks@
 \expandafter{\the\noexpandtoks@\let#1\relax}}
\def\prevanish@{\saveskip@\z@\ifhmode\saveskip@\lastskip\unskip\fi}
\def\postvanish@{\ifdim\saveskip@>\z@\hskip\saveskip@\fi\FN@\postvanish@@}
\def\postvanish@@{\DN@.{}%
 \ifx\next\space@\ifdim\saveskip@>\z@\DN@. {}\fi\fi\next@.}
\def\invisible#1{\prevanish@\ignorespaces#1\unskip\postvanish@}
\def\vanishlist@{\\\invisible}
\let\noindent@\noindent
\def\noindent{\par\noindent@\FN@\pretendspace@}
\def\pretendspace@{\ismember@\vanishlist@\next
 \iftest@\nobreak\hskip-\p@\hskip\p@\fi}

\newtoks\everypartoks@
\def\noindent@@{\par\everypartoks@\expandafter{\the\everypar}\everypar{}%
 \noindent@\everypar\expandafter{\the\everypartoks@}}
\def\page{\Err@{\noexpand\page has no meaning by itself}}
\let\page@C\pageno
\let\page@P\empty
\let\page@Q\empty
\def\page@S#1{#1\/}
\def\page@F{\rm}
\def\page@N{\arabic}   
\newif\ifindexing@
\def\indexfile{\ifindexing@\else
 \alloc@@7\write\chardef\sixt@@n\ndx@
 \immediate\openout\ndx@=\jobname.ndx
 \global\indexing@true\fi}
\global\advance\insc@unt\m@ne
\ch@ck0\insc@unt\count
\ch@ck1\insc@unt\dimen
\ch@ck2\insc@unt\skip
\ch@ck4\insc@unt\box
\allocationnumber\insc@unt
\global\chardef\margin@\allocationnumber
\dimen\margin@\maxdimen
\count\margin@\z@
\skip\margin@\z@
\newif\ifindexproofing@
\def\indexproofing{\indexproofing@true}
\def\noindexproofing{\indexproofing@false}
\def\unmacro@#1:#2->#3\unmacro@{\def\macpar@{#2}\def\macdef@{#3}}
\def\starparts@#1{\def\stari@{#1}\def\starii@{#1}\let\stariii@\empty
 \test@false
 \DN@##1*##2##3\next@{\ifx\starparts@##2\test@false\else\test@true\fi}%
 \next@#1*\starparts@\next@
 \iftest@\DN@{\starparts@@#1\starparts@@}\else\let\next@\relax\fi\next@}
\def\starparts@@#1*#2\starparts@@{\def\starii@{#1}\def\stariii@{*#2}}
\def\windex@{\ifindexing@
 \expandafter\unmacro@\meaning\stari@\unmacro@
 \edef\macdef@{\string"\macdef@\string"}%
 \edef\next@{\write\ndx@{\macdef@}}\next@
 \write\ndx@{{\number\pageno}{\page@N}{\page@P}{\page@Q}}%
 \fi
 \ifindexproofing@
  \ifx\stariii@\empty\else
   \expandafter\unmacro@\meaning\stariii@\unmacro@\fi
  \insert\margin@{\hbox{\rm\vrule\height9\p@\depth2\p@\width\z@\starii@
  \ifx\stariii@\empty\else\tt\macdef@\fi}}\fi}
\catcode`\"=\active
\def"{\FN@\quote@}
\def\quote@{\ifx\next"\expandafter\quote@@\else\expandafter\quote@@@\fi}
\def\quote@@@#1"{\starparts@{#1}\starii@\windex@}
\def\quote@@"#1"{\prevanish@\starparts@{#1}\windex@\FN@\quote@@@@}
\def\quote@@@@{\ifx\next"\DN@"{\postvanish@}\else
 \let\next@\postvanish@\fi\next@}
\rightadd@"\to\vanishlist@
\def\idefine#1{\DN@{#1}\DNii@{\noexpand#1}%
 \afterassignment\idefine@\def\nextiii@}
\def\idefine@{\ifindexing@
 \expandafter\let\next@\nextiii@
 \expandafter\unmacro@\meaning\nextiii@\unmacro@
 \immediate\write\ndx@{\noexpand\define\nextii@\macpar@{\macdef@}}\fi}
\def\iabbrev*#1#2{\ifindexing@\toks@{#2}%
 \immediate\write\ndx@{\noexpand\abbrev*\noexpand#1{\the\toks@}}\fi}
\newread\laxread@
\newwrite\laxwrite@
\let\fnpages@\empty
\def\Finit@#1#2\Finit@{\let\nextii@#1\def\nextiii@{#2}}
\catcode`\~=11
\def\getparts@ @#1~#2~#3~#4~#5~#6{\def\nextiv@{#1}%
 \def\nextiii@{#2~#3~#4~#5~}\count@#6\relax}
\newif\ifdocument@
\def\document{\ifdocument@\else\global\document@true
 \let\fontlist@\empty
 \immediate\openin\laxread@=\jobname.lax\relax
 {\endlinechar\m@ne\noexpands@\catcode`\@=11 \catcode`\~=11
  \loop\ifeof\laxread@\else
   \read\laxread@ to\next@
   \ifx\next@\empty
   \else
    \expandafter\Finit@\next@\Finit@
    \if\nextii@ F%
     \expandafter\rightadd@\nextiii@\to\fnpages@
    \else
     \expandafter\getparts@\next@
     \edef\next@{\gdef\csname\nextiv@ @L\endcsname{\nextiii@\number\count@}}%
     \next@
    \fi
   \fi
  \repeat}%
 \immediate\closein\laxread@
 \immediate\openout\laxwrite@=\jobname.lax\relax\fi}
\let\thelabel@\relax
\def\thelabels@{\thelabel@ ~\thelabel@@ ~\thelabel@@@ ~\thelabel@@@@ ~}
\def\label#1{\prevanish@
 \ifx\thelabel@\relax
  \Err@{There's nothing here to be labelled}%
 \else
  {\noexpands@
  \expandafter\ifx\csname#1@L\endcsname\relax
   \expandafter\xdef\csname#1@L\endcsname{\thelabels@0}%
   \immediate\write\laxwrite@{@#1~\thelabels@1}%
  \else
   \edef\next@{@~\csname#1@L\endcsname}%
    \expandafter\getparts@\next@
    \ifodd\count@
    \expandafter\xdef\csname#1@L\endcsname{\thelabels@0}%
    \immediate\write\laxwrite@{@#1~\thelabels@1}%
   \else
    \Err@{Label #1 already used}%
   \fi
  \fi
  }%
 \fi
 \postvanish@}
\rightadd@\label\to\vanishlist@
\def\thepages@{\page@N{\number\page@C}~%
 \page@S{\page@P\page@N{\number\page@C}\page@Q}~%
 \number\page@C ~\page@P\page@N{\number\page@C}\page@Q ~}
\def\pagelabel#1{\prevanish@
 \expandafter\ifx\csname#1@L\endcsname\relax
  {\noexpands@
  \expandafter\xdef\csname#1@L\endcsname{\thepages@2}}%
  \write\laxwrite@{@#1~\thepages@3}%
 \else
  {\noexpands@
  \edef\next@{@~\csname#1@L\endcsname}%
  \expandafter\getparts@\next@
  \ifodd\count@
   \ifnum\count@=\@ne
    \expandafter\xdef\csname#1@L\endcsname{\thelabels@2}%
   \fi
   \write\laxwrite@{@#1~\thepages@3}%
  \else
   \Err@{Label #1 already used}%
  \fi
  }%
 \fi
 \postvanish@}
\rightadd@\pagelabel\to\vanishlist@
\newif\ifreferr@
\referr@true
\def\RefErrors{\global\referr@true}
\def\RefWarnings{\global\referr@false}
\setbox\z@\hbox{\global\count@=`^^30}
\ifnum\count@=48 \let\versionthree@\relax\fi
\def\nolabel@#1#2#3{\expandafter\ifx\csname#2@L\endcsname\relax
 \ifreferr@\Err@{No \noexpand\label found for #2}\else
 \W@{Warning: No \noexpand\label found for #2.}%
 \ifx\versionthree@\relax\W@{l.\number\inputlineno\space ... \string#1{#2}}\fi
 \fi#3\else}
\def\csL@#1{{\noexpands@\xdef\Next@{\csname#1@L\endcsname}}}
\def\ref#1{\nolabel@\ref{#1}\relax
 \DNii@##1~##2\nextii@{##1}%
 \csL@{#1}\expandafter\nextii@\Next@\nextii@\fi}
\def\Ref#1{\nolabel@\Ref{#1}\relax
 \DNii@##1~##2~##3\nextii@{##2}%
 \csL@{#1}\expandafter\nextii@\Next@\nextii@\fi}
\def\nref#1{\nolabel@\nref{#1}\relax
 \DNii@##1~##2~##3~##4\nextii@{##3}%
 \csL@{#1}\expandafter\nextii@\Next@\nextii@\fi}
\def\pref#1{\nolabel@\pref{#1}\relax
 \DNii@##1~##2~##3~##4~##5\nextii@{##4}%
 \csL@{#1}\expandafter\nextii@\Next@\nextii@\fi}
\let\pref@\pref
\def\Evaluatenref#1{\nolabel@\Evaluatenref{#1}{\gdef\Nref{-10000 }}%
 \DNii@##1~##2~##3~##4\nextii@{\DNii@{##3}}%
 \csL@{#1}\expandafter\nextii@\Next@\nextii@
 \xdef\Nref{\nextii@}\fi}
\def\Evaluatepref#1{\nolabel@\Evaluatepref{#1}{\global\let\Pref\empty}%
 \DNii@##1~##2~##3~##4~##5\nextii@{\DNii@{##4}}%
 \csL@{#1}\expandafter\nextii@\Next@\nextii@
 \xdef\Pref{\nextii@}\fi}
\def\readlax#1{\immediate\openin\laxread@=#1.lax\relax
 \ifeof\laxread@\W@{}\W@{File #1.lax not found.}\W@{}\fi
 {\endlinechar\m@ne\noexpands@\catcode`\@=11 \catcode`\~=11
  \loop\ifeof\laxread@\else
   \read\laxread@ to\nextv@
   \ifx\nextv@\empty
   \else
    \expandafter\Finit@\nextv@\Finit@
    \ifx\nextii@ F%
    \else
     \expandafter\getparts@\nextv@
     \expandafter\ifx\csname\nextiv@ @L\endcsname\relax
      \edef\next@{\gdef\csname\nextiv@ @L\endcsname
       {\nextiii@\ifnum\count@=\@ne0\else2\fi}}%
      \next@
     \else
      \Err@{Label \nextiv@\space in #1.lax already used}%
     \fi
    \fi
   \fi
  \repeat}%
 \immediate\closein\laxread@}
\catcode`\~=\active

\def\FNSSP@{\FNSS@\pretendspace@}
\everydisplay{\csname displaymath \endcsname}
\expandafter\def\csname displaymath \endcsname#1$${#1$$\FNSSP@}
\def\locallabel@{\let\thelabel@\Thelabel@\let\thelabel@@\Thelabel@@
 \let\thelabel@@@\Thelabel@@@\let\thelabel@@@@\Thelabel@@@@}
\newcount\tag@C
\tag@C\z@
\let\tag@P\empty
\let\tag@Q\empty
\def\tag@S#1{{\rm(}{#1\/}{\rm)}}
\let\tag@N\arabic
\def\tag@F{\rm}
\def\maketag@{\FN@\maketag@@}
\def\maketag@@{\ifx\next\relax\DN@\relax{\FN@\maketag@@}\else
 \ifx\next"\let\next@\maketag@@@\else
 \let\next@\maketag@@@@\fi\fi\next@}
\def\xdefThelabel@#1{\xdef\Thelabel@{#1{\Thelabel@@@}}}
\def\xdefThelabel@@#1{\xdef\Thelabel@@{#1{\Thelabel@@@@}}}
\def\maketag@@@@#1\maketag@{\global\advance\tag@C\@ne
 {\noexpands@
  \xdef\Thelabel@@@{\number\tag@C}%
  \xdefThelabel@\tag@N
  \xdef\Thelabel@@@@{\ifmathtags@$\tag@P\Thelabel@\tag@Q$\else
   \tag@P\Thelabel@\tag@Q\fi}%
  \xdefThelabel@@\tag@S
  }%
 \locallabel@
 \hbox{\tag@F\thelabel@@}%
 #1}
\def\Qlabel@#1{{\noexpands@\xdef\Thelabel@@{#1}%
 \let\style\empty\xdef\Thelabel@@@@{#1}%
 \let\pre\empty\let\post\empty\xdef\Thelabel@{#1}%
 \let\numstyle\empty\xdef\Thelabel@@@{#1}}}
\def\maketag@@@"#1"#2\maketag@{%
 {\let\pre\tag@P\let\post\tag@Q\let\style\tag@S\let\numstyle\tag@N
  \hbox{\tag@F#1}%
  \noexpands@
  \Qlabel@{#1}%
  }%
 \locallabel@
 #2}
\def\align@{\inalign@true\inany@true
 \vspace@\allowdisplaybreak@\displaybreak@\intertext@
 \def\tag{\global\tag@true\ifnum\and@=\z@
  \DN@{&\omit\global\rwidth@\z@&\relax}\else
  \DN@{&\relax}\fi\next@}%
 \iftagsleft@\DN@{\csname align \endcsname}\else
  \DN@{\csname align \space\endcsname}\fi\next@}
\def\noset@{\def\Offset##1##2{\prevanish@\postvanish@}%
 \def\Reset##1##2{\prevanish@\postvanish@}}
\def\measure@#1\endalign{\global\lwidth@\z@\global\rwidth@\z@
 \global\maxlwidth@\z@\global\maxrwidth@\z@
 \global\and@\z@
 \setbox\z@\vbox
  {\noset@\everycr{\noalign{\global\tag@false\global\and@\z@}}\Let@
  \halign{\setboxz@h{$\m@th\displaystyle{\@lign##}$}%
   \global\lwidth@\wdz@
   \ifdim\lwidth@>\maxlwidth@\global\maxlwidth@\lwidth@\fi
   \global\advance\and@\@ne
   &\setboxz@h{$\m@th\displaystyle{{}\@lign##}$}\global\rwidth@\wdz@
   \ifdim\rwidth@>\maxrwidth@\global\maxrwidth@\rwidth@\fi
   \global\advance\and@\@ne
   &\Tag@\eat@{##}\crcr#1\crcr}}%
 \totwidth@\maxlwidth@\advance\totwidth@\maxrwidth@}
\def\prepost@{\global\let\tag@P@\tag@P\global\let\tag@Q@\tag@Q}
\def\reprepost@{\let\tag@P\tag@P@\let\tag@Q\tag@Q@}
\expandafter\def\csname align \space\endcsname#1\endalign
 {\measure@#1\endalign\global\and@\z@
 \ifingather@\everycr{\noalign{\global\and@\z@}}\else\displ@y@\fi
 \Let@\tabskip\centering@
 \halign to\displaywidth
  {\hfil\strut@\setboxz@h{$\m@th\displaystyle{\@lign##\prepost@}$}%
  \boxz@\global\advance\and@\@ne
  \tabskip\z@skip
  &\setboxz@h{$\m@th\displaystyle{{}\@lign##\prepost@}$}%
  \global\rwidth@\wdz@\boxz@\hfil\global\advance\and@\@ne
  \tabskip\centering@
  &\setboxz@h{\@lign\strut@\reprepost@\maketag@##\maketag@}%
  \dimen@\displaywidth\advance\dimen@-\totwidth@
  \divide\dimen@\tw@\advance\dimen@\maxrwidth@\advance\dimen@-\rwidth@
  \ifdim\dimen@<\tw@\wdz@\llap{\vtop{\normalbaselines\null\boxz@}}%
  \else\llap{\boxz@}\fi
  \tabskip\z@skip
  \crcr#1\crcr
  \black@\totwidth@}}
\expandafter\def\csname align \endcsname#1\endalign{\measure@#1\endalign
 \global\and@\z@
 \ifdim\totwidth@>\displaywidth\let\displaywidth@\totwidth@\else
  \let\displaywidth@\displaywidth\fi
 \ifingather@\everycr{\noalign{\global\and@\z@}}\else\displ@y@\fi
 \Let@\tabskip\centering@\halign to\displaywidth
  {\hfil\strut@\setboxz@h{$\m@th\displaystyle{\@lign##\prepost@}$}%
  \global\lwidth@\wdz@\global\lineht@\ht\z@
  \boxz@\global\advance\and@\@ne
  \tabskip\z@skip&\setboxz@h{$\m@th\displaystyle{{}\@lign##\prepost@}$}%
  \ifdim\ht\z@>\lineht@\global\lineht@\ht\z@\fi
  \boxz@\hfil\global\advance\and@\@ne
  \tabskip\centering@&\kern-\displaywidth@
  \setboxz@h{\@lign\strut@\reprepost@\maketag@##\maketag@}%
  \dimen@\displaywidth\advance\dimen@-\totwidth@
  \divide\dimen@\tw@\advance\dimen@\maxlwidth@\advance\dimen@-\lwidth@
  \ifdim\dimen@<\tw@\wdz@
   \rlap{\vbox{\normalbaselines\boxz@\vbox to\lineht@{}}}\else
   \rlap{\boxz@}\fi
  \tabskip\displaywidth@\crcr#1\crcr\black@\totwidth@}}
\def\attag@#1{\let\Maketag@\maketag@\let\TAG@\Tag@
 \let\Prepost@\prepost@\let\Reprepost@\reprepost@
 \let\Tag@\relax\let\maketag@\relax
 \let\prepost@\relax\let\reprepost@\relax
 \ifmeasuring@
  \def\llap@##1{\setboxz@h{##1}\hbox to\tw@\wdz@{}}%
  \def\rlap@##1{\setboxz@h{##1}\hbox to\tw@\wdz@{}}%
 \else\let\llap@\llap\let\rlap@\rlap\fi
 \toks@{\hfil\strut@
  $\m@th\displaystyle{\@lign\the\hashtoks@\prepost@}$%
  \tabskip\z@skip\global\advance\and@\@ne&
  $\m@th\displaystyle{{}\@lign\the\hashtoks@\prepost@}$\hfil
  \ifxat@\tabskip\centering@\fi\global\advance\and@\@ne}%
 \iftagsleft@
  \toks@@{\tabskip\centering@&\Tag@\kern-\displaywidth
   \rlap@{\@lign\reprepost@\maketag@\the\hashtoks@\maketag@}%
   \global\advance\and@\@ne\tabskip\displaywidth}\else
  \toks@@{\tabskip\centering@&\Tag@\llap@{\@lign\reprepost@\maketag@
   \the\hashtoks@\maketag@}\global\advance\and@\@ne\tabskip\z@skip}\fi
 \atcount@#1\relax\advance\atcount@\m@ne
 \loop\ifnum\atcount@>\z@
  \toks@\expandafter{\the\toks@&\hfil$\m@th\displaystyle{\@lign
  \the\hashtoks@\prepost@}$\global\advance\and@\@ne
  \tabskip\z@skip
  &$\m@th\displaystyle{{}\@lign\the\hashtoks@\prepost@}$\hfil\ifxat@
  \tabskip\centering@\fi\global\advance\and@\@ne}\advance\atcount@\m@ne
 \repeat
 \edef\preamble@{\the\toks@\the\toks@@}%
 \edef\preamble@@{\preamble@}%
 \let\maketag@\Maketag@\let\Tag@\TAG@
 \let\prepost@\Prepost@\let\reprepost@\Reprepost@}
\def\unlabel@{\def\label##1{\prevanish@\postvanish@}%
 \def\pagelabel##1{\prevanish@\postvanish@}}
\newcount\tag@CC
\expandafter\def\csname alignat \endcsname#1#2\endalignat
 {\inany@true\xat@false
 \def\tag{\global\tag@true
  \count@#1\relax\multiply\count@\tw@\advance\count@\m@ne
  \gdef\tag@{&}%
  \loop\ifnum\count@>\and@\xdef\tag@{&\omit\tag@}%
  \advance\count@\m@ne\repeat
  \tag@\relax}%
 \vspace@\allowdisplaybreak@\displaybreak@\intertext@
 \displ@y@\measuring@true\tag@CC\tag@C
 \setbox\savealignat@\hbox{\noset@\unlabel@$\m@th\displaystyle\Let@
  \attag@{#1}\vbox{\halign{\span\preamble@@\crcr#2\crcr}}$}%
 \measuring@false
 \Let@\attag@{#1}\tag@C\tag@CC
 \tabskip\centering@\halign to\displaywidth
  {\span\preamble@@\crcr#2\crcr\black@{\wd\savealignat@}}}
\expandafter\def\csname xalignat \endcsname#1#2\endxalignat
 {\inany@true\xat@true
 \def\tag{\global\tag@true
  \count@#1\relax\multiply\count@\tw@\advance\count@\m@ne
  \gdef\tag@{&}%
  \loop\ifnum\count@>\and@\xdef\tag@{&\omit\tag@}%
  \advance\count@\m@ne\repeat
  \tag@\relax}%
 \vspace@\allowdisplaybreak@\displaybreak@\intertext@
 \displ@y@\measuring@true\tag@CC\tag@C
 \setbox\savealignat@\hbox{\noset@\unlabel@$\m@th\displaystyle\Let@
  \attag@{#1}\vbox{\halign{\span\preamble@@\crcr#2\crcr}}$}%
 \measuring@false\Let@\attag@{#1}\tag@C\tag@CC
 \tabskip\centering@\halign to\displaywidth
 {\span\preamble@@\crcr#2\crcr\black@{\wd\savealignat@}}}
\def\gather{\RIfMIfI@\DN@{\onlydmatherr@\gather}\else
 \ingather@true\inany@true\def\tag{&\relax}%
 \vspace@\allowdisplaybreak@\displaybreak@\intertext@
 \displ@y\Let@
 \iftagsleft@\DN@{\csname gather \endcsname}\else
  \DN@{\csname gather \space\endcsname}\fi\fi
 \else\DN@{\onlydmatherr@\gather}\fi\next@}
\def\exstring@{\expandafter\eat@\string}
\def\newcounter#1{\define#1{}%
 \edef\next@{\def\noexpand#1{\futurelet\noexpand\next
  \csname\exstring@#1@Z\endcsname}}\next@
 \edef\next@{\def\csname\exstring@#1@Z\endcsname
  {\global\advance\csname\exstring@#1@C\endcsname\@ne
  {\csname\exstring@#1@F\endcsname\csname\exstring@#1@S\endcsname
   {\csname\exstring@#1@P\endcsname\csname\exstring@#1@N\endcsname
   {\noexpand\number\csname\exstring@#1@C\endcsname}%
   \csname\exstring@#1@Q\endcsname}}%
  \noexpand\ifx\noexpand\next\noexpand\label
   \def\noexpand\next@\noexpand\label########1{{\noexpand\noexpands@
    \xdef\noexpand\Thelabel@{\csname\exstring@#1@N\endcsname
     {\noexpand\number\csname\exstring@#1@C\endcsname}}%
    \xdef\noexpand\Thelabel@@@{\noexpand\number
     \csname\exstring@#1@C\endcsname}%
    \xdef\noexpand\Thelabel@@{\csname\exstring@#1@S\endcsname
     {\csname\exstring@#1@P\endcsname
     \csname\exstring@#1@N\endcsname
     {\noexpand\number\csname\exstring@#1@C\endcsname}%
     \csname\exstring@#1@Q\endcsname}}%
    \xdef\noexpand\Thelabel@@@@{\csname\exstring@#1@P\endcsname
     \csname\exstring@#1@N\endcsname
     {\noexpand\number\csname\exstring@#1@C\endcsname}%
     \csname\exstring@#1@Q\endcsname}}%
    {\noexpand\locallabel@\noexpand\label{########1}}}%
   \noexpand\else\let\noexpand\next@\relax\noexpand\fi\noexpand\next@}}\next@
 \expandafter\newcount@\csname\exstring@#1@C\endcsname
 \expandafter\let\csname\exstring@#1@N\endcsname\arabic
 \expandafter\def\csname\exstring@#1@S\endcsname##1{##1\/}%
 \expandafter\let\csname\exstring@#1@P\endcsname\empty
 \expandafter\let\csname\exstring@#1@Q\endcsname\empty
 \expandafter\def\csname\exstring@#1@F\endcsname{\rm}%
 }
\def\HASH@#1#2{\ifnum#2=\z@\else
 \edef\next@{\toks@{\the\toks@\the\hashtoks@#2}%
 \toks@@{\the\toks@@{\the\hashtoks@#2}}}\next@\expandafter\HASH@\fi}
\def\HASH@@{\toks@{}\toks@@{}\expandafter\HASH@\macpar@00}
\def\usecounter#1#2{\expandafter\ifx\csname\exstring@#1@Z\endcsname
 \relax\Err@{\noexpand#1not created with \string\newcounter}\fi
 \expandafter\let\csname\exstring@#1@@Z\endcsname\relax
 \expandafter\let\csname\exstring@#1@@Z@\endcsname\relax
 \expandafter\let\csname\exstring@#1@@Z@@\endcsname\relax
 \edef\next@{\def\noexpand#2{\futurelet\noexpand\next
  \csname\exstring@#1@@Z\endcsname}}\next@
 \edef\next@{\def\csname\exstring@#1@@Z\endcsname{\noexpand\ifx
  \noexpand\next\noexpand\label\def\noexpand\next@\noexpand\label
   ########1{\csname\exstring@#1@@Z@\endcsname
   {\noexpand#1\noexpand\label{########1}}}%
   \noexpand\else\noexpand\ifx\noexpand\next
   \noexpand"\def\noexpand\next@\noexpand"########1\noexpand"%
   {\csname\exstring@#1@@Z@\endcsname{{\expandafter\noexpand
   \csname\exstring@#1@F\endcsname
   \let\noexpand\pre\expandafter\noexpand\csname\exstring@#1@P\endcsname
   \let\noexpand\post\expandafter\noexpand\csname\exstring@#1@Q\endcsname
   \let\noexpand\style\expandafter\noexpand\csname\exstring@#1@S\endcsname
   \let\noexpand\numstyle\expandafter\noexpand\csname\exstring@#1@N\endcsname
   ########1}}}\noexpand\else
   \def\noexpand\next@{\csname\exstring@#1@@Z@\endcsname{\noexpand#1}}%
   \noexpand\fi\noexpand\fi\noexpand\next@}}\next@
 \def\next@{\expandafter\expandafter\expandafter\unmacro@\expandafter
  \meaning\csname\exstring@#1@@Z@@\endcsname\unmacro@
  \HASH@@
  \edef\next@{\def\csname\exstring@#1@@Z@\endcsname\the\toks@{%
   \expandafter\noexpand\csname\exstring@#1@@Z@@\endcsname\the\toks@@
   \noexpand\FNSSP@}}\next@}%
 \afterassignment\next@
 \expandafter\def\csname\exstring@#1@@Z@@\endcsname}
\def\listbi@{\penalty50 \medskip}
\def\listbii@{\penalty100 \smallskip}
\let\listbiii@\relax
\let\listbiv@\relax
\let\listbv@\relax
\def\listmi@{\advance\leftskip30\p@\relax}
\let\listmii@\listmi@
\let\listmiii@\listmi@
\let\listmiv@\listmi@
\let\listmv@\listmi@
\def\itemi@#1{\noindent@@\llap{#1\hskip5\p@}}
\let\itemii@\itemi@
\let\itemiii@\itemi@
\let\itemiv@\itemi@
\let\itemv@\itemi@
\def\liste@{\penalty-50 \medskip}
\def\listei@{\penalty-100 \smallskip}
\let\listeii@\relax
\let\listeiii@\relax
\let\listeiv@\relax
\expandafter\newcount\csname list@C1\endcsname
\csname list@C1\endcsname\z@
\expandafter\newcount\csname list@C2\endcsname
\csname list@C2\endcsname\z@
\expandafter\newcount\csname list@C3\endcsname
\csname list@C3\endcsname\z@
\expandafter\newcount\csname list@C4\endcsname
\csname list@C4\endcsname\z@
\expandafter\newcount\csname list@C5\endcsname
\csname list@C5\endcsname\z@
\expandafter\let\csname list@P1\endcsname\empty
\expandafter\let\csname list@P2\endcsname\empty
\expandafter\let\csname list@P3\endcsname\empty
\expandafter\let\csname list@P4\endcsname\empty
\expandafter\let\csname list@P5\endcsname\empty
\expandafter\let\csname list@Q1\endcsname\empty
\expandafter\let\csname list@Q2\endcsname\empty
\expandafter\let\csname list@Q3\endcsname\empty
\expandafter\let\csname list@Q4\endcsname\empty
\expandafter\let\csname list@Q5\endcsname\empty
\expandafter\def\csname list@S1\endcsname#1{{\rm(}{#1\/}{\rm)}}
\expandafter\def\csname list@S2\endcsname#1{{\rm(}{#1\/}{\rm)}}
\expandafter\def\csname list@S3\endcsname#1{{\rm(}{#1\/}{\rm)}}
\expandafter\def\csname list@S4\endcsname#1{{\rm(}{#1\/}{\rm)}}
\expandafter\def\csname list@S5\endcsname#1{{\rm(}{#1\/}{\rm)}}
\expandafter\let\csname list@N1\endcsname\arabic
\expandafter\let\csname list@N2\endcsname\arabic
\expandafter\let\csname list@N3\endcsname\arabic
\expandafter\let\csname list@N4\endcsname\arabic
\expandafter\let\csname list@N5\endcsname\arabic
\expandafter\def\csname list@F1\endcsname{\rm}
\expandafter\def\csname list@F2\endcsname{\rm}
\expandafter\def\csname list@F3\endcsname{\rm}
\expandafter\def\csname list@F4\endcsname{\rm}
\expandafter\def\csname list@F5\endcsname{\rm}
\newcount\listlevel@
\listlevel@\z@
\def\list@@C{\csname list@C\number\listlevel@\endcsname}
\def\list@@P{\csname list@P\number\listlevel@\endcsname}
\def\list@@Q{\csname list@Q\number\listlevel@\endcsname}
\def\list@@S{\csname list@S\number\listlevel@\endcsname}
\def\list@@N{\csname list@N\number\listlevel@\endcsname}
\def\list@@F{\csname list@F\number\listlevel@\endcsname}
\newif\iffirstitemi@
\newif\iffirstitemii@
\newif\iffirstitemiii@
\newif\iffirstitemiv@
\newif\iffirstitemv@
\def\Firstitem@true{\csname firstitem\romannumeral\listlevel@
 @true\endcsname}
\def\Firstitem@false{\csname firstitem\romannumeral\listlevel@
 @false\endcsname}
\def\Listm@{\csname listm\romannumeral\listlevel@ @\endcsname}
\def\Item@{\csname item\romannumeral\listlevel@ @\endcsname}
\def\Liste@{\csname liste\romannumeral\listlevel@ @\endcsname}
\newif\iflistcontinue@
\def\keepitem{\listcontinue@true}
\newcount\list@C@
\def\list{%
 \iflistcontinue@\csname list@C1\endcsname\csname list@C@\endcsname\fi
 \global\csname list@C2\endcsname\z@
 \global\csname list@C3\endcsname\z@
 \global\csname list@C4\endcsname\z@
 \global\csname list@C5\endcsname\z@
 \begingroup
 \firstitemi@true
 \listlevel@\@ne
 \def\item{\FN@\item@}%
 \FN@\list@}
\Invalid@\runinitem
\def\list@{\ifx\next\par
 \DN@\par{\FN@\list@}\else
 \ifx\next\runinitem
  \DN@\runinitem{\FN@\runinitem@}\else
  \DN@{\par\dimen@\parskip\parskip\dimen@}\fi\fi\next@}
\newif\ifoutlevel@
\newif\ifrunin@
\def\item@{%
 \ifoutlevel@\Liste@\outlevel@false\fi
 \ifrunin@\runin@false\par
  \dimen@\parskip\parskip\dimen@
  \Listm@\fi
 \iffirstitemi@\listbi@\listmi@\firstitemi@false\else\par\fi
 \iffirstitemii@\listbii@\listmii@\firstitemii@false\else\par\fi
 \iffirstitemiii@\listbiii@\listmiii@\firstitemiii@false\else\par\fi
 \iffirstitemiv@\listbiv@\listmiv@\firstitemiv@false\else\par\fi
 \iffirstitemv@\listbv@\listmv@\firstitemv@false\else\par\fi
 \DN@"##1"{{\let\pre\list@@P\let\post\list@@Q
  \let\style\list@@S\let\numstyle\list@@N
  \vskip-\parskip
  \Item@{\list@@F##1}%
  \noexpands@
  \Qlabel@{##1}}%
  \locallabel@
  \FNSSP@}%
 \DNii@{\global\advance\list@@C\@ne
  {\noexpands@
   \xdef\Thelabel@@@{\number\list@@C}%
   \xdefThelabel@\list@@N
   \xdef\Thelabel@@@@{\list@@P\Thelabel@\list@@Q}%
   \xdefThelabel@@\list@@S
  }%
  \locallabel@
  \vskip-\parskip
  \Item@{\list@@F\thelabel@@}%
  \FN@\pretendspace@}%
 \ifx\next"\expandafter\next@\else\expandafter\nextii@\fi}
\def\runinitem@{%
  \runin@true
  \Firstitem@false
  \DN@"##1"{{\let\pre\list@@P\let\post\list@@Q
   \let\style\list@@S\let\numstyle\list@@N
   \unskip\space{\list@@F##1} %
   \noexpands@
   \Qlabel@{##1}}%
   \locallabel@
   \ignorespaces}%
  \DNii@{\global\advance\list@@C\@ne
   {\noexpands@
    \xdef\Thelabel@@@{\number\list@@C}%
    \xdefThelabel@\list@@N
    \xdef\Thelabel@@@@{\list@@P\Thelabel@\list@@Q}%
    \xdefThelabel@@\list@@S
   }%
   \locallabel@
   \unskip\space{\list@@F\thelabel@@} }%
  \ifx\next"\expandafter\next@\else\expandafter\nextii@\fi}
\def\inlevel{\ifnum\listlevel@=5
 \DN@{\Err@{Already 5 levels down}}\else
 \DN@{\begingroup\advance\listlevel@\@ne
 \Firstitem@true\FN@\inlevel@}\fi\next@}
\def\inlevel@{\ifx\next\par
 \DN@\par{\FN@\inlevel@}\else
 \ifx\next\runinitem
  \DN@\runinitem{\FN@\runinitem@}\else
  \let\next@\relax\fi\fi\next@}
\def\outlevel{\ifnum\listlevel@=\@ne
 \Err@{At top level}\else
 \par\global\list@@C\z@\endgroup\outlevel@true\fi}
\def\endlist{%
 \expandafter\global\csname list@C@\endcsname\csname list@C1\endcsname
 \par
 \global\toks\@ne{}\count@\listlevel@
 {\loop
  \ifnum\count@>\z@\global\toks\@ne\expandafter{\the\toks\@ne\endgroup}%
  \advance\count@\m@ne
  \repeat}%
 \the\toks\@ne
 \liste@
 \listcontinue@false\global\csname list@C1\endcsname\z@
 \vskip-\parskip
 \noindent@@
 \FN@\pretendspace@}
\newif\iffirstdescribe@
\def\describe{\par
 \begingroup\firstdescribe@true
 \def\item##1{%
  \iffirstdescribe@\penalty50 \medskip\vskip-\parskip
  \firstdescribe@false\else\par\fi
  \noindent@@\hangindent2pc\hangafter\@ne
  {\bf##1}\hskip.5em}}

\Invalid@\pullin
\Invalid@\pullinmore
\newif\iffirstpull@
\def\margins{\par\begingroup\firstpull@true
 \def\pullin##1##2{\par
  \iffirstpull@\firstpull@false\else\endgroup\fi
  \begingroup\DN@{##1}%
  \ifx\next@\empty\leftskip\z@\else\ifx\next@\space\leftskip\z@
  \else\leftskip##1\fi\fi
  \DN@{##2}\ifx\next@\empty\rightskip\z@\else\ifx\next@\space
  \rightskip\z@\else\rightskip##2\fi\fi\ignorespaces}%
 \def\pullinmore##1##2{\par
  \xdef\Next@{\leftskip\the\leftskip\relax\rightskip\the\rightskip\relax}%
  \iffirstpull@\firstpull@false\else\endgroup\fi
  \begingroup\Next@
  \DN@{##1}%
  \ifx\next@\empty\else\ifx\next@\space\else\advance\leftskip##1\fi\fi
  \DN@{##2}\ifx\next@\empty\else\ifx\next@\space\else
  \advance\rightskip##2\fi\fi\ignorespaces}}

\newif\ifnopunct@
\newif\ifnospace@
\newif\ifoverlong@
\let\nofrillslist@\empty
\let\overlonglist@\empty
\def\nopunct{\nopunct@true\FN@\nopunct@}
\def\nospace{\nospace@true\FN@\nospace@}
\def\overlong{\overlong@true\FN@\overlong@}
\def\nopunct@{\ifx\next\nospace
 \DN@\nospace{\nospace@true\FN@\nopnos@}\else\ifx\next\overlong
 \DN@\overlong{\overlong@true\FN@\nopol@}\else
 \let\next@\nopunct@@\fi\fi\next@}
\def\nopunct@@#1{\ismember@\nofrillslist@#1%
 \iftest@\let\next@#1\else
 \DN@{\nopunct@false\Err@{\noexpand\nopunct can't be used with
 \string#1}#1}\fi\next@}
\def\nospace@{\ifx\next\nopunct
 \DN@\nopunct{\nopunct@true\FN@\nopnos@}\else\ifx\next\overlong
 \DN@\overlong{\overlong@true\FN@\nosol@}\else
 \let\next@\nospace@@\fi\fi\next@}
\def\nospace@@#1{\ismember@\nofrillslist@#1%
 \iftest@\let\next@#1\else
 \DN@{\nospace@false\Err@{\noexpand\nospace can't be used with
 \string#1}#1}\fi\next@}
\def\overlong@{\ifx\next\nopunct
 \DN@\nopunct{\nopunct@true\FN@\nopol@}\else\ifx\next\nospace
 \DN@\nospace{\nospace@true\FN@\nosol@}\else
 \let\next@\overlong@@\fi\fi\next@}
\def\overlong@@#1{\ismember@\overlonglist@#1%
 \iftest@\let\next@#1\else
 \DN@{\overlong@false\Err@{\noexpand\overlong can't be used with
 \string#1}#1}\fi\next@}
\def\nopnos@{\ifx\next\overlong
 \DN@\overlong{\overlong@true\nopnosol@}\else
 \let\next@\nopnos@@\fi\next@}
\def\nopol@{\ifx\next\nospace
 \DN@\nospace{\nospace@true\nopnosol@}\else
 \let\next@\nopol@@\fi\next@}
\def\nosol@{\ifx\next\nopunct
 \DN@\nopunct{\nopunct@true\nopnosol@}\else
 \let\next@\nosol@@\fi\next@}
\def\nopnos@@#1{\ismember@\nofrillslist@#1%
 \iftest@\let\next@#1\else
 \DN@{\nopunct@false\nospace@false
  \Err@{\noexpand\nopunct\noexpand\nospace
   can't be used with \string#1}#1}\fi\next@}
\def\testii@#1{\ismember@\nofrillslist@#1%
 \iftest@\let\nextiii@ T\else\let\nextiii@ F\fi
 \ismember@\overlonglist@#1%
 \iftest@\let\nextiv@ T\else\let\nextiv@ F\fi
 \test@false\if\nextiii@ T\if\nextiv@ T\test@true\fi\fi}
\def\nopol@@#1{\testii@{#1}%
 \iftest@\let\next@#1%
 \else\DN@{\if\nextiii@ T\else\nopunct@false\fi
  \if\nextiv@ T\else\overlong@false\fi
  \Err@{\if\nextiii@ T\else\noexpand\nopunct\fi
  \if\nextiv@ T\else\noexpand\overlong\fi can't be used
  with \string#1}#1}\fi\next@}
\def\nosol@@#1{\testii@{#1}%
 \iftest@\let\next@#1%
 \else\DN@{\if\nextiii@ T\else\nospace@false\fi
  \if\nextiv@ T\else\overlong@false\fi
  \Err@{\if\nextiii@ T\else\noexpand\nospace\fi
  \if\nextiv@ T\else\noexpand\overlong\fi can't be used
  with \string#1}#1}\fi\next@}
\def\nopnosol@#1{\testii@{#1}%
 \iftest@\let\next@#1%
 \else\DN@{\if\nextiii@ T\else\nopunct@false\nospace@false\fi
  \if\nextiv@ T\else\overlong@false\fi
  \Err@{\if\nextiii@ T\else\noexpand\nopunct\noexpand\nospace\fi
  \if\nextiv@ T\else\noexpand\overlong\fi can't be used
  with \string#1}#1}\fi\next@}
\def\punct@#1{\ifnopunct@\else#1\fi}
\def\addspace@#1{\ifnospace@\else#1\fi}
\def\hss@{\ifoverlong@\z@ plus\@m\p@ minus\@m\p@
 \else \z@ plus\@m\p@\fi}
\rightadd@\demo\to\nofrillslist@
\newif\ifclaim@
\def\exxx@{\expandafter\expandafter\expandafter\eat@\expandafter\string}
\let\colon@:
\def\demo#1{\ifclaim@
 \Err@{Previous \expandafter\noexpand\claimtype@ has
  no matching \string\end\exxx@\claimtype@}%
 \let\next@\relax
 \else
  \par
  \ifdim\lastskip<\smallskipamount\removelastskip\smallskip\fi
  \begingroup
  \noindent@@{\smc\ignorespaces#1\unskip
   \punct@{\null\colon@}\addspace@\enspace}%
  \nopunct@false\nospace@false
  \rm
  \DN@{\FNSSP@}%
 \fi
 \next@}
\def\enddemo{\par\endgroup\nopunct@false\nospace@false\smallskip}
\rightadd@\claim\to\nofrillslist@
\def\claim@F{\smc}
\def\claim@@@F{\csname\exxx@\claimtype@ @F\endcsname}
\def\claimformat@#1#2#3{%
 \medbreak\noindent@@{\smc#1 {\claim@@@F#2} #3%
 \punct@{\null.}\addspace@\enspace}\sl}
\def\claimformat@@#1#2{\claimformat@{\ignorespaces#1\unskip}%
 {\ifx\thelabel@@\empty\unskip\else\thelabel@@\fi}%
 {\ignorespaces#2\unskip}%
 \let\Claimformat@@\claimformat@@\FNSSP@}
\let\Claimformat@@\claimformat@@
\def\claim@@@P{\csname\exxx@\claimtype@ @P\endcsname}
\def\claim@@@Q{\csname\exxx@\claimtype@ @Q\endcsname}
\def\claim@@@S{\csname\exxx@\claimtype@ @S\endcsname}
\def\claim@@@N{\csname\exxx@\claimtype@ @N\endcsname}
\def\claim@@@C{\csname claim@C\claimclass@\endcsname}
\newcount\claim@C
\claim@C\z@
\let\claim@P\empty
\let\claim@Q\empty
\def\claim@S#1{#1\/}
\let\claim@N\arabic
\def\claim{\claim@true\let\claimclass@\empty
 \def\claimtype@{\claim}\FN@\claim@}
\def\claim@{%
 \ifx\next\c
  \let\next@\claim@c
 \else
  \ifx\next"%
   \let\next@\claim@q
  \else
   \begingroup\global\advance\claim@C\@ne
   {\noexpands@
    \xdef\Thelabel@@@{\number\claim@C}%
    \xdefThelabel@\claim@N
    \xdef\Thelabel@@@@{\claim@P\Thelabel@\claim@Q}%
    \xdefThelabel@@\claim@S
   }%
   \locallabel@
   \let\next@\Claimformat@@
  \fi
 \fi
 \next@}
\def\claim@c\c#1{\claim@true\begingroup
 \expandafter
 \ifx\csname claim@C#1\endcsname\relax
  \expandafter\newcount@\csname claim@C#1\endcsname
  \global\csname claim@C#1\endcsname\@ne
 \else
  \global\advance\csname claim@C#1\endcsname\@ne
 \fi
 \def\claimclass@{#1}%
 {\noexpands@
  \xdef\Thelabel@@@{\number\claim@@@C}%
  \xdefThelabel@\claim@@@N
  \xdef\Thelabel@@@@{\claim@@@P\Thelabel@\claim@@@Q}%
  \xdefThelabel@@\claim@@@S
 }%
 \locallabel@
 \FNSS@\claim@c@}
\def\claim@q"#1"{\begingroup
 {\let\pre\claim@@@P\let\post\claim@@@Q
  \let\style\claim@@@S\let\numstyle\claim@@@N
  \noexpands@
  \Qlabel@{#1}}%
 \locallabel@
 \FNSS@\claim@q@}
\def\claim@c@{\ifx\next"%
 \global\advance\claim@@@C\m@ne\let\next@\claim@cq
 \else\let\next@\Claimformat@@\fi\next@}
\def\claim@cq"#1"{{\let\pre\claim@@@P\let\post\claim@@@Q
 \let\style\claim@@@S\let\numstyle\claim@@@N
 \noexpands@
 \Qlabel@{#1}}%
 \locallabel@
 \FNSS@\Claimformat@@}
\def\claim@q@{\ifx\next\c\expandafter\claim@qc
 \else\expandafter\Claimformat@@\fi}
\def\claim@qc\c#1{\expandafter\ifx\csname claim@C#1\endcsname\relax
 \expandafter\newcount@\csname claim@C#1\endcsname
 \global\csname claim@C#1\endcsname\z@\fi
 \FNSS@\Claimformat@@}
\def\endclaim{\endgroup\claim@false\nopunct@false\nospace@false
 \let\Claimformat@@\claimformat@@\medbreak}
\Invalid@\claimclause
\def\newclaim{\FN@\newclaim@}
\def\newclaim@{\ifx\next\claimclause
 \DN@\claimclause##1{\newclaim@@{##1}}\else
 \DN@{\newclaim@@\relax}\fi\next@}
\def\claimlist@{\\\claim}
\newtoks\claim@i
\newtoks\claim@v
\let\noclaimclause@=F
\def\newclaim@@#1#2#3\c#4#5{\define#2{}%
 \rightadd@#2\to\claimlist@\rightadd@#2\to\nofrillslist@%
 \expandafter\def\csname\exstring@#2@P\endcsname{\claim@P}%
 \expandafter\def\csname\exstring@#2@Q\endcsname{\claim@Q}%
 \expandafter\def\csname\exstring@#2@S\endcsname{\claim@S}%
 \expandafter\def\csname\exstring@#2@N\endcsname{\claim@N}%
 \expandafter\def\csname\exstring@#2@F\endcsname{\claim@F}%
 \expandafter\def\csname end\exstring@#2\endcsname{\endclaim}%
 \expandafter\ifx\csname claim@C#4\endcsname\relax
  \expandafter\newcount@\csname claim@C#4\endcsname
  \global\csname claim@C#4\endcsname\z@\fi
 \edef\next@{\let\csname\exstring@#2@C\endcsname
   \csname claim@C#4\endcsname}\next@
 \def#2{\ifx\noclaimclause@ T\else#1\fi
  \global\claim@i{#1}\gdef\claim@iv{#4}\global\claim@v{#5}%
  \def\claimtype@{#2}\def\Claimformat@@{\claimformat@@{#5}}\claim@c\c{#4}}}
\def\shortenclaim#1#2{\define#2{}%
 \ismember@\claimlist@#1%
 \iftest@
  \rightadd@#2\to\nofrillslist@%
  \expandafter\def\csname\exstring@#2@P\endcsname
   {\csname\exstring@#1@P\endcsname}%
  \expandafter\def\csname\exstring@#2@Q\endcsname
   {\csname\exstring@#1@Q\endcsname}%
  \expandafter\def\csname\exstring@#2@S\endcsname
   {\csname\exstring@#1@S\endcsname}%
  \expandafter\def\csname\exstring@#2@N\endcsname
   {\csname\exstring@#1@N\endcsname}%
  \expandafter\def\csname\exstring@#2@F\endcsname
   {\csname\exstring@#1@F\endcsname}%
  \expandafter\def\csname end\exstring@#2\endcsname{\endclaim}%
  \edef\next@{\let\csname\exstring@#2@C\endcsname
    \csname claim\exstring@#1C\endcsname}\next@
  \setbox\z@\vbox{\let\noclaimclause@ T#1""\relax\endgroup}%
  \edef#2{\the\claim@i
   \def\noexpand\claimtype@{\noexpand#2}%
   \def\noexpand\Claimformat@@{\noexpand\claimformat@@{\the\claim@v}\relax}%
   \noexpand\claim@c\noexpand\c{\claim@iv}}%
 \else
  \Err@{\noexpand#1not yet created by \string\newclaim}%
 \fi}
\def\classtest@#1{\DN@{#1}\ifx\next@\claimclass@
 \test@true\else\test@false\fi}
\def\typetest@#1{\DN@{#1}\ifx\next@\claimtype@\test@true\else
  \test@false\fi}
\newif\iftoc@
\def\tocfile{\iftoc@\else\alloc@@7\write\chardef\sixt@@n\toc@
 \immediate\openout\toc@=\jobname.toc
 \alloc@@7\write\chardef\sixt@@n\tic@
 \immediate\openout\tic@=\jobname.tic
 \global\toc@true\fi}
\rightadd@\hl\to\nofrillslist@
\rightadd@\HL\to\overlonglist@
\def\HL@@C{\csname HL@C\HLlevel@\endcsname}
\def\HL@@P{\csname HL@P\HLlevel@\endcsname}
\def\HL@@Q{\csname HL@Q\HLlevel@\endcsname}
\def\HL@@S{\csname HL@S\HLlevel@\endcsname}
\def\HL@@N{\csname HL@N\HLlevel@\endcsname}
\def\HL@@F{\csname HL@F\HLlevel@\endcsname}
\def\HL@@@C{\csname\exxx@\HLtype@ @C\endcsname}
\def\HL@@@P{\csname\exxx@\HLtype@ @P\endcsname}
\def\HL@@@Q{\csname\exxx@\HLtype@ @Q\endcsname}
\def\HL@@@S{\csname\exxx@\HLtype@ @S\endcsname}
\def\HL@@@N{\csname\exxx@\HLtype@ @N\endcsname}
\def\HL#1{\expandafter
 \ifx\csname HL@C#1\endcsname\relax
  \DN@{\Err@{\string\HL#1 not defined in this style}}%
 \else
  \DN@{\gdef\HLlevel@{#1}\def\HLname@{\HL{#1}}\let\HLtype@\relax\FNSS@\HL@}%
 \fi
 \next@}%
\newif\ifquoted@
\let\aftertoc@\relax
\def\HL@{%
 \DN@"##1"##2\endHL{\def\entry@{##2}\quoted@true
  {\noexpands@
  \ifx\HLtype@\relax
   \let\pre\HL@@P\let\post\HL@@Q\let\style\HL@@S\let\numstyle\HL@@N
  \else
   \let\pre\HL@@@P\let\post\HL@@@Q\let\style\HL@@@S\let\numstyle\HL@@@N
  \fi
  \Qlabel@{##1}\let\style\relax\xdef\Qlabel@@@@{##1}%
  \xdef\Thepref@{\Thelabel@@@@}}%
  \csname HL@\HLlevel@\endcsname##2\endHL
  \let\pref\Thepref@
  \csname HL@I\HLlevel@\endcsname
  \csname HL@J\HLlevel@\endcsname
  \let\pref\pref@
  \HLtoc@
  \aftertoc@
  \let\aftertoc@\relax\overlong@false}%
 \DNii@##1\endHL{\def\entry@{##1}\quoted@false
  {\noexpands@
  \ifx\HLtype@\relax
   \global\advance\HL@@C\@ne
   \xdef\Thelabel@@@{\number\HL@@C}%
   \xdefThelabel@{\HL@@N}%
   \xdef\Thelabel@@@@{\HL@@P\Thelabel@\HL@@Q}%
   \xdefThelabel@@{\HL@@S}%
  \else
   \global\advance\HL@@@C\@ne
   \xdef\Thelabel@@@{\number\HL@@@C}%
   \xdefThelabel@{\HL@@@N}%
   \xdef\Thelabel@@@@{\HL@@@P\Thelabel@\HL@@@Q}%
   \xdefThelabel@@{\HL@@@S}%
  \fi
  \xdef\Thepref@{\Thelabel@@@@}}%
  \csname HL@\HLlevel@\endcsname##1\endHL
  \let\pref\Thepref@
  \csname HL@I\HLlevel@\endcsname
  \csname HL@J\HLlevel@\endcsname
  \let\pref\pref@
  \HLtoc@
  \aftertoc@
  \let\aftertoc@\relax\overlong@false}%
 \ifx\next"\expandafter\next@\else\expandafter\nextii@\fi}%
\Invalid@\endHL
\def\hl@@C{\csname hl@C\hllevel@\endcsname}
\def\hl@@P{\csname hl@P\hllevel@\endcsname}
\def\hl@@Q{\csname hl@Q\hllevel@\endcsname}
\def\hl@@S{\csname hl@S\hllevel@\endcsname}
\def\hl@@N{\csname hl@N\hllevel@\endcsname}
\def\hl@@F{\csname hl@F\hllevel@\endcsname}
\def\hl@@@C{\csname\exxx@\hltype@ @C\endcsname}
\def\hl@@@P{\csname\exxx@\hltype@ @P\endcsname}
\def\hl@@@Q{\csname\exxx@\hltype@ @Q\endcsname}
\def\hl@@@S{\csname\exxx@\hltype@ @S\endcsname}
\def\hl@@@N{\csname\exxx@\hltype@ @N\endcsname}
\def\hl#1{\expandafter
 \ifx\csname hl@C#1\endcsname\relax
  \DN@{\Err@{\string\hl#1 not defined in this style}}%
 \else
  \DN@{\gdef\hllevel@{#1}\def\hlname@{\hl{#1}}\let\hltype@\relax\FNSS@\hl@}%
 \fi
 \next@}
\def\hl@{%
 \DN@"##1"##2{\def\entry@{##2}\quoted@true
  {\noexpands@
  \ifx\hltype@\relax
   \let\pre\hl@@P\let\post\hl@@Q\let\style\hl@@S\let\numstyle\hl@@N
  \else
   \let\pre\hl@@@P\let\post\hl@@@Q\let\style\hl@@@S\let\numstyle\hl@@@N
  \fi
  \Qlabel@{##1}\let\style\relax\xdef\Qlabel@@@@{##1}%
  \xdef\Thepref@{\Thelabel@@@@}}%
  \csname hl@\hllevel@\endcsname{##2}%
  \let\pref\Thepref@
  \csname hl@I\hllevel@\endcsname
  \csname hl@J\hllevel@\endcsname
  \let\pref\pref@
  \hltoc@
  \aftertoc@
  \let\aftertoc@\relax\nopunct@false\nospace@false\FNSSP@}%
 \DNii@##1{\def\entry@{##1}\quoted@false
  {\noexpands@
  \ifx\hltype@\relax
   \global\advance\hl@@C\@ne
   \xdef\Thelabel@@@{\number\hl@@C}%
   \xdefThelabel@{\hl@@N}%
   \xdef\Thelabel@@@@{\hl@@P\Thelabel@\hl@@Q}%
   \xdefThelabel@@{\hl@@S}%
  \else
   \global\advance\hl@@@C\@ne
   \xdef\Thelabel@@@{\number\hl@@@C}%
   \xdefThelabel@{\hl@@@N}%
   \xdef\Thelabel@@@@{\hl@@@P\Thelabel@\hl@@@Q}%
   \xdefThelabel@@{\hl@@@S}%
  \fi
  \xdef\Thepref@{\Thelabel@@@@}}%
  \csname hl@\hllevel@\endcsname{##1}%
  \let\pref\Thepref@
  \csname hl@I\hllevel@\endcsname
  \csname hl@J\hllevel@\endcsname
  \let\pref\pref@
  \hltoc@
  \aftertoc@
  \let\aftertoc@\relax\nopunct@false\nospace@false\FNSSP@}%
 \ifx\next"\expandafter\next@\else\expandafter\nextii@\fi}%
\def\six@#1#2 #3 #4 #5 #6 #7 {\DN@{#2}\ifx\next@\empty
 \DN@##1\six@{}\else
 \write#1{ #2 #3 #4 #5 #6 #7}\DN@{\six@#1}\fi
 \next@}
\def\Sixtoc@{\ifx\macdef@\empty\else
 \DN@##1##2\next@{\def\macdef@{##1##2}}%
 \expandafter\next@\macdef@\next@
 \edef\next@
  {\noexpand\six@\toc@\macdef@
  \space\space\space\space\space\space\space\space\space\space\space\space
  \noexpand\six@}%
 \next@\let\macdef@\relax\fi}
\def\QorThelabel@@@@{\ifquoted@
 \noexpand\noexpand\noexpand"\Qlabel@@@@\noexpand\noexpand\noexpand"\else
 \Thelabel@@@@\fi}
\def\HLtoc@{%
 \iftoc@
 \expandafter\expandafter\expandafter\unmacro@
  \expandafter\meaning\csname HL@W\HLlevel@\endcsname\unmacro@
  {\noexpands@\let\style\relax
   \edef\next@{\write\toc@{\noexpand\noexpand\expandafter\noexpand\HLname@
   {\macdef@}{\QorThelabel@@@@}}}%
  \next@}%
  \expandafter\unmacro@\meaning\entry@\unmacro@
  \Sixtoc@
  \write\toc@{\noexpand\Page{\number\pageno}{\page@N}%
   {\page@P}{\page@Q}^^J}%
 \fi}
\def\hltoc@{%
 \iftoc@
 \expandafter\expandafter\expandafter\unmacro@
  \expandafter\meaning\csname hl@W\hllevel@\endcsname\unmacro@
  {\noexpands@\let\style\relax
  \edef\next@{\write\toc@{%
   \ifnopunct@\noexpand\noexpand\noexpand\nopunct\fi
   \ifnospace@\noexpand\noexpand\noexpand\nospace\fi
   \noexpand\noexpand\expandafter\noexpand\hlname@
   {\macdef@}{\QorThelabel@@@@}}}%
  \next@}%
  \expandafter\unmacro@\meaning\entry@\unmacro@
  \Sixtoc@
  \write\toc@{\noexpand\Page{\number\pageno}{\page@N}%
   {\page@P}{\page@Q}^^J}%
 \fi}
\def\mainfile#1{\def\mainfile@{#1}}
\def\checkmainfile@{\ifx\mainfile@\undefined
 \Err@{No \noexpand\mainfile specified}\fi}
\expandafter\newcount@\csname HL@C1\endcsname
\csname HL@C1\endcsname\z@
\expandafter\def\csname HL@S1\endcsname#1{#1\null.}
\expandafter\let\csname HL@N1\endcsname\arabic
\expandafter\let\csname HL@P1\endcsname\empty
\expandafter\let\csname HL@Q1\endcsname\empty
\expandafter\def\csname HL@F1\endcsname{\bf}
\expandafter\let\csname HL@W1\endcsname\empty
\expandafter\newcount@\csname hl@C1\endcsname
\csname hl@C1\endcsname\z@
\expandafter\def\csname hl@S1\endcsname#1{#1\/}
\expandafter\let\csname hl@N1\endcsname\arabic
\expandafter\let\csname hl@P1\endcsname\empty
\expandafter\let\csname hl@Q1\endcsname\empty
\expandafter\def\csname hl@F1\endcsname{\bf}
\expandafter\let\csname hl@W1\endcsname\empty
\expandafter\def\csname HL@1\endcsname#1\endHL{\bigbreak
 {\locallabel@
  \global\setbox\@ne\vbox{\Let@\tabskip\hss@
  \halign to\hsize{\bf\hfil\ignorespaces##\unskip\hfil\cr
  \expandafter\ifx\csname HL@W1\endcsname\empty\else
   \csname HL@W1\endcsname\space\fi
  {\HL@@F\ifx\thelabel@@\empty\else\thelabel@@\space\fi}%
  \ignorespaces#1\crcr}}%
  }%
 \unvbox\@ne\nobreak\medskip}
\expandafter\def\csname hl@1\endcsname#1{\medbreak\noindent@@
 {\locallabel@
 \bf{\hl@@F\ifx\thelabel@@\empty\else\thelabel@@\space\fi}%
 \ignorespaces#1\unskip\punct@{\null.}\addspace@\enspace}}
\expandafter\def\csname HL@I1\endcsname{\Reset\hl1{1}%
 \ifx\pref\empty\newpre\hl1{}\else\newpre\hl1{\pref.}\fi}
\def\NameHL#1#2{\define#2{}%
 \expandafter\ifx\csname HL@R#1\endcsname\relax
 \else
  \def\nextiv@{\let\nextiii@}%
  \expandafter\nextiv@\csname HL@R#1\endcsname
  \expandafter\let\nextiii@\undefined
  \expandafter\let\csname\exxx@\nextiii@ @C\endcsname\relax
  \expandafter\let\csname\exxx@\nextiii@ @P\endcsname\relax
  \expandafter\let\csname\exxx@\nextiii@ @Q\endcsname\relax
  \expandafter\let\csname\exxx@\nextiii@ @S\endcsname\relax
  \expandafter\let\csname\exxx@\nextiii@ @N\endcsname\relax
  \expandafter\let\csname\exxx@\nextiii@ @F\endcsname\relax
  \expandafter\let\csname\exxx@\nextiii@ @W\endcsname\relax
  \expandafter\let\csname end\exxx@\nextiii@\endcsname\undefined
 \fi
 \expandafter\gdef\csname HL@R#1\endcsname{#2}%
 \expandafter\gdef\csname\exstring@#2@R\endcsname{{HL}{#1}}%
 \iftoc@\write\toc@{\noexpand\NameHL#1\noexpand#2^^J}\fi
 \rightadd@#2\to\overlonglist@
 \edef\next@{\let\csname\exstring@#2@C\endcsname\expandafter\noexpand
  \csname HL@C#1\endcsname}\next@
 \edef\next@{\let\csname\exstring@#2@P\endcsname\expandafter\noexpand
  \csname HL@P#1\endcsname}\next@
 \edef\next@{\let\csname\exstring@#2@Q\endcsname\expandafter\noexpand
  \csname HL@Q#1\endcsname}\next@
 \edef\next@{\let\csname\exstring@#2@S\endcsname\expandafter\noexpand
  \csname HL@S#1\endcsname}\next@
 \edef\next@{\let\csname\exstring@#2@N\endcsname\expandafter\noexpand
  \csname HL@N#1\endcsname}\next@
 \edef\next@{\let\csname\exstring@#2@F\endcsname\expandafter\noexpand
  \csname HL@F#1\endcsname}\next@
 \edef\next@{\let\csname\exstring@#2@W\endcsname\expandafter\noexpand
  \csname HL@W#1\endcsname}\next@
 \edef\next@{\def\noexpand#2####1\expandafter\noexpand
  \csname end\exstring@#2\endcsname
  {\def\noexpand\HLtype@{\noexpand#2}%
   \def\noexpand\HLname@{\noexpand#2}%
   \gdef\noexpand\HLlevel@{#1}%
   \noexpand\FNSS@\noexpand\HL@####1\noexpand\endHL}}%
  \next@
 \edef\next@{\noexpand\Invalid@\expandafter\noexpand
  \csname end\exstring@#2\endcsname}%
 \next@}
\def\Namehl#1#2{\define#2{}%
 \expandafter\ifx\csname hl@R#1\endcsname\relax
 \else
  \def\nextiv@{\let\nextiii@}%
  \expandafter\nextiv@\csname hl@R#1\endcsname
  \expandafter\let\nextiii@\undefined
  \expandafter\let\csname\exxx@\nextiii@ @C\endcsname\relax
  \expandafter\let\csname\exxx@\nextiii@ @P\endcsname\relax
  \expandafter\let\csname\exxx@\nextiii@ @Q\endcsname\relax
  \expandafter\let\csname\exxx@\nextiii@ @S\endcsname\relax
  \expandafter\let\csname\exxx@\nextiii@ @N\endcsname\relax
  \expandafter\let\csname\exxx@\nextiii@ @F\endcsname\relax
  \expandafter\let\csname\exxx@\nextiii@ @W\endcsname\relax
 \fi
 \expandafter\gdef\csname hl@R#1\endcsname{#2}%
 \expandafter\gdef\csname\exstring@#2@R\endcsname{{hl}{#1}}%
 \iftoc@\write\toc@{\noexpand\Namehl#1\noexpand#2^^J}\fi
 \rightadd@#2\to\nofrillslist@%
 \edef\next@{\let\csname\exstring@#2@C\endcsname\expandafter\noexpand
  \csname hl@C#1\endcsname}\next@
 \edef\next@{\let\csname\exstring@#2@P\endcsname\expandafter\noexpand
  \csname hl@P#1\endcsname}\next@
 \edef\next@{\let\csname\exstring@#2@Q\endcsname\expandafter\noexpand
  \csname hl@Q#1\endcsname}\next@
 \edef\next@{\let\csname\exstring@#2@S\endcsname\expandafter\noexpand
  \csname hl@S#1\endcsname}\next@
 \edef\next@{\let\csname\exstring@#2@N\endcsname\expandafter\noexpand
  \csname hl@N#1\endcsname}\next@
 \edef\next@{\let\csname\exstring@#2@F\endcsname\expandafter\noexpand
  \csname hl@F#1\endcsname}\next@
 \edef\next@{\let\csname\exstring@#2@W\endcsname\expandafter\noexpand
  \csname hl@W#1\endcsname}\next@
 \edef\next@{\def\noexpand#2{%
  \def\noexpand\hltype@{\noexpand#2}%
  \def\noexpand\hlname@{\noexpand#2}%
  \gdef\noexpand\hllevel@{#1}%
  \noexpand\FNSS@\noexpand\hl@}}%
 \next@}%
\def\Initialize{\FN@\Init@}
\def\Init@{\ifx\next\HL\let\next@\InitH@\else\ifx\next\hl\let\next@\InitH@
  \else\let\next@\InitS@\fi\fi\next@}
\def\InitH@#1#2{\expandafter\ifx\csname\exstring@#1@C#2\endcsname\relax
 \DN@{\Err@{\noexpand#1level #2 not defined in this style}}\else
 \DN@{\expandafter\gdef\csname\exstring@#1@J#2\endcsname}\fi\next@}
\def\InitC@#1#2{\edef\nextii@{\expandafter\noexpand\csname#1\endcsname{#2}}}
\def\InitS@#1{\expandafter\ifx\csname\exstring@#1@R\endcsname\relax
 \Err@{\noexpand#1not defined in this style}\let\next@\relax\else
 \DN@{\let\next@}\expandafter\next@\csname\exstring@#1@R\endcsname
 \expandafter\InitC@\next@
 \DN@{\expandafter\InitH@\nextii@}\fi\next@}
\def\value#1{\expandafter
 \ifx\csname\exstring@#1@C\endcsname\relax
  \expandafter\ifx\csname\exstring@#1@C1\endcsname\relax
   \DN@{\Err@{\noexpand\value can't be used with \string#1}}%
  \else
   \DN@{\value@#1}%
  \fi
 \else
  \DN@{\number\csname\exstring@#1@C\endcsname\relax}%
 \fi
 \next@}
\def\value@#1#2{\expandafter
 \ifx\csname\exstring@#1@C#2\endcsname\relax
  \DN@{\Err@{\string\value\string#1 can't be followed by \string#2}}%
 \else
  \DN@{\number\csname\exstring@#1@C#2\endcsname\relax}%
 \fi
 \next@}
\newcount\Value
\def\Evaluate#1{\expandafter
 \ifx\csname\exstring@#1@C\endcsname\relax
  \expandafter\ifx\csname\exstring@#1@C1\endcsname\relax
   \DN@{\Err@{\noexpand\Evaluate can't be used with \string#1}}%
  \else
   \DN@{\Evaluate@#1}%
  \fi
 \else
  \DN@{\global\Value\csname\exstring@#1@C\endcsname}%
 \fi
 \next@}
\def\Evaluate@#1#2{\expandafter
 \ifx\csname\exstring@#1@C#2\endcsname\relax
  \DN@{\Err@{\string\Evaluate\string#1 can't be followed by \string#2}}%
 \else
  \DN@{\global\Value\csname\exstring@#1@C#2\endcsname}%
 \fi\next@}
\def\pre#1{\expandafter
 \ifx\csname\exstring@#1@P\endcsname\relax
  \expandafter\ifx\csname\exstring@#1@P1\endcsname\relax
   \DN@{\Err@{\noexpand\pre can't be used with \string#1}}%
  \else
   \DN@{\pre@#1}%
  \fi
 \else
  \DN@{{\csname\exstring@#1@P\endcsname}}%
 \fi
 \next@}
\def\pre@#1#2{\expandafter
 \ifx\csname\exstring@#1@P#2\endcsname\relax
  \DN@{\Err@{\string\pre\string#1 can't be followed by \string#2}}%
 \else
  \DN@{{\csname\exstring@#1@P#2\endcsname}}%
 \fi
 \next@}
\def\post#1{\expandafter
 \ifx\csname\exstring@#1@Q\endcsname\relax
  \expandafter\ifx\csname\exstring@#1@Q1\endcsname\relax
   \DN@{\Err@{\noexpand\post can't be used with \string#1}}%
  \else
   \DN@{\post@#1}%
  \fi
 \else
  \DN@{{\csname\exstring@#1@Q\endcsname}}%
 \fi
 \next@}
\def\post@#1#2{\expandafter
 \ifx\csname\exstring@#1@Q#2\endcsname\relax
  \DN@{\Err@{\string\post\string#1 can't be followed by \string#2}}%
 \else
  \DN@{{\csname\exstring@#1@Q#2\endcsname}}%
 \fi
 \next@}
\def\style#1{\expandafter
 \ifx\csname\exstring@#1@S\endcsname\relax
  \expandafter\ifx\csname\exstring@#1@S1\endcsname\relax
   \DN@{\Err@{\noexpand\style can't be used with \string#1}}%
  \else
   \DN@{\style@#1}%
  \fi
 \else
  \DN@{\csname\exstring@#1@S\endcsname}%
 \fi
 \next@}
\def\style@#1#2{\expandafter
 \ifx\csname\exstring@#1@S#2\endcsname\relax
  \DN@{\Err@{\string\style\string#1 can't be followed by \string#2}}%
 \else
  \DN@{\csname\exstring@#1@S#2\endcsname}%
 \fi
 \next@}
\def\fontstyle#1{\expandafter
 \ifx\csname\exstring@#1@F\endcsname\relax
  \expandafter\ifx\csname\exstring@#1@F1\endcsname\relax
   \DN@{\Err@{\noexpand\fontstyle can't be used with \string#1}}%
  \else
   \DN@{\fontstyle@#1}%
  \fi
 \else
  \DN@##1{{\csname\exstring@#1@F\endcsname##1}}%
 \fi
 \next@}
\def\fontstyle@#1#2{\expandafter
 \ifx\csname\exstring@#1@F#2\endcsname\relax
  \DN@{\Err@{\string\fontstyle\string#1 can't be followed by \string#2}}%
 \else
  \DN@##1{{\csname\exstring@#1@F#2\endcsname##1}}%
 \fi
 \next@}
\def\Reset#1{\expandafter
 \ifx\csname\exstring@#1@C\endcsname\relax
  \expandafter\ifx\csname\exstring@#1@C1\endcsname\relax
   \DN@{\Err@{\noexpand\Reset can't be used with \string#1}}%
  \else
   \DN@{\Reset@#1}%
  \fi
 \else
  \DN@##1{\count@##1\relax\ifx#1\page\else\advance\count@\m@ne\fi
   \global\csname\exstring@#1@C\endcsname\count@}%
 \fi
 \next@}
\def\Reset@#1#2{\expandafter
 \ifx\csname\exstring@#1@C#2\endcsname\relax
  \DN@{\Err@{\string\Reset\string#1 can't be followed by \string#2}}%
 \else
  \DN@##1{\count@##1\relax\advance\count@\m@ne
   \global\csname\exstring@#1@C#2\endcsname\count@}%
 \fi
 \next@}
\def\Offset#1{\expandafter
 \ifx\csname\exstring@#1@C\endcsname\relax
  \expandafter\ifx\csname\exstring@#1@C1\endcsname\relax
   \DN@{\Err@{\noexpand\Offset can't be used with \string#1}}%
  \else
   \DN@{\Offset@#1}%
  \fi
 \else
  \DN@##1{\count@##1\relax\advance\count@\m@ne\global\advance
   \csname\exstring@#1@C\endcsname\count@}%
 \fi
 \next@}
\def\Offset@#1#2{\expandafter
 \ifx\csname\exstring@#1@C#2\endcsname\relax
  \DN@{\Err@{\string\Offset\string#1 can't be followed by \string#2}}%
 \else
  \DN@##1{\count@##1\relax\advance\count@\m@ne
   \global\advance\csname\exstring@#1@C#2\endcsname\count@}%
 \fi
 \next@}
\def\getR@#1#2{\def\nextiv@{\let\nextiii@}\expandafter\nextiv@
 \csname\exstring@#1@R#2\endcsname}
\def\letR@#1#2#3{\expandafter\let\csname#1@#3#2\endcsname\Next@}
\def\letR@@#1#2{\expandafter\let\csname\exstring@#1@#2\endcsname\Next@}
\def\newpre#1{\expandafter
 \ifx\csname\exstring@#1@P\endcsname\relax
  \expandafter\ifx\csname\exstring@#1@P1\endcsname\relax
   \DN@{\Err@{\noexpand\newpre can't be used with \string#1}}%
  \else
   \DN@{\newpre@#1}%
  \fi
 \else
  \DN@{%
   \DNii@{%
    \endgroup
    \expandafter\let\csname\exstring@#1@P\endcsname\Next@
    \expandafter\ifx\csname\exstring@#1@R\endcsname\relax\else
    \getR@#1{}\expandafter\letR@\nextiii@ P\fi
    }%
   \begingroup\noexpands@\afterassignment\nextii@\xdef\Next@}%
 \fi
 \next@}
\def\newpre@#1#2{\expandafter
 \ifx\csname\exstring@#1@P#2\endcsname\relax
  \DN@{\Err@{\string\newpre\string#1 can't be followed by \string#2}}%
 \else
  \DN@{%
   \DNii@{%
    \endgroup
    \expandafter\let\csname\exstring@#1@P#2\endcsname\Next@
    \expandafter\ifx\csname\exstring@#1@R#2\endcsname\relax\else
    \getR@#1{#2}\expandafter\letR@@\nextiii@ P\fi
    }%
   \begingroup\noexpands@\afterassignment\nextii@\xdef\Next@}%
 \fi
 \next@}
\def\newpost#1{\expandafter
 \ifx\csname\exstring@#1@Q\endcsname\relax
  \expandafter\ifx\csname\exstring@#1@Q1\endcsname\relax
   \DN@{\Err@{\noexpand\newpost can't be used with \string#1}}%
  \else
   \DN@{\newpost@#1}%
  \fi
 \else
  \DN@{%
   \DNii@{%
    \endgroup
    \expandafter\let\csname\exstring@#1@Q\endcsname\Next@
    \expandafter\ifx\csname\exstring@#1@R\endcsname\relax\else
    \getR@#1{}\expandafter\letR@\nextiii@ Q\fi
    }%
   \begingroup\noexpands@\afterassignment\nextii@\xdef\Next@}%
 \fi
 \next@}
\def\newpost@#1#2{\expandafter
 \ifx\csname\exstring@#1@Q#2\endcsname\relax
  \DN@{\Err@{\string\newpost\string#1 can't be followed by \string#2}}%
 \else
  \DN@{%
   \DNii@{%
    \endgroup
    \expandafter\let\csname\exstring@#1@Q#2\endcsname\Next@
    \expandafter\ifx\csname\exstring@#1@R#2\endcsname\relax\else
    \getR@#1{#2}\expandafter\letR@@\nextiii@ Q\fi
    }%
   \begingroup\noexpands@\afterassignment\nextii@\xdef\Next@}%
 \fi
 \next@}
\def\newstyle#1{\expandafter
 \ifx\csname\exstring@#1@S\endcsname\relax
  \expandafter\ifx\csname\exstring@#1@S1\endcsname\relax
   \DN@{\Err@{\noexpand\newstyle can't be used
    with \string#1}}%
  \else
   \DN@{\newstyle@#1}%
  \fi
 \else
  \DN@{%
   \DNii@{%
    \expandafter\let\csname\exstring@#1@S\endcsname\Next@
    \expandafter\ifx\csname\exstring@#1@R\endcsname\relax\else
    \getR@#1{}\expandafter\letR@\nextiii@ S\fi
    }%
   \afterassignment\nextii@\gdef\Next@}%
 \fi
 \next@}
\def\newstyle@#1#2{\expandafter
 \ifx\csname\exstring@#1@S#2\endcsname\relax
  \DN@{\Err@{\string\newstyle\string#1 can't be followed by
   \string#2}}%
 \else
  \DN@{%
   \DNii@{%
    \expandafter\let\csname\exstring@#1@S#2\endcsname\Next@
    \expandafter\ifx\csname\exstring@#1@R#2\endcsname\relax\else
    \getR@#1{#2}\expandafter\letR@@\nextiii@ S\fi
    }%
   \afterassignment\nextii@\gdef\Next@}%
 \fi
 \next@}
\def\newnumstyle#1{\expandafter
 \ifx\csname\exstring@#1@N\endcsname\relax
  \expandafter\ifx\csname\exstring@#1@N1\endcsname\relax
   \DN@{\Err@{\noexpand\newnumstyle can't be used with
    \string#1}}%
  \else
   \DN@{\newnumstyle@#1}%
  \fi
 \else
  \DN@##1{%
   \gdef\Next@{##1}%
    \expandafter\let\csname\exstring@#1@N\endcsname\Next@
    \expandafter\ifx\csname\exstring@#1@R\endcsname\relax\else
    \getR@#1{}\expandafter\letR@\nextiii@ N\fi
    }%
 \fi
 \next@}
\def\newnumstyle@#1#2{\expandafter
 \ifx\csname\exstring@#1@N#2\endcsname\relax
  \DN@{\Err@{\string\newnumstyle\string#1 can't be followed by
   \string#2}}%
 \else
  \DN@##1{%
   \gdef\Next@{##1}%
    \expandafter\let\csname\exstring@#1@N#2\endcsname\Next@
    \expandafter\ifx\csname\exstring@#1@R#2\endcsname\relax\else
    \getR@#1{#2}\expandafter\letR@@\nextiii@ N\fi
    }%
  \fi
 \next@}
\def\newfontstyle#1{\expandafter
 \ifx\csname\exstring@#1@F\endcsname\relax
  \expandafter\ifx\csname\exstring@#1@F1\endcsname\relax
   \DN@{\Err@{\noexpand\newfontstyle can't be used with
    \string#1}}%
  \else
   \DN@{\newfontstyle@#1}%
  \fi
 \else
  \DN@##1{%
   \gdef\Next@{##1}%
    \expandafter\let\csname\exstring@#1@F\endcsname\Next@
    \expandafter\ifx\csname\exstring@#1@R\endcsname\relax\else
    \getR@#1{}\expandafter\letR@\nextiii@ F\fi
    }%
 \fi
 \next@}
\def\newfontstyle@#1#2{\expandafter
 \ifx\csname\exstring@#1@F#2\endcsname\relax
  \DN@{\Err@{\string\newfontstyle\string#1 can't be followed by
   \string#2}}%
 \else
  \DN@##1{%
   \gdef\Next@{##1}%
    \expandafter\let\csname\exstring@#1@F#2\endcsname\Next@
    \expandafter\ifx\csname\exstring@#1@R#2\endcsname\relax\else
    \getR@#1{#2}\expandafter\letR@@\nextiii@ F\fi
    }%
 \fi
 \next@}
\def\word#1{\expandafter
 \ifx\csname\exstring@#1@W\endcsname\relax
  \expandafter\ifx\csname\exstring@#1@W1\endcsname\relax
   \DN@{\Err@{\noexpand\word can't be used with \string#1}}%
  \else
   \DN@{\word@#1}%
  \fi
 \else
  \DN@{{\csname\exstring@#1@W\endcsname}}%
 \fi
 \next@}
\def\word@#1#2{\expandafter
 \ifx\csname\exstring@#1@W#2\endcsname\relax
  \DN@{\Err@{\string\word\noexpand#1can't be followed by \string#2}}%
 \else
  \DN@{{\csname\exstring@#1@W#2\endcsname}}%
 \fi
 \next@}
\def\newword#1{\expandafter
 \ifx\csname\exstring@#1@W\endcsname\relax
  \expandafter\ifx\csname\exstring@#1@W1\endcsname\relax
   \DN@{\Err@{\noexpand\newword can't be used  with \string#1}}%
  \else
   \DN@{\newword@#1}%
  \fi
 \else
  \DN@{%
   \DNii@{%
    \expandafter\let\csname\exstring@#1@W\endcsname\Next@
    \expandafter\ifx\csname\exstring@#1@R\endcsname\relax\else
     \getR@#1{}\expandafter\letR@\nextiii@ W\fi
    }%
   \afterassignment\nextii@\gdef\Next@}%
 \fi
 \next@}
\def\newword@#1#2{\expandafter
 \ifx\csname\exstring@#1@W#2\endcsname\relax
  \DN@{\Err@{\string\newword\noexpand#1can't be followed by \string#2}}%
 \else
  \DN@{%
   \DNii@{%
    \expandafter\let\csname\exstring@#1@W#2\endcsname\Next@
    \expandafter\ifx\csname\exstring@#1@R#2\endcsname\relax\else
     \getR@#1{#2}\expandafter\letR@@\nextiii@ W\fi
    }%
   \afterassignment\nextii@\gdef\Next@}%
 \fi
 \next@}
\newif\iffn@
\newcount\footmark@C
\footmark@C\z@
\def\footmark@S#1{$^{#1}$}
\let\footmark@N\arabic
\def\footmark@F{\rm}
\def\foottext@S#1{$^{#1}$}
\def\foottext@F{\rm}
\let\modifyfootnote@\relax
\def\modifyfootnote#1{\def\modifyfootnote@{#1}}
\def\vfootnote@#1{\insert\footins
 \bgroup
 \floatingpenalty\@MM\interlinepenalty\interfootnotelinepenalty
 \leftskip\z@\rightskip\z@\spaceskip\z@\xspaceskip\z@
 \rm\splittopskip\ht\strutbox\splitmaxdepth\dp\strutbox
 \locallabel@\noindent@@{\foottext@F#1}\modifyfootnote@
 \footstrut\FN@\fo@t}
\def\fo@t{\ifcat\bgroup\noexpand\next\expandafter\f@@t\else
 \expandafter\f@t\fi}
\def\f@t#1{#1\@foot}
\def\f@@t{\bgroup\aftergroup\@foot\afterassignment\FNSSP@\let\next@}
\def\@foot{\unskip\lower\dp\strutbox\vbox to\dp\strutbox{}\egroup
 \iffn@\expandafter\fn@false\else
 \expandafter\postvanish@\fi}
\newif\ifplainfn@
\plainfn@true
\def\fancyfootnotes{\plainfn@false}
\newcount\fancyfootmarkcount@
\fancyfootmarkcount@\z@
\newcount\lastfnpage@
\lastfnpage@-\@M
\let\justfootmarklist@\empty
\def\footmark{\let\@sf\empty
 \ifhmode\edef\@sf{\spacefactor\the\spacefactor}\/\fi
 \DN@{\ifx"\next\expandafter\nextii@\else\expandafter\footmark@\fi}%
 \DNii@"##1"{%
  \iffirstchoice@
   {\let\style\footmark@S\let\numstyle\footmark@N
   \footmark@F##1%
   \noexpands@
   \let\style\foottext@S
   \Qlabel@{##1}%
   }%
   \iffn@\else
    {\noexpands@
    \xdef\Next@{{\Thelabel@}{\Thelabel@@}{\Thelabel@@@}{\Thelabel@@@@}}%
    }%
    \expandafter\rightappend@\Next@\to\justfootmarklist@
   \fi
  \fi
  \@sf\relax}%
 \FN@\next@}
\def\footmark@{%
 \iffirstchoice@
  \global\advance\footmark@C\@ne
  \ifplainfn@
   \xdef\adjustedfootmark@{\number\footmark@C}%
  \else
   {\let\\\or\xdef\Next@{\ifcase\number\footmark@C\fnpages@\else
     -\@M\fi}}%
   \ifnum\Next@=-\@M
    \xdef\adjustedfootmark@{\number\footmark@C}%
   \else
    \ifnum\Next@=\lastfnpage@
     \global\advance\fancyfootmarkcount@\@ne
    \else
     \global\fancyfootmarkcount@\@ne
     \global\lastfnpage@\Next@
    \fi
    \xdef\adjustedfootmark@{\number\fancyfootmarkcount@}%
   \fi
  \fi
  {\noexpands@
  \xdef\Thelabel@@@{\adjustedfootmark@}%
  \xdefThelabel@\footmark@N
  \xdef\Thelabel@@@@{\Thelabel@}%
  \xdefThelabel@@\foottext@S
  }%
  \iffn@\else
   {\noexpands@
   \xdef\Next@{{\Thelabel@}{\Thelabel@@}{\Thelabel@@@}{\Thelabel@@@@}}%
   }%
   \expandafter\rightappend@\Next@\to\justfootmarklist@
  \fi
  \ifplainfn@
  \else
   \edef\next@{\write\laxwrite@{F\noexpand\the\pageno}}\next@
  \fi
 \fi
 \footmark@S{\footmark@N{\adjustedfootmark@}}%
 \@sf\relax}
\def\foottext{\prevanish@
 \ifx\justfootmarklist@\empty
  \Err@{There is no \noexpand\footmark for this \string\foottext}\fi
 \DN@\\##1##2\next@{\DN@{##1}\gdef\justfootmarklist@{##2}}%
 \expandafter\next@\justfootmarklist@\next@
 \expandafter\foottext@\next@}
\def\foottext@#1#2#3#4{{\noexpands@
  \xdef\Thelabel@{#1}\xdef\Thelabel@@{#2}%
  \xdef\Thelabel@@@{#3}\xdef\Thelabel@@@@{#4}}%
  \vfootnote@{\thelabel@@}}
\rightadd@\foottext\to\vanishlist@
\def\footnote{\fn@true
 \let\@sf\empty
 \ifhmode\edef\@sf{\spacefactor\the\spacefactor}\/\fi
 \DN@{\ifx"\next\expandafter\nextii@\else\expandafter\nextiii@\fi}%
 \DNii@"##1"{\footmark"##1"\vfootnote@{\let\style\foottext@S
  \let\numstyle\footmark@N##1}}%
 \def\nextiii@{\footmark\vfootnote@{\foottext@S{\footmark@N
  {\adjustedfootmark@}}}}%
 \FN@\next@}
\newdimen\litindent
\litindent20\p@
\newbox\litbox@
\newbox\Litbox@
\newcount\interlitpenalty@
\interlitpenalty@\@M
\newcount\litlines@
{\obeyspaces\gdef\defspace@{\def {\allowbreak\hskip.5emminus.15em}}}
{\obeylines\gdef\letM@{\let^^M\CtrlM@}}
\def\CtrlM@{\egroup
 \ifcase\litlines@\advance\litlines@\@ne\or
 \box\litbox@\advance\litlines@\@ne\else
 \penalty\interlitpenalty@\box\litbox@\fi
 \Lit@}
\def\Lit@{\setbox\litbox@\hbox\bgroup\litdefs@\hskip\litindent}
\newcount\littab@
\littab@8
\def\littab#1{\littab@#1\relax}
{\catcode`\^^I=\active\gdef\letTAB@{\let^^I\TAB@}}
\def\TAB@{\egroup
 \dimen@\wd\litbox@
 \advance\dimen@-\litindent
 \setboxz@h{\tt0}%
 \dimen@ii\littab@\wdz@
 \divide\dimen@\dimen@ii
 \multiply\dimen@\dimen@ii
 \advance\dimen@\littab@\wdz@
 \advance\dimen@\litindent
 \setbox\litbox@\hbox\bgroup\litdefs@\hbox to\dimen@{\unhbox\litbox@\hfil}}
{\catcode`\`=\active\gdef`{\relax\lq}}
\let\litbs@\relax
\let\litbs@@\relax
\def\litbackslash#1{%
 \edef\litbs@{\catcode`\string#1=\z@
 \def\noexpand\litbs@@{\def\expandafter\noexpand\csname\string#1\endcsname
  {\char`\string#1}}}}
\def\litcodes@{\catcode`\\=12
 \catcode`\{=12 \catcode`\}=12
 \catcode`\$=12 \catcode`\&=12
 \catcode`\#=12
 \catcode`\^=12 \catcode`\_=12
 \catcode`\@=12 \catcode`\~=12 \catcode`\"=12
 \catcode`\;=12 \catcode`\:=12 \catcode`\!=12 \catcode`\?=12
 \catcode`\%=12 \litbs@\catcode`\`=\active\obeyspaces\defspace@}
\def\activate@#1#2{{\lccode`\~=`#2%
 \lowercase{%
  \if0#1%
  \gdef\Next@{\def~{\egroup\endgroup\bigskip\vskip-\parskip
   \def\next@{\noindent@@\FN@\pretendspace@}\FNSS@\next@}}\else
  \gdef\Next@{\def~{\egroup\egroup\endgroup}}\fi
  }%
 }}
\def\litdefs@{\let\0\empty\let\1\litdelim@\def\ {\char32 }\litbs@@}%
\def\litdelimiter#1{%
 \edef\litdelim@{\char`#1}%
 \def\lit#1{\leavevmode\begingroup\litcodes@\litdefs@
  \tt\hyphenchar\tentt\m@ne\lit@}%
 \def\lit@##1#1{##1\endgroup\null}%
 \def\Lit#1{\ifhmode$$\abovedisplayskip\bigskipamount
  \abovedisplayshortskip\bigskipamount
  \belowdisplayskip\z@\belowdisplayshortskip\z@
  \postdisplaypenalty\@M
  $$\vskip-\baselineskip\else\bigskip\fi
  \begingroup\litlines@\z@
  \catcode`#1=\active\activate@0#1\Next@
  \def\displaybreak{\egroup\break\litlines@\z@\Lit@}%
  \def\allowdisplaybreak{\egroup\allowbreak\litlines@\z@\Lit@}%
  \def\allowdisplaybreaks{\egroup\allowbreak\interlitpenalty@\z@
   \litlines@\z@\Lit@}%
  \litcodes@\tt\catcode`\^^I=\active\letTAB@
  \obeylines\letM@\Lit@}%
 \def\Litbox##1=#1{\begingroup\ifodd##1\relax\aftergroup\global\fi
  \aftergroup\setbox\aftergroup##1\aftergroup\box\aftergroup\Litbox@
  \def\allowdisplaybreak{\egroup\allowbreak\litlines@\z@\Lit@}%
  \def\allowdisplaybreaks{\egroup\allowbreak\interlitpenalty@\z@
   \litlines@\z@\Lit@}%
  \catcode`#1=\active\activate@1#1\Next@
  \litcodes@\tt\catcode`\^^I=\active\letTAB@
  \obeylines\letM@\global\setbox\Litbox@\vbox\bgroup\litindent\z@%
  \litlines@\z@\Lit@}%
}
\newbox\titlebox@
\setbox\titlebox@\vbox{}
\rightadd@\title\to\overlonglist@
\def\title{\begingroup\Let@
 \global\setbox\titlebox@\vbox\bgroup\tabskip\hss@
 \halign to\hsize\bgroup\bf\hfil\ignorespaces##\unskip\hfil\cr}
\def\endtitle{\crcr\egroup\egroup\endgroup\overlong@false}
\newbox\authorbox@
\rightadd@\author\to\overlonglist@
\def\author{\begingroup\Let@
 \global\setbox\authorbox@\vbox\bgroup\tabskip\hss@
 \halign to\hsize\bgroup\rm\hfil\ignorespaces##\unskip\hfil\cr}
\def\endauthor{\crcr\egroup\egroup\endgroup\overlong@false}
\newbox\affilbox@
\def\affil{\begingroup\Let@
 \global\setbox\affilbox@\vbox\bgroup\tabskip\hss@
 \halign to\hsize\bgroup\rm\hfil\ignorespaces##\unskip\hfil\cr}%
\def\endaffil{\crcr\egroup\egroup\endgroup\overlong@false}
\let\date@\relax
\def\date#1{\gdef\date@{\ignorespaces#1\unskip}}
\def\today{\ifcase\month\or January\or February\or March\or April\or May\or
 June\or July\or August\or September\or October\or November\or December\fi
 \space\number\day, \number\year}
\def\maketitle{\hrule\height\z@\vskip-\topskip
 \vskip24\p@ plus12\p@ minus12\p@
 \unvbox\titlebox@
 \ifvoid\authorbox@\else\vskip12\p@ plus6\p@ minus3\p@\unvbox\authorbox@\fi
 \ifvoid\affilbox@\else\vskip10\p@ plus5\p@ minus2\p@\unvbox\affilbox@\fi
 \ifx\date@\relax\else\vskip6\p@ plus2\p@ minus\p@\centerline{\rm\date@}\fi
 \vskip18\p@ plus12\p@ minus6\p@}
\def\cite{%
 \DNii@(##1)##2{{\rm[}{##2}, {##1\/}{\rm]}}%
 \def\nextiii@##1{{\rm[}{##1\/}{\rm]}}%
 \DN@{\ifx\next(\expandafter\nextii@\else\expandafter\nextiii@\fi}%
 \FN@\next@}
\def\makebib@W{Bibliography}
\def\makebibp@W{References}
\def\makebib{\begingroup\rm\bigbreak\centerline{\smc\makebib@W}%
 \nobreak\medskip
 \sfcode`\.=\@m\everypar{}\parindent\z@
 \def\nopunct{\nopunct@true}\def\nospace{\nospace@true}%
 \nopunct@false\nospace@false
 \def\lkerns@{\null\kern\m@ne sp\kern\@ne sp}%
 \def\nkerns@{\null\kern-\tw@ sp\kern\tw@ sp}%
}

\newif\ifnoprepunct@
\newif\ifnoprespace@
\newif\ifnoquotes@
\def\noprepunct{\noprepunct@true}
\def\noprespace{\noprespace@true}
\def\noquotes{\noquotes@true}
\newbox\nobox@
\newbox\keybox@
\newbox\bybox@
\newbox\paperbox@
\newbox\paperinfobox@
\newbox\jourbox@
\newbox\volbox@
\newbox\issuebox@
\newbox\yrbox@
\newbox\pgbox@
\newbox\ppbox@
\newbox\bookbox@
\newbox\inbookbox@
\newbox\bookinfobox@
\newbox\publbox@
\newbox\publaddrbox@
\newbox\edbox@
\newbox\edsbox@
\newbox\langbox@
\newbox\translbox@
\newbox\finalinfobox@
\def\setbibinfo@#1{\edef\next@{\ifnopunct@1\else0\fi
 \ifnospace@1\else0\fi\ifnoprepunct@1\else0\fi\ifnoprespace@1\else0\fi
 \ifnoquotes@1\else0\fi}%
 \DNii@{00000}%
 \ifx\next@\nextii@\else\xdef\bibinfo@{\bibinfo@\the#1,\next@}%
 \fi}
\def\getbibinfo@#1{\ifx\bibinfo@\empty
 \let\next@0\let\nextii@0\let\nextiii@0\let\nextiv@0\let\nextv@0\else
 \edef\next@{\def
  \noexpand\next@####1\the#1,####2####3####4####5####6####7\noexpand\next@
  {\let\noexpand\next@####2\let\noexpand\nextii@####3%
  \let\noexpand\nextiii@####4\let\noexpand\nextiv@####5%
  \let\noexpand\nextv@####6}%
  \noexpand\next@\bibinfo@\the#1,00000\noexpand\next@}\next@
 \fi}
\newif\ifbookinquotes@
\def\bookinquotes{\bookinquotes@true}
\newif\ifpaperinquotes@
\def\paperinquotes{\paperinquotes@true}
\newif\ifininbook@
\def\ininbook{\ininbook@true}
\newif\ifopenquotes@
\def\closequotes@{\ifopenquotes@''\openquotes@false\fi}
\newif\ifbeginbib@
\newif\ifendbib@
\newif\ifprevjour@
\newif\ifprevbook@
\newdimen\bibindent@
\bibindent@20\p@
\def\bib{\global\let\bibinfo@\empty\global\let\translinfo@\relax\beginbib@true
 \begingroup\noindent@
 \hangindent\bibindent@\hangafter\@ne\bib@}
\def\v@id#1{\setbox#1\box\voidb@x}
\def\bib@{\v@id\nobox@\v@id\keybox@\v@id\bybox@\v@id\paperbox@
 \v@id\paperinfobox@\v@id\jourbox@\v@id\volbox@\v@id\issuebox@
 \v@id\yrbox@\v@id\pgbox@\v@id\ppbox@\v@id\bookbox@\v@id\inbookbox@
 \v@id\bookinfobox@\v@id\publbox@\v@id\publaddrbox@\v@id\edbox@
 \v@id\edsbox@\v@id\langbox@\v@id\translbox@\v@id\finalinfobox@
 \bgroup}
\def\Setnonemptybox@#1#2{\unskip\setbibinfo@#1\egroup#2%
 \def\aftergroup@{\ifdim\wd#1=\z@\setbox#1\box\voidb@x\fi}%
 \setbox#1\vbox\bgroup\aftergroup\aftergroup@\hsize\maxdimen\leftskip\z@
 \rightskip\z@\hbadness\@M\hfuzz\maxdimen\noindent}
\def\setnonemptybox@#1{\Setnonemptybox@#1\relax}
\def\no{\setnonemptybox@\nobox@}
\def\key{\setnonemptybox@\keybox@\bf}
\def\by{\setnonemptybox@\bybox@}
\def\bysame{\setnonemptybox@\bybox@\leaders\hrule\hskip3em\null}
\def\paper{\setnonemptybox@\paperbox@
 \ifpaperinquotes@\getbibinfo@\paperbox@
 \if\nextv@1\else``\fi\else\it\fi}
\def\paperinfo{\setnonemptybox@\paperinfobox@}
\def\jour{\Setnonemptybox@\jourbox@\prevjour@true}
\def\vol{\setnonemptybox@\volbox@\bf}
\def\issue{\setnonemptybox@\issuebox@}
\def\yr{\setnonemptybox@\yrbox@}

\def\pg{\setnonemptybox@\pgbox@}
\def\pp{\setnonemptybox@\ppbox@}
\def\book{\Setnonemptybox@\bookbox@\prevbook@true
 \ifbookinquotes@\getbibinfo@\bookbox@
 \if\nextv@1\else``\fi\else\it\fi}
\def\inbook{\Setnonemptybox@\inbookbox@\prevbook@true
 \ifininbook@ in \fi\ifbookinquotes@\getbibinfo@\inbookbox@
 \if\nextv@1\else``\fi\fi}
\def\bookinfo{\setnonemptybox@\bookinfobox@}
\def\publ{\setnonemptybox@\publbox@}
\def\publaddr{\setnonemptybox@\publaddrbox@}
\def\ed{\setnonemptybox@\edbox@}
\def\eds{\setnonemptybox@\edsbox@}
\def\lang{\setnonemptybox@\langbox@}
\def\finalinfo{\setnonemptybox@\finalinfobox@}
\def\setboxzl@{\setbox\z@\lastbox}
\def\getbox@#1{\setbox\z@\vbox{\vskip-\@M\p@
 \unvbox#1%
 \setboxzl@
 \global\setbox\@ne\hbox{\unhbox\z@\unskip\unskip\unpenalty}%
 \ifdim\lastskip=-\@M\p@\else
 \loop\ifdim\lastskip=-\@M\p@
 \else\unskip\unpenalty\setboxzl@
 \global\setbox\@ne\hbox{\unhbox\z@\unhbox\@ne}%
 \repeat\fi}%
 \unhbox\@ne}
\def\adjustpunct@#1{\count@\lastkern
 \ifnum\count@=\z@#1\closequotes@\else
 \ifnum\count@>\tw@#1\closequotes@\else
 \ifnum\count@<-\tw@#1\closequotes@\else
  \unkern\unkern\setboxzl@
  \skip@\lastskip\unskip
  \count@@\lastpenalty\unpenalty
  \ifnum\count@=\tw@\unskip\setboxzl@\fi
  \ifdim\skip@=\z@\else\hskip\skip@\fi
  #1\closequotes@
  \ifnum\count@=\tw@\null\hfill\fi
  \penalty\count@@
 \fi\fi\fi}
\def\prepunct@#1#2{\getbibinfo@#2%
 \ifnopunct@
 \else
  \if\nextiii@0\adjustpunct@#1\fi
 \fi
 \closequotes@
 \ifnospace@
 \else
  \if\nextiv@0\space\else\fi
 \fi
 \nopunct@false\nospace@false
 \if\next@1\nopunct@true\fi
 \if\nextii@1\nospace@true\fi}
\def\ppunbox@#1#2{\prepunct@{#1}#2%
 \getbox@#2}
\let\semicolon@;
\def\endbib@{%
 \ifbeginbib@
  \ifvoid\nobox@
   \ifvoid\keybox@\else\hbox to\bibindent@{[\getbox@\keybox@]\hss}\fi
  \else\hbox to\bibindent@{\hss\getbox@\nobox@. }\fi
  \ifvoid\bybox@\else\getbox@\bybox@\fi
 \else
  \nopunct@true
  \ifvoid\bybox@\else\ppunbox@\relax\bybox@\fi
 \fi
 \ifvoid\translbox@\else\ppunbox@,\translbox@\fi
 \ifvoid\paperbox@\else\ppunbox@,\paperbox@\ifpaperinquotes@
  \if\nextv@1\else\openquotes@true\fi\fi
 \fi
 \ifvoid\paperinfobox@\else\ppunbox@,\paperinfobox@\fi
 \test@false
 \ifvoid\jourbox@\else\test@true\ppunbox@,\jourbox@\fi
 \ifprevjour@\test@true\fi
 \iftest@
  \ifvoid\volbox@\else\ppunbox@\relax\volbox@\fi
  \ifvoid\issuebox@
   \else\prepunct@\relax\issuebox@ no.~\getbox@\issuebox@\fi
  \ifvoid\yrbox@\else\prepunct@\relax\yrbox@(\getbox@\yrbox@)\fi
  \ifvoid\ppbox@\else\ppunbox@,\ppbox@\fi
  \ifvoid\pgbox@\else\prepunct@,\pgbox@ p.~\getbox@\pgbox@\fi
 \fi
 \test@false
 \ifvoid\bookbox@\else\test@true\ppunbox@,\bookbox@\ifbookinquotes@
  \if\nextv@1\else\openquotes@true\fi\fi\fi
 \ifvoid\inbookbox@\else\test@true\ppunbox@,\inbookbox@\ifbookinquotes@
  \if\nextv@1\else\openquotes@true\fi\fi\fi
 \ifprevbook@\test@true\fi
 \iftest@
  \ifvoid\edbox@\else\prepunct@\relax\edbox@(\getbox@\edbox@, ed.)\fi
  \ifvoid\edsbox@\else\prepunct@\relax\edsbox@(\getbox@\edsbox@, eds.)\fi
  \ifvoid\bookinfobox@\else\ppunbox@,\bookinfobox@\fi
  \ifvoid\publbox@\else\ppunbox@,\publbox@\fi
  \ifvoid\publaddrbox@\else\ppunbox@,\publaddrbox@\fi
  \ifvoid\yrbox@\else\ppunbox@,\yrbox@\fi
  \ifvoid\ppbox@\else\prepunct@,\ppbox@ pp.~\getbox@\ppbox@\fi
  \ifvoid\pgbox@\else\prepunct@,\pgbox@ p.~\getbox@\pgbox@\fi
 \fi
 \ifvoid\finalinfobox@
  \ifendbib@
   \ifnopunct@\else.\closequotes@\fi
  \else
  \ifvoid\langbox@\else\space(\getbox@\langbox@)\fi
   \/\semicolon@\closequotes@
  \fi
 \else
  \ifendbib@
   \ppunbox@{.\spacefactor3000\relax}\finalinfobox@
    \ifnopunct@\else.\fi
  \else
   \ppunbox@,\finalinfobox@\/\semicolon@\fi
 \fi
 \ifvoid\langbox@\else\space(\getbox@\langbox@)\fi
}
\def\endbib{\unskip\egroup\endbib@true\endbib@\par\endgroup}
\def\morebib{\unskip\egroup
 \endbib@false\endbib@
 \global\let\bibinfo@\empty\beginbib@false
 \bib@}
\def\anotherbib{\unskip\egroup
 \endbib@false\endbib@
 \global\let\bibinfo@\empty\beginbib@false
 \prevjour@false\prevbook@false\bib@}
\def\transl{\unskip
 \xdef\translinfo@{\the\translbox@,\ifnopunct@1\else0\fi
 \ifnospace@1\else0\fi\ifnoprepunct@1\else0\fi\ifnoprespace@1\else0\fi0}%
 \egroup\endbib@false\endbib@
 \global\let\bibinfo@\translinfo@\beginbib@false
 \bib@
 \egroup
 \def\aftergroup@{\ifdim\wd\translbox@=\z@\setbox\translbox@\box\voidb@x\fi}%
 \setbox\translbox@\vbox\bgroup\aftergroup\aftergroup@
 \hsize\maxdimen\leftskip\z@\rightskip\z@\hbadness\@M\hfuzz\maxdimen
 \noindent}
\newwrite\auxwrite@
\newread\bbl@
\def\UseBibTeX{\immediate\openout\auxwrite@=\jobname.aux
 \let\cite\BTcite@
 \def\nocite##1{\immediate\write\auxwrite@{\string\citation{##1}}}%
 \def\bibliographystyle##1{\immediate\write\auxwrite@{\string
  \bibstyle{##1}}}%
 \def\bibliography@W{Bibliography}%
 \def\bibliography##1{\immediate\write\auxwrite@{\string\bibdata{##1}}%
  \immediate\openin\bbl@=\jobname.bbl
  \ifeof\bbl@
   \W@{No .bbl file}%
  \else
   \immediate\closein\bbl@
   \begingroup\input bibtex \input\jobname.bbl \endgroup
  \fi}%
 }
\def\BTcite@{%
 \DNii@(##1)##2{{\rm[}\BTcite@@##2,\BTcite@@{\rm, }{##1\/}{\rm]}%
  \immediate\write\auxwrite@{\string\citation{##2}}}%
 \def\nextiii@##1{{\rm[}\BTcite@@##1,\BTcite@@\/{\rm]}%
  \immediate\write\auxwrite@{\string\citation{##1}}}%
 \DN@{\ifx\next(\expandafter\nextii@\else\expandafter\nextiii@\fi}%
 \FN@\next@}%
\def\BTcite@@#1,{\BTcite@@@{#1}\FN@\BTcite@@@@}
\def\BTcite@@@@{\ifx\next\BTcite@@
 \expandafter\eat@\else{\rm, }\expandafter\BTcite@@\fi}
\catcode`\~=11
\def\BTcite@@@#1{\nolabel@\cite{#1}\relax
 \DNii@##1~##2\nextii@{##1}%
 \csL@{#1}\expandafter\nextii@\Next@\nextii@\fi}
\catcode`\~=\active

\def\beginthebibliography@#1{\rm\setboxz@h{#1\ }\bibindent@\wdz@
 \bigbreak\centerline{\smc\bibliography@W}\nobreak\medskip
 \sfcode`\.=\@m\everypar{}\parindent\z@}

\newif\iffigproofing@
\def\Figureproofing{\figproofing@true}
\def\noFigureproofing{\figproofing@false}
\newif\ifHby@
\def\Hbyw#1{\global\Hby@true\hbyw\vsize{#1}}
\def\hbyw#1#2{%
 \hbox{%
  \ifHby@
  \else
   \iffigproofing@
    \setbox\z@\vbox{\hrule\width5\p@}\ht\z@\z@
    \vbox to#1{\hrule\height5\p@\width.4\p@\vfil\hrule\height5\p@\width.4\p@}%
    \kern-.4\p@\rlap{\copy\z@}\raise#1\hbox{\rlap{\copy\z@}}%
   \fi
  \fi
  \vbox to#1{\hbox to#2{}\vfil}%
  \ifHby@
  \else
   \iffigproofing@
    \vbox to#1{\hrule\height5\p@\width.4\p@\vfil\hrule\height5\p@\width.4\p@}%
    \kern-.4\p@\llap{\copy\z@}\raise#1\hbox{\llap{\boxz@}}%
   \fi
  \fi}}
\newcount\island@C
\let\island@P\empty
\let\island@Q\empty
\def\island@S#1{#1\null.}
\let\island@N\arabic
\def\island@F{\rm}
\def\island@@@P{\csname\exxx@\islandtype@ @P\endcsname}
\def\island@@@Q{\csname\exxx@\islandtype@ @Q\endcsname}
\def\island@@@S{\csname\exxx@\islandtype@ @S\endcsname}
\def\island@@@N{\csname\exxx@\islandtype@ @N\endcsname}
\def\island@@@F{\csname\exxx@\islandtype@ @F\endcsname}
\def\island@@@C{\csname island@C\islandclass@\endcsname}
\newif\ifplace@
\newif\ifisland@
\def\island{%
 \ifplace@
  \DN@{\let\islandclass@\empty\def\islandtype@{\island}\FN@\island@}%
 \else
  \long\DN@##1\endisland{\Err@{\noexpand\island must be used after some
   type of \string\...place}}%
 \fi
 \next@}
\def\island@{\ifx\next\c\let\next@\island@c\else
 \DN@{\FN@\island@@}\fi\next@}
\def\island@@{\ifcat\bgroup\noexpand\next\let\next@\island@@@\else
 \DN@{\Err@{\noexpand\island must be followed by a {prefix} for
 \string\caption's}}\fi\next@}
\newbox\islandbox@
\newcount\captioncount@
\def\island@@@#1{\def\captionprefix@{#1}\captioncount@\z@
 \global\setbox\islandbox@\vbox\bgroup}
\def\island@c\c#1{%
 \ifplace@
 \DN@{\def\islandclass@{#1}%
  \expandafter\ifx\csname island@C#1\endcsname\relax
  \expandafter\newcount@\csname island@C#1\endcsname
   \global\csname island@C#1\endcsname\z@\fi
  \FNSS@\island@c@}%
 \else
 \DN@{\edef\next@{\long\def\noexpand\next@########1\expandafter\noexpand
  \csname end\exxx@\islandtype@\endcsname{\noexpand\Err@{\noexpand\noexpand
  \expandafter\noexpand
  \islandtype@ must be used after some type of \noexpand\string
   \noexpand\...place}}}\next@\next@}%
 \fi
 \next@}
\def\island@c@{%
 \ifcat\bgroup\noexpand\next
  \let\next@\island@c@@
 \else
  \DN@{\Err@{\noexpand\island\string\c{\expandafter\string\islandclass@} must
   be followed by a {prefix} for \string\caption's}}%
 \fi\next@}
\def\island@c@@#1{\def\captionprefix@{#1}%
 \captioncount@\z@\global\setbox\islandbox@\vbox\bgroup}
\rightadd@\caption\to\nofrillslist@
\newbox\captionbox@
\newbox\Captionbox@
\def\caption{%
 \ifnum\captioncount@=\z@
  \ifnopunct@
   \DN@{\egroup\nopunct@true}%
  \else
   \let\next@\egroup
  \fi
 \else
  \let\next@\relax
 \fi
 \next@
 \advance\captioncount@\@ne
 \FN@\caption@}
\def\caption@{\ifx\next"\expandafter\caption@q\else\expandafter\caption@@\fi}
\def\caption@q"#1"{\quoted@true
 {\noexpands@
 \let\pre\island@@@P\let\post\island@@@Q
 \let\style\island@@@S\let\numstyle\island@@@N
 \Qlabel@{#1}\let\style\relax\xdef\Qlabel@@@@{#1}}%
 \finishcaption@}
\def\caption@@{\quoted@false
 \global\advance\island@@@C\@ne
 {\noexpands@
 \xdef\Thelabel@@@{\number\island@@@C}%
 \xdefThelabel@\island@@@N
 \xdef\Thelabel@@@@{\island@@@P\Thelabel@\island@@@Q}%
 \xdefThelabel@@\island@@@S
 \xdef\Thepref@{\Thelabel@@@@}}%
 \finishcaption@}
\long\def\captionformat@#1#2#3{\rm\strut#1 {\island@@@F#2} #3%
 \punct@.\strut}
\long\def\widerthanisland@#1#2#3{\test@true\setbox\z@\vbox{\hsize\maxdimen
 \noindent@@\captionformat@{#1}{#2}{#3}\par\setboxzl@}%
 \ifdim\wdz@=\z@
  \global\setbox\captionbox@\hbox{\noset@\unlabel@
   \captionformat@{#1}{#2}{#3}}%
  \ifdim\wd\captionbox@>\wd\islandbox@\else\test@false\fi
 \fi}
\long\def\captionformat@@#1#2#3{\widerthanisland@{#1}{#2}{#3}%
 \iftest@
  \global\setbox\captionbox@\vbox{\hsize\wd\islandbox@
   \vskip-\parskip\noindent@@\noset@\unlabel@
   \captionformat@{#1}{#2}{#3}\par}%
 \else
  \global\setbox\captionbox@
   \hbox to\wd\islandbox@{\hfil\box\captionbox@\hfil}%
 \fi}
\long\def\finishcaption@#1{\def\entry@{#1}%
 {\locallabel@
 \captionformat@@
  {\expandafter\ignorespaces\captionprefix@\unskip}%
  {\ifx\thelabel@@\empty\unskip\else\thelabel@@\fi}%
  {\ignorespaces#1\unskip}%
 \ifnum\captioncount@=\@ne
  \global\setbox\islandbox@\vbox{\ticwrite@\vbox{\box\islandbox@}}%
  \global\setbox\Captionbox@\vbox{\box\captionbox@}%
 \else
  \global\setbox\islandbox@\vbox{\unvbox\islandbox@\setboxzl@
   \ticwrite@\boxz@}%
  \global\setbox\Captionbox@\vbox{\unvbox\Captionbox@
   \smallskip\box\captionbox@}%
 \fi}%
 \nopunct@false\nospace@false\ignorespaces}
\def\Sixtic@{\ifx\macdef@\empty\else
 \DN@##1##2\next@{\def\macdef@{##1##2}}%
 \expandafter\next@\macdef@\next@
 \edef\next@
  {\noexpand\six@\tic@\macdef@
  \space\space\space\space\space\space\space\space\space\space\space\space
  \noexpand\six@}%
 \next@\let\macdef@\relax\fi}
\def\ticwrite@{%
 \iftoc@
  {\noexpands@\let\style\relax
  \DN@{\island}%
  \edef\next@{\write\tic@{%
   \ifnopunct@\noexpand\noexpand\noexpand\nopunct\fi
   \ifx\islandtype@\next@\noexpand\noexpand\noexpand\island
    \noexpand\string\noexpand\c{\islandclass@}{\captionprefix@}%
     {\QorThelabel@@@@}\else\noexpand\noexpand\expandafter\noexpand
     \islandtype@{\QorThelabel@@@@}}\fi}%
  \next@}%
  \expandafter\unmacro@\meaning\entry@\unmacro@
  \Sixtic@
  \write\tic@{\noexpand\Page{\number\pageno}{\page@N}{\page@P}{\page@Q}^^J}%
 \fi}
\def\Htrim@#1{%
 \ifHby@
  \dimen@\vsize
  \ifnum\captioncount@=\z@
  \else
   \advance\dimen@-\ht\Captionbox@
   \advance\dimen@-#1%
  \fi
  \global\Hby@false
  \dimen@ii\wd\islandbox@
  \global\setbox\islandbox@\vbox
   {\unvbox\islandbox@\setboxzl@
   \vbox to\z@{\vss\boxz@}\nointerlineskip\hbyw\dimen@\dimen@ii}%
  \global\Hby@true
 \fi}
\newif\ifdata@
\def\iclasstest@#1{\DN@{#1}\ifx\next@\islandclass@
 \test@true\else\test@false\fi}
\skipdef\skipi@=1
\def\endisland{\ifnum\captioncount@=\z@\expandafter\egroup\fi
 \ifdata@
 \else
  \iclasstest@{T}%
  \iftest@
   {\rm\global\skipi@-\dp\strutbox}\global\advance\skipi@\bigskipamount
   \Htrim@\skipi@
   \global\setbox\islandbox@\vbox
    {\ifnum\captioncount@=\z@\else
     \box\Captionbox@
     \nointerlineskip
     \vskip\skipi@\fi
     \box\islandbox@}%
  \else
   {\rm\global\skipi@\dp\strutbox}\global\advance\skipi@\medskipamount
   \Htrim@\skipi@
   \global\setbox\islandbox@\vbox
    {\box\islandbox@
     \ifnum\captioncount@=\z@\else
     \nointerlineskip
     \vskip\skipi@
     \box\Captionbox@
     \fi}%
  \fi
  \ifHby@
  \else
   \dimen@\ht\islandbox@\advance\dimen@\dp\islandbox@
   \ifdim\dimen@>\vsize
    \DN@{\island}%
    \Err@{%
     \ifx\islandtype@\next@\noexpand\island\else
      \expandafter\noexpand\islandtype@\fi
     \ifnum\captioncount@=\z@\else
       with \noexpand\caption\fi
      is larger than page}%
     \ht\islandbox@=\vsize
   \fi
  \fi
 \fi
 \global\Hby@false\island@true}
\def\newisland#1\c#2#3{\define#1{}%
 \iftoc@\immediate\write\tic@{\noexpand\newisland\noexpand#1%
  \string\c{#2}{#3}^^J}\fi
 \expandafter\def\csname\exstring@#1@S\endcsname{\island@S}%
 \expandafter\def\csname\exstring@#1@N\endcsname{\island@N}%
 \expandafter\def\csname\exstring@#1@P\endcsname{\island@P}%
 \expandafter\def\csname\exstring@#1@Q\endcsname{\island@Q}%
 \expandafter\def\csname\exstring@#1@F\endcsname{\island@F}%
 \expandafter\def\csname end\exstring@#1\endcsname{\endisland}%
 \expandafter
 \ifx\csname island@C#2\endcsname\relax
  \expandafter\newcount@\csname island@C#2\endcsname
  \global\csname island@C#2\endcsname\z@
 \fi
 \edef\next@{\noexpand\expandafter\noexpand\let\noexpand
  \csname\exstring@#1@C\noexpand\endcsname
  \csname island@C#2\endcsname}%
 \next@
 \def#1{\def\islandtype@{#1}\island@c\c{#2}{#3}}}
\newisland\Figure\c{F}{Figure}
\newisland\Table\c{T}{Table}
\newbox\islandboxi
\newbox\islandboxii
\newbox\islandboxiii
\newbox\captionboxi
\newbox\captionboxii
\newbox\captionboxiii
\long\def\islandpairdata#1#2{{\data@true
 \place@true
 #1%
 \global\setbox\islandboxi\box\islandbox@
 \global\setbox\captionboxi\box\Captionbox@
 #2%
 \global\setbox\islandboxii\box\islandbox@
 \global\setbox\captionboxii\box\Captionbox@
 }}
\long\def\islandpairbox#1#2{\islandpairdata{#1}{#2}%
 \dimen@\ht\captionboxi
 \ifdim\ht\captionboxii>\dimen@\dimen@\ht\captionboxii\fi
 \ifdim\dimen@>\z@
  \ifdim\ht\captionboxi<\dimen@
   \global\setbox\captionboxi\vbox to\dimen@{\unvbox\captionboxi\vfil}\fi
  \ifdim\ht\captionboxii<\dimen@
   \global\setbox\captionboxii\vbox to\dimen@{\unvbox\captionboxii\vfil}\fi
 \fi
 \global\setbox\islandbox@\vbox
 {\hbox to\hsize{\hfil\box\islandboxi\hfil\box\islandboxii\hfil}%
 \ifdim\dimen@>\z@\nointerlineskip
 {\rm\global\skipi@\dp\strutbox}\global\advance\skipi@\medskipamount
  \vskip\skipi@
  \hbox to\hsize{\hfil\box\captionboxi\hfil\box\captionboxii\hfil}\fi}}
\long\def\islandpairboxa#1#2{\islandpairdata{#1}{#2}%
 \dimen@\ht\captionboxi
 \ifdim\ht\captionboxii>\dimen@\dimen@\ht\captionboxii\fi
 \ifdim\dimen@>\z@
  \ifdim\ht\captionboxi<\dimen@
   \global\setbox\captionboxi\vbox to\dimen@{\vfil\unvbox\captionboxi}\fi
  \ifdim\ht\captionboxii<\dimen@
   \global\setbox\captionboxii\vbox to\dimen@{\vfil\unvbox\captionboxii}\fi
 \fi
 \dimen@ii\ht\islandboxi
 \ifdim\ht\islandboxii>\dimen@ii \dimen@ii\ht\islandboxii\fi
 \ifdim\dimen@ii>\z@
  \ifdim\ht\islandboxi<\dimen@ii
   \global\setbox\islandboxi\vbox to\dimen@ii{\box\islandboxi\vfil}\fi
  \ifdim\ht\islandboxii<\dimen@ii
   \global\setbox\islandboxii\vbox to\dimen@ii{\box\islandboxii\vfil}\fi
 \fi
 \global\setbox\islandbox@\vbox{\ifdim\dimen@>\z@
  \hbox to\hsize{\hfil\box\captionboxi\hfil\box\captionboxii\hfil}%
  \nointerlineskip{\rm\global\skipi@-\dp\strutbox}%
  \global\advance\skipi@\bigskipamount\vskip\skipi@\fi
  \hbox to\hsize{\hfil\box\islandboxi\hfil\box\islandboxii\hfil}}}
\long\def\islandtripledata#1#2#3{{\data@true\place@true
 #1%
 \global\setbox\islandboxi\box\islandbox@
 \global\setbox\captionboxi\box\Captionbox@
 #2%
 \global\setbox\islandboxii\box\islandbox@
 \global\setbox\captionboxii\box\Captionbox@
 #3%
 \global\setbox\islandboxiii\box\islandbox@
 \global\setbox\captionboxiii\box\Captionbox@
 }}
\long\def\islandtriplebox#1#2#3{\islandtripledata{#1}{#2}{#3}%
 \dimen@\ht\captionboxi
 \ifdim\ht\captionboxii>\dimen@ \dimen@\ht\captionboxii\fi
 \ifdim\ht\captionboxiii>\dimen@ \dimen@\ht\captionboxiii\fi
 \ifdim\dimen@>\z@
  \ifdim\ht\captionboxi<\dimen@
   \global\setbox\captionboxi\vbox to\dimen@{\unvbox\captionboxi\vfil}\fi
  \ifdim\ht\captionboxii<\dimen@
   \global\setbox\captionboxii\vbox to\dimen@{\unvbox\captionboxii\vfil}\fi
  \ifdim\ht\captionboxiii<\dimen@
   \global\setbox\captionboxiii\vbox to\dimen@{\unvbox\captionboxiii\vfil}\fi
 \fi
 \global\setbox\islandbox@\vbox
  {\hbox to\hsize{\hfil\box\islandboxi\hfil\box\islandboxii\hfil
   \box\islandboxiii\hfil}%
 \ifdim\dimen@>\z@\nointerlineskip
  {\rm\global\skipi@\dp\strutbox}\global\advance\skipi@\medskipamount
  \vskip\skipi@
  \hbox to\hsize{\hfil\box\captionboxi\hfil\box\captionboxii\hfil
   \box\captionboxiii\hfil}\fi}}
\def\islandtripleboxa#1#2#3{\islandtripledata{#1}{#2}{#3}%
 \dimen@\ht\captionboxi
 \ifdim\ht\captionboxii>\dimen@ \dimen@\ht\captionboxii\fi
 \ifdim\ht\captionboxiii>\dimen@ \dimen@\ht\captionboxiii\fi
 \ifdim\dimen@>\z@
  \ifdim\ht\captionboxi<\dimen@
   \global\setbox\captionboxi\vbox to\dimen@{\vfil\unvbox\captionboxi}\fi
  \ifdim\ht\captionboxii<\dimen@
   \global\setbox\captionboxii\vbox to\dimen@{\vfil\unvbox\captionboxii}\fi
  \ifdim\ht\captionboxiii<\dimen@
   \global\setbox\captionboxiii\vbox to\dimen@{\vfil\unvbox\captionboxiii}\fi
 \fi
 \dimen@ii\ht\islandboxi
 \ifdim\ht\islandboxii>\dimen@ii \dimen@ii\ht\islandboxii\fi
 \ifdim\ht\islandboxiii>\dimen@ii \dimen@ii\ht\islandboxiii\fi
 \ifdim\dimen@ii>\z@
  \ifdim\ht\islandboxi<\dimen@ii
   \global\setbox\islandboxi\vbox to\dimen@ii{\box\islandboxi\vfil}\fi
  \ifdim\ht\islandboxii<\dimen@ii
   \global\setbox\islandboxii\vbox to\dimen@ii{\box\islandboxii\vfil}\fi
  \ifdim\ht\islandboxiii<\dimen@ii
   \global\setbox\islandboxiii\vbox to\dimen@ii{\box\islandboxiii\vfil}\fi
 \fi
 \global\setbox\islandbox@\vbox
  {\ifdim\dimen@>\z@
  \hbox to\hsize{\hfil\box\captionboxi\hfil\box\captionboxii\hfil
   \box\captionboxiii\hfil}%
  \nointerlineskip{\rm\global\skipi@-\dp\strutbox}%
  \global\advance\skipi@\bigskipamount\vskip\skipi@\fi
  \hbox to\hsize{\hfil\box\islandboxi\hfil\box\islandboxii\hfil
   \box\islandboxiii\hfil}}}
\def\Figurepair#1\and#2\endFigurepair{\island@true
 \islandpairbox{\Figure#1\endFigure}{\Figure#2\endFigure}}
\def\Figuretriple#1\and#2\and#3\endFiguretriple{\island@true
 \islandtriplebox{\Figure#1\endFigure}{\Figure#2\endFigure}%
  {\Figure#3\endFigure}}
\def\Tablepair#1\and#2\endTablepair{\island@true
 \islandpairboxa{\Table#1\endTable}{\Table#2\endTable}}
\def\Tabletriple#1\and#2\and#3\endTabletriple{\island@true
 \islandtripleboxa{\Table#1\endTable}{\Table#2\endTable}%
 {\Table#3\endTable}}
\def\place#1{\place@true\island@false
 #1%
 \ifisland@
  \box\islandbox@
 \else
  \Err@{Whoa ... there's no \string\Figure, \string\Table,
   etc., here}%
 \fi
 \place@false}
\newskip\belowtopfigskip
\belowtopfigskip 15\p@ plus 5\p@ minus5\p@
\newskip\abovebotfigskip
\abovebotfigskip 18\p@ plus 6\p@ minus6\p@
\newdimen\minpagesize
\minpagesize 5pc
\dimen@\belowtopfigskip
\advance\dimen@-\abovebotfigskip
\skip\topins\dimen@
\dimen\topins\z@
\newcount\topinscount@
\newbox\topinsdims@
\def\storedim@{\global\setbox\topinsdims@
 \vbox{\hbox to\dimen@{}\unvbox\topinsdims@}}
\def\advancedimtopins@{%
 \ifnum\pageno=\@ne
 \else
   \advance\dimen@\dimen\topins
   \global\dimen\topins\dimen@
 \fi}
\newcount\flipcount@
\def\fliptopins@{%
 \global\flipcount@\z@
 \ifvoid\topins\else
 \setbox\z@\vbox
  {\vskip\p@
   \unvbox\topins
   \global\setbox\topins\vbox{}%
   \loop
    \test@false
    \ifdim\lastskip=\z@\unskip
     \ifdim\lastskip=\z@
      \test@true\fi\fi
    \iftest@
    \global\advance\flipcount@\@ne
    \setboxzl@
    \global\setbox\topins\vbox{\unvbox\topins\boxz@}%
    \unpenalty
   \repeat}\fi}
\newif\ifPar@
\newcount\Parcount@
\newbox\Parbox@
\expandafter\newbox\csname Parfigbox1\endcsname
\expandafter\newbox\csname Parfigbox2\endcsname
\expandafter\newbox\csname Parfigbox3\endcsname
\expandafter\newbox\csname Parfigbox4\endcsname
\expandafter\newbox\csname Parfigbox5\endcsname
\expandafter\newdimen\csname Parprev1\endcsname
\expandafter\newdimen\csname Parprev2\endcsname
\expandafter\newdimen\csname Parprev3\endcsname
\expandafter\newdimen\csname Parprev4\endcsname
\expandafter\newdimen\csname Parprev5\endcsname
\expandafter\newdimen\csname Parprev6\endcsname
\def\Par{\par\global\csname Parprev1\endcsname\prevdepth
 \global\Parcount@\@ne
 \global\Par@true\global\let\Parlist@\empty
 \global\setbox\Parbox@\vbox\bgroup\break}
\def\place@#1#2{%
 \ifisland@
  \ifhmode
   \ifPar@
    \ifnum\Parcount@>5
     \Err@{Only 5 \string\place's allowed per
      \string\Par...\noexpand\endPar paragraph}%
    \else
     \expandafter\expandafter\expandafter
      \global\expandafter\setbox
       \csname Parfigbox\number\Parcount@\endcsname\box\islandbox@
     \global\advance\Parcount@\@ne
     \xdef\Parlist@{\Parlist@#1}%
    \fi
   \else
    \vadjust{#2}%
   \fi
  \else
   #2%
  \fi
 \else
  \Err@{Whoa ... there's no \string\Figure,
   \string\Table, etc., here}%
 \fi
 \place@false}
\long\def\Aplace#1{\prevanish@
 \place@true\island@false
 #1%
 \place@ a\Aplace@
 \postvanish@}
\long\def\AAplace#1{\prevanish@\place@true\island@false
 #1%
 \place@ A\AAplace@
 \postvanish@}
\newif\ifAA@
\def\AAplace@{\AA@true\Aplace@\AA@false}
\let\AAlist@\empty
\def\Aplace@{\allowbreak
 \dimen@=\ht\islandbox@
 \advance\dimen@\abovebotfigskip
 \ht\islandbox@\dimen@
 \advance\dimen@\dp\islandbox@
 \storedim@
 \ifAA@
  \xdef\AAlist@{\AAlist@1}%
  \advancedimtopins@
 \else
  \xdef\AAlist@{\AAlist@0}%
  \ifnum\topinscount@>\@ne\else\advancedimtopins@\fi
 \fi
 \insert\topins{\penalty\z@\splittopskip\z@\floatingpenalty\z@
  \box\islandbox@}%
 \global\advance\topinscount@\@ne}
\long\def\Bplace#1{\prevanish@\place@true\island@false
 #1%
 \place@ b\Bplace@
 \postvanish@}
\def\Bplace@{\allowbreak
 \ifnum\topinscount@=\z@
  \setbox\z@\vbox{\vbox to-\belowtopfigskip{}}%
  \dimen@-\skip\topins
  \ht\z@\dimen@
  \storedim@
  \advancedimtopins@
  \insert\topins{\boxz@}%
  \global\advance\topinscount@\@ne
  \xdef\AAlist@{\AAlist@0}%
 \fi
 \dimen@\ht\islandbox@
 \advance\dimen@\abovebotfigskip
 \ht\islandbox@\dimen@
 \advance\dimen@\dp\islandbox@
 \storedim@
 \xdef\AAlist@{\AAlist@0}%
 \ifnum\topinscount@>\@ne\else\advancedimtopins@\fi
 \insert\topins{\penalty\z@\splittopskip\z@
  \floatingpenalty\z@
  \box\islandbox@}%
 \global\advance\topinscount@\@ne}
\def\breakisland@{\global\setbox\@ne\lastbox\global\skipi@\lastskip\unskip
 \global\setbox\thr@@\lastbox}%
\def\printisland@{\centerline{\box\thr@@}\nobreak\nointerlineskip
 \vskip\skipi@
 \ifdim\ht\@ne<\z@\box\@ne\else\centerline{\box\@ne}\fi}
\def\bottomfigs@{%
 \count@\@ne
 \loop
  \ifnum\count@<\flipcount@
  \nointerlineskip
  \vskip\abovebotfigskip
  \global\setbox\topins\vbox{\unvbox\topins\setboxzl@
   \unvbox\z@
   \breakisland@}%
  \printisland@
  \advance\count@\@ne
  \repeat}
\def\resetdimtopins@{%
 \global\advance\topinscount@-\flipcount@
 \global\setbox\topinsdims@\vbox
  {\unvbox\topinsdims@
   \count@\z@
   \DN@##1##2\next@{\gdef\AAlist@{##2}}%
   \loop
    \ifnum\count@<\flipcount@\setboxzl@
    \expandafter\next@\AAlist@\next@
    \advance\count@\@ne
    \repeat
   \dimen@\z@
   \count@\z@
   \setbox\tw@\vbox{}%
   \edef\nextiii@{\AAlist@}%
   \DN@##1##2\next@{\DNii@{##1}\def\nextiii@{##2}}%
   \loop
    \test@false
    \ifnum\count@<\topinscount@
    \expandafter\next@\nextiii@\next@
     \ifnum\count@<\tw@
      \test@true
     \else
      \if\nextii@ 1\test@true\fi
     \fi
    \fi
    \iftest@
     \setboxzl@
     \advance\dimen@\wdz@
     \setbox\tw@\vbox{\boxz@\unvbox\tw@}%
     \advance\count@\@ne
    \repeat
    \unvbox\tw@
    \global\dimen\topins\dimen@}}
\def\Place@#1#2{%
 \ifisland@
  \ifhmode
   \ifPar@
    \ifnum\Parcount@>5
     \Err@{Only 5 \string\place's allowed per
       \string\Par...\noexpand\endPar paragraph}%
    \else
     \expandafter\expandafter\expandafter\global\expandafter\setbox
      \csname Parfigbox\number\Parcount@\endcsname\box\islandbox@
     \global\advance\Parcount@\@ne
     \xdef\Parlist@{\Parlist@#1}%
     \vadjust{\break}%
    \fi
   \else
    \Err@{\noexpand#2allowed only in a \string\Par...\noexpand\endPar
     paragraph}%
   \fi
  \else
   #2%
  \fi
 \else
  \Err@{Who ... there's no \string\Figure, \string\Table,
   etc., here}%
 \fi
 \place@false}
\newif\ifC@
\newdimen\Cdim@
\long\def\Cplace#1{\prevanish@\place@true\island@false
 #1%
 \Place@ c\Cplace@
 \postvanish@}
\def\Cplace@{\allowbreak
 \ifnum\topinscount@>\z@\else
  \global\C@true\global\Cdim@\pagetotal\fi
 \Aplace@}
\long\def\Mplace#1{\prevanish@\place@true\island@false
 #1%
 \Place@ m\Mplace@
 \postvanish@}
\long\def\MXplace#1{\prevanish@\place@true\island@false
 #1%
 \Place@ M\MXplace@
 \postvanish@}
\newif\ifMX@
\def\MXplace@{\MX@true\Mplace@\MX@false}
\def\Mplace@{\allowbreak
 \dimen@\ht\islandbox@\advance\dimen@\dp\islandbox@
 \ifdim\pagetotal=\z@\else
  \ifdim\lastskip<\abovebotfigskip\advance\dimen@\abovebotfigskip
  \advance\dimen@-\lastskip\fi
 \fi
 \advance\dimen@\pagetotal
 \ifdim\dimen@>\pagegoal
  \Aplace@
 \else
  \nointerlineskip
  \ifdim\lastskip<\abovebotfigskip\removelastskip\vskip\abovebotfigskip\fi
  \setbox\z@\vbox{\unvbox\islandbox@
   \breakisland@}%
  \printisland@
  \ifnum\topinscount@=\z@
   \setbox\z@\vbox{\vbox to-\belowtopfigskip{}}%
   \dimen@-\skip\topins
   \ht\z@\dimen@
   \storedim@
   \advancedimtopins@
   \insert\topins{\boxz@}%
   \global\advance\topinscount@\@ne
   \xdef\AAlist@{\AAlist@0}%
  \fi
  \ifMX@
   \ifnum\topinscount@=\@ne
    \setbox\z@\vbox{\vbox to-\abovebotfigskip{}}%
    \ht\z@\z@
    \dimen@\z@
    \storedim@
    \advancedimtopins@
    \insert\topins{\boxz@}%
    \global\advance\topinscount@\@ne
    \xdef\AAlist@{\AAlist@0}%
   \fi
  \fi
  \nointerlineskip
  \vskip\belowtopfigskip
 \fi}
\expandafter\newbox\csname Parbox1\endcsname
\expandafter\newbox\csname Parbox2\endcsname
\expandafter\newbox\csname Parbox3\endcsname
\expandafter\newbox\csname Parbox4\endcsname
\expandafter\newbox\csname Parbox5\endcsname
\def\endPar{\egroup
 \count@\@ne
 {\vbadness\@M\vfuzz\maxdimen\splitmaxdepth\maxdimen\splittopskip\ht\strutbox
 \setbox\z@\vsplit\Parbox@ to\ht\Parbox@
 \loop
  \ifnum\count@<\Parcount@
  \expandafter\expandafter\expandafter\global\expandafter\setbox
   \csname Parbox\number\count@\endcsname\vsplit\Parbox@ to\ht\Parbox@
  \count@@\count@\advance\count@@\@ne
  \global\csname Parprev\number\count@@\endcsname
   \dp\csname Parbox\number\count@\endcsname
  \advance\count@\@ne
  \repeat}%
 \vskip\parskip
 \count@\@ne
 \def\nextv@##1##2\nextv@{\DN@{##1}\gdef\Parlist@{##2}}%
 \loop
  \ifnum\count@<\Parcount@
   \dimen@\csname Parprev\number\count@\endcsname
   \advance\dimen@\ht\strutbox
   \ifdim\dimen@<\baselineskip
    \advance\dimen@-\baselineskip\vskip-\dimen@
   \else
    \vskip\lineskip
   \fi
   \unvbox\csname Parbox\number\count@\endcsname
   \global\setbox\islandbox@\box\csname Parfigbox\number\count@\endcsname
   \expandafter\nextv@\Parlist@\nextv@
   \if a\next@\Aplace@\else
   \if A\next@\AAplace@\else
   \if b\next@\Bplace@\else
   \if c\next@\Cplace@\else
   \if m\next@\Mplace@\else
   \if M\next@\MXplace@\fi\fi\fi\fi\fi\fi
  \advance\count@\@ne
  \repeat
 \global\Par@false
 \ifvoid\Parbox@
  \prevdepth\csname Parprev\number\count@\endcsname
 \else
  \dimen@\csname Parprev\number\count@\endcsname\advance\dimen@\ht\strutbox
  \ifdim\dimen@<\baselineskip
    \advance\dimen@-\baselineskip\vskip-\dimen@
  \else
    \vskip\lineskip
  \fi
  \dimen@\dp\Parbox@
  \unvbox\Parbox@
  \prevdepth\dimen@
 \fi}
\def\folio{{\page@F\page@S{\page@P\page@N{\number\page@C}\page@Q}}}
\def\advancepageno{\global\advance\pageno\@ne}
\newif\ifspecialsplit@
\newbox\outbox@
\let\shipout@\shipout
\def\plainoutput{\specialsplit@false\ifvoid\topins\else\ifdim\ht\topins=\z@
 \specialsplit@true\advance\minpagesize-\skip\topins\fi\fi
 \fliptopins@
 \setbox\outbox@\vbox{\makeheadline\pagebody\makefootline}%
 {\noexpands@\let\style\relax
 \shipout@\box\outbox@}%
 \advancepageno
 \resetdimtopins@
 \ifvoid\@cclv\else\unvbox\@cclv\penalty\outputpenalty\fi
 \ifnum\outputpenalty>-\@MM\else\dosupereject\fi}
\def\pagebody{\vbox to\vsize{\boxmaxdepth\maxdepth
 \ifvoid\margin@\else
 \rlap{\kern\hsize\vbox to\z@{\kern4\p@\box\margin@\vss}}\fi
 \pagecontents}}
\newif\ifonlytop@
\def\pagecontents{%
 \onlytop@false
 \ifdim\ht\@cclv<\minpagesize\ifnum\flipcount@<\tw@\ifvoid\footins
  \onlytop@true\fi\fi\fi
 \test@false
 \ifC@
  \ifnum\flipcount@=\@ne
   \global\multiply\Cdim@\tw@
   \ifdim\Cdim@>\ht\@cclv
    \test@true
   \fi
  \fi
 \fi
 \global\C@false
 \iftest@
  \dimen@\ht\@cclv
  \advance\dimen@\skip\topins
  {\vfuzz\maxdimen\vbadness\@M
  \splitmaxdepth\maxdepth\splittopskip\topskip
  \setbox\z@\vsplit\@cclv to\dimen@
  \unvbox\z@}%
  \global\setbox\topins\vbox{\unvbox\topins
   \global\setbox\@ne\lastbox}%
  \setbox\z@\vbox{\unvbox\@ne
   \breakisland@}%
  \nointerlineskip
  \vskip\abovebotfigskip
  \printisland@
 \else
  \ifnum\flipcount@>\z@
   \global\setbox\topins\vbox{\unvbox\topins\global\setbox\@ne\lastbox}%
   \setbox\z@\vbox{\unvbox\@ne
    \breakisland@}%
   \printisland@
   \ifonlytop@\kern-\prevdepth\vfill\else\vskip\belowtopfigskip\fi
  \fi
 \fi
 \ifdim\ht\@cclv<\minpagesize
  \ifonlytop@\else\vfill\fi
 \else
  \ifspecialsplit@
   {\vfuzz\maxdimen\vbadness\@M
   \splitmaxdepth\maxdepth\splittopskip\topskip
   \dimen@ii\ht\@cclv \advance\dimen@ii\skip\topins
   \setbox\z@\vsplit\@cclv to\dimen@ii
   \unvbox\z@}%
  \else
   \unvbox\@cclv
  \fi
 \fi
 \bottomfigs@
 \ifvoid\footins\else\vskip\skip\footins\footnoterule\unvbox\footins\fi}
\newread\readdata@
\def\readthedata@#1{\expandafter
 \ifx\csname#1@D\endcsname\relax
  \immediate\openin\readdata@=#1.dat
  \ifeof\readdata@
   \Err@{No file #1.dat}%
  \else
   {\endlinechar\m@ne\gdef\Next@{}%
   \DNii@##1 ##2 ##3pt{\global\data@ht##1\global\data@dp##2%
    \global\data@wd##3pt}%
   \loop
    \ifeof\readdata@
    \else
    \read\readdata@ to\next@
    \ifx\next@\empty\else
     \edef\next@{\expandafter\nextii@\next@}%
     \expandafter\rightadd@\next@\to\Next@
    \fi
    \repeat}%
   \immediate\closein\readdata@
   \expandafter\expandafter\expandafter\global\expandafter
    \let\csname#1@D\endcsname\Next@\global\let\Next@\relax
  \fi
 \fi}
\newdimen\data@ht
\newdimen\data@dp
\newdimen\data@wd
\newif\ifgetdata@
\def\getdata@#1#2{\global\getdata@true\count@#2\relax
 {\let\\\or\xdef\Next@{\ifcase\number\count@#1\else
 \global\noexpand\getdata@false\fi}}\Next@}
\def\paste#1#2{\readthedata@{#1}%
 \getdata@{\csname#1@D\endcsname}{#2}%
 \ifgetdata@
 \dimen@\data@ht \advance\dimen@\data@dp
  \hbox{\special{dvipaste: #1 #2}%
   \lower\data@dp\vbox to\dimen@{\hbox to\data@wd{}\vfil}}%
 \else
  {\lccode`\Z=`\#\lccode`\N=`\N\lccode`\F=`\F%
   \lowercase{\Err@{No data for File [#1], Z#2}}}%
 \fi}
\newdimen\httable
\newdimen\dptable
\newdimen\wdtable
\def\measuretable#1#2{\readthedata@{#1}%
 \getdata@{\csname#1@D\endcsname}{#2}%
 \ifgetdata@
  \httable\data@ht \dptable\data@dp \wdtable\data@wd
 \else
  {\lccode`\Z=`\#\lccode`\N=`\N\lccode`\F=`\F%
  \lowercase{\Err@{No data for File [#1], Z#2}}}%
 \fi}
\def\East#1#2{\setboxz@h{$\m@th\ssize\;{#1}\;\;$}%
 \setbox\tw@\hbox{$\m@th\ssize\;{#2}\;\;$}\setbox4=\hbox{$\m@th#2$}%
 \dimen@\minaw@
 \ifdim\wdz@>\dimen@\dimen@\wdz@\fi\ifdim\wd\tw@>\dimen@\dimen@\wd\tw@\fi
 \ifdim\wd4 >\z@
  \mathrel{\mathop{\hbox to\dimen@{\rightarrowfill}}\limits^{#1}_{#2}}%
 \else
  \mathrel{\mathop{\hbox to\dimen@{\rightarrowfill}}\limits^{#1}}%
 \fi}
\def\West#1#2{\setboxz@h{$\m@th\ssize\;\;{#1}\;$}%
 \setbox\tw@\hbox{$\m@th\ssize\;\;{#2}\;$}\setbox4=\hbox{$\m@th#2$}%
 \dimen@\minaw@
 \ifdim\wdz@>\dimen@\dimen@\wdz@\fi\ifdim\wd\tw@>\dimen@\dimen@\wd\tw@\fi
 \ifdim\wd4 >\z@
  \mathrel{\mathop{\hbox to\dimen@{\leftarrowfill}}\limits^{#1}_{#2}}%
 \else
  \mathrel{\mathop{\hbox to\dimen@{\leftarrowfill}}\limits^{#1}}%
 \fi}
\font\arrow@i=lams1
\font\arrow@ii=lams2
\font\arrow@iii=lams3
\font\arrow@iv=lams4
\font\arrow@v=lams5
\newdimen\standardcgap
\standardcgap40\p@
\newdimen\hunit
\hunit\tw@\p@
\newdimen\standardrgap
\standardrgap32\p@
\newdimen\vunit
\vunit1.6\p@
\def\Cgaps#1{\RIfM@
 \standardcgap#1\standardcgap\relax\hunit#1\hunit\relax
 \else\nonmatherr@\Cgaps\fi}
\def\Rgaps#1{\RIfM@
 \standardrgap#1\standardrgap\relax\vunit#1\vunit\relax
 \else\nonmatherr@\Rgaps\fi}
\newdimen\getdim@
\def\getcgap@#1{\ifcase#1\or\getdim@\z@\else\getdim@\standardcgap\fi}
\def\getrgap@#1{\ifcase#1\getdim@\z@\else\getdim@\standardrgap\fi}
\def\cgaps{\RIfM@\expandafter\cgaps@\else\expandafter\nonmatherr@
 \expandafter\cgaps\fi}
\def\cgaps@{\ifnum\catcode`\;=\active\expandafter\cgapsA@\else
 \expandafter\cgapsO@\fi}
\def\cgapsO@#1{\toks@{\ifcase\i@\or\getdim@=\z@}%
 \gaps@@\standardcgap#1;\gaps@@\gaps@@
 \edef\next@{\the\toks@\noexpand\else\noexpand\getdim@\noexpand\standardcgap
  \noexpand\fi}%
 \toks@=\expandafter{\next@}%
 \edef\getcgap@##1{\i@##1\relax\the\toks@}\toks@{}}
{\catcode`\;=\active
 \gdef\cgapsA@#1{\toks@{\ifcase\i@\or\getdim@=\z@}%
 \gaps@@\standardcgap#1;\gaps@@\gaps@@
 \edef\next@{\the\toks@\noexpand\else\noexpand\getdim@\noexpand\standardcgap
  \noexpand\fi}%
 \toks@=\expandafter{\next@}%
 \edef\getcgap@##1{\i@##1\relax\the\toks@}\toks@{}}
}
\def\Gaps@@{\gaps@@}
\def\gaps@@#1#2;#3{\mgaps@#1#2\mgaps@
 \edef\next@{\the\toks@\noexpand\or\noexpand\getdim@
  \noexpand#1\the\mgapstoks@@}%
 \toks@\expandafter{\next@}%
 \DN@{#3}%
 \ifx\next@\Gaps@@\def\next@##1\gaps@@{}\else
  \def\next@{\gaps@@#1#3}\fi\next@}
{\catcode`\;=\active
 \gdef\rgaps#1{\RIfM@{\ifnum\catcode`\;=\active\def;{\string;}\fi
   \xdef\Next@{\noexpand\rgaps@{#1}}}%
  \Next@\edef\getrgap@##1{\i@##1\relax\the\toks@}\toks@{}\else
  \nonmatherr@\rgaps\fi}
}
\def\rgaps@#1{\toks@{\ifcase\i@\getdim@=\z@}%
 \gaps@@\standardrgap#1;\gaps@@\gaps@@
 \edef\next@{\the\toks@\noexpand\else\noexpand\getdim@\noexpand\standardrgap
  \noexpand\fi}%
 \toks@=\expandafter{\next@}}
\newbox\ZER@
\def\mgaps@#1{\let\mgapsnext@#1\FNSS@\mgaps@@}
\def\mgaps@@{\ifx\next\w\expandafter\mgaps@@@\else
 \expandafter\mgaps@@@@\fi}
\newtoks\mgapstoks@@
\def\mgaps@@@@#1\mgaps@{\getdim@\mgapsnext@\getdim@#1\getdim@
 \edef\next@{\noexpand\getdim@\the\getdim@}%
 \mgapstoks@@\expandafter{\next@}}
\def\mgaps@@@\w#1#2\mgaps@{\mgaps@@@@#2\mgaps@
 \setbox\ZER@\hbox{$\m@th\hskip15\p@\tsize@#1$}%
 \dimen@\wd\ZER@
 \ifdim\dimen@>\getdim@\getdim@\dimen@\fi
 \edef\next@{\noexpand\getdim@\the\getdim@}%
 \mgapstoks@@\expandafter{\next@}}
\def\changewidth#1#2{\setbox\ZER@{$\m@th#2}%
 \hbox to\wd\ZER@{\hss$\m@th#1$\hss}}
\atdef@({\FN@\ARROW@}
\def\ARROW@{\ifx\next)\let\next@\OPTIONS@\else
 \DN@{\csname\string @(\endcsname}\fi\next@}
\newif\ifoptions@
\def\OPTIONS@){\ifoptions@\let\next@\relax\else
 \DN@{\global\options@true\begingroup\optioncodes@}\fi\next@}
\newif\ifN@
\newif\ifE@
\newif\ifNESW@
\newif\ifH@
\newif\ifV@
\newif\ifHshort@
\expandafter\def\csname\string @(\endcsname #1,#2){%
 \ifoptions@\expandafter\endgroup\fi
 \N@false\E@false\H@false\V@false\Hshort@false
 \ifnum#1>\z@\E@true\fi
 \ifnum#1=\z@\V@true\global\tX@false\global\tY@false\global\a@false\fi
 \ifnum#2>\z@\N@true\fi
 \ifnum#2=\z@\H@true\global\tX@false\global\tY@false\global\a@false
  \ifshort@\Hshort@true\fi\fi
 \NESW@false
 \ifN@\ifE@\NESW@true\fi\else\ifE@\else\NESW@true\fi\fi
 \arrow@{#1}{#2}%
 \global\options@false
 \global\scount@\z@\global\tcount@\z@\global\arrcount@\z@
 \global\s@false\global\sxdimen@\z@\global\sydimen@\z@
 \global\tX@false\global\tXdimen@i\z@\global\tXdimen@ii\z@
 \global\tY@false\global\tYdimen@i\z@\global\tYdimen@ii\z@
 \global\a@false\global\exacount@\z@
 \global\x@false\global\xdimen@\z@
 \global\X@false\global\Xdimen@\z@
 \global\y@false\global\ydimen@\z@
 \global\Y@false\global\Ydimen@\z@
 \global\p@false\global\pdimen@\z@
 \global\label@ifalse\global\label@iifalse
 \global\dl@ifalse\global\ldimen@i\z@
 \global\dl@iifalse\global\ldimen@ii\z@
 \global\short@false\global\unshort@false}
\newif\iflabel@i
\newif\iflabel@ii
\newcount\scount@
\newcount\tcount@
\newcount\arrcount@
\newif\ifs@
\newdimen\sxdimen@
\newdimen\sydimen@
\newif\iftX@
\newdimen\tXdimen@i
\newdimen\tXdimen@ii
\newif\iftY@
\newdimen\tYdimen@i
\newdimen\tYdimen@ii
\newif\ifa@
\newcount\exacount@
\newif\ifx@
\newdimen\xdimen@
\newif\ifX@
\newdimen\Xdimen@
\newif\ify@
\newdimen\ydimen@
\newif\ifY@
\newdimen\Ydimen@
\newif\ifp@
\newdimen\pdimen@
\newif\ifdl@i
\newif\ifdl@ii
\newdimen\ldimen@i
\newdimen\ldimen@ii
\newif\ifshort@
\newif\ifunshort@
\def\zero@#1{\ifnum\scount@=\z@
 \if#1e\global\scount@\m@ne\else
 \if#1t\global\scount@\tw@\else
 \if#1h\global\scount@\thr@@\else
 \if#1'\global\scount@6 \else
 \if#1`\global\scount@7 \else
 \if#1(\global\scount@8 \else
 \if#1)\global\scount@9 \else
 \if#1s\global\scount@12 \else
 \if#1H\global\scount@13 \else
 \Err@{\Invalid@@ option \string\0}\fi\fi\fi\fi\fi\fi\fi\fi\fi
 \fi}
\def\one@#1{\ifnum\tcount@=\z@
 \if#1e\global\tcount@\m@ne\else
 \if#1h\global\tcount@\tw@\else
 \if#1t\global\tcount@\thr@@\else
 \if#1'\global\tcount@4 \else
 \if#1`\global\tcount@5 \else
 \if#1(\global\tcount@\ten@ \else
 \if#1)\global\tcount@11 \else
 \if#1s\global\tcount@12 \else
 \if#1H\global\tcount@13 \else
 \Err@{\Invalid@@ option \string\1}\fi\fi\fi\fi\fi\fi\fi\fi\fi
 \fi}
\def\a@#1{\ifnum\arrcount@=\z@
 \if#10\global\arrcount@\m@ne\else
 \if#1+\global\arrcount@\@ne\else
 \if#1-\global\arrcount@\tw@\else
 \if#1=\global\arrcount@\thr@@\else
 \Err@{\Invalid@@ option \string\a}\fi\fi\fi\fi
 \fi}
\def\ds@{\ifnum\catcode`\;=\active\expandafter\dsA@\else
 \expandafter\dsO@\fi}
\def\dsO@(#1;#2){\ds@@{#1}{#2}}
\def\ds@@#1#2{\ifs@\else
 \global\s@true
 \global\sxdimen@\hunit\global\sxdimen@#1\sxdimen@\relax
 \global\sydimen@\vunit\global\sydimen@#2\sydimen@\relax
 \fi}
\def\dtX@{\ifnum\catcode`\;=\active\expandafter\dtXA@\else
 \expandafter\dtXO@\fi}
\def\dtXO@(#1;#2){\dtX@@{#1}{#2}}
\def\dtX@@#1#2{\iftX@\else
 \global\tX@true
 \global\tXdimen@i\hunit\global\tXdimen@i#1\tXdimen@i\relax
 \global\tXdimen@ii\vunit\global\tXdimen@ii#2\tXdimen@ii\relax
 \fi}
\def\dtY@{\ifnum\catcode`\;=\active\expandafter\dtYA@\else
 \expandafter\dtYO@\fi}
\def\dtYO@(#1;#2){\dtY@@{#1}{#2}}
\def\dtY@@#1#2{\iftY@\else
 \global\tY@true
 \global\tYdimen@i\hunit\global\tYdimen@i#1\tYdimen@i\relax
 \global\tYdimen@ii\vunit\global\tYdimen@ii#2\tYdimen@ii\relax
 \fi}
{\catcode`\;=\active
 \gdef\dsA@(#1;#2){\ds@@{#1}{#2}}
 \gdef\dtXA@(#1;#2){\dtX@@{#1}{#2}}
 \gdef\dtYA@(#1;#2){\dtY@@{#1}{#2}}
}
\def\da@#1{\ifa@\else\global\a@true\global\exacount@#1\relax\fi}
\def\dx@#1{\ifx@\else
 \global\x@true
 \global\xdimen@\hunit\global\xdimen@#1\xdimen@\relax
 \fi}
\def\dX@#1{\ifX@\else
 \global\X@true
 \global\Xdimen@\hunit\global\Xdimen@#1\Xdimen@\relax
 \fi}
\def\dy@#1{\ify@\else
 \global\y@true
 \global\ydimen@\vunit\global\ydimen@#1\ydimen@\relax
 \fi}
\def\dY@#1{\ifY@\else
 \global\Y@true
 \global\Ydimen@\vunit\global\Ydimen@#1\Ydimen@\relax
 \fi}
\def\p@@#1{\ifp@\else
 \global\p@true
 \global\pdimen@\hunit\global\divide\pdimen@\tw@
 \global\pdimen@#1\pdimen@\relax
 \fi}
\def\L@#1{\iflabel@i\else
 \global\label@itrue\gdef\label@i{#1}%
 \fi}
\def\l@#1{\iflabel@ii\else
 \global\label@iitrue\gdef\label@ii{#1}%
 \fi}
\def\dL@#1{\ifdl@i\else
 \global\dl@itrue\global\ldimen@i\hunit\global\ldimen@i#1\ldimen@i\relax
 \fi}
\def\dl@#1{\ifdl@ii\else
 \global\dl@iitrue\global\ldimen@ii\hunit\global\ldimen@ii#1\ldimen@ii\relax
 \fi}
\def\s@{\ifunshort@\else\global\short@true\fi}
\def\uns@{\ifshort@\else\global\unshort@true\global\short@false\fi}
\def\optioncodes@{\let\0\zero@\let\1\one@\let\a\a@\let\ds\ds@\let\dtX\dtX@
 \let\dtY\dtY@\let\da\da@\let\dx\dx@\let\dX\dX@\let\dY\dY@\let\dy\dy@
 \let\p\p@@\let\L\L@\let\l\l@\let\dL\dL@\let\dl\dl@\let\s\s@\let\uns\uns@}
\def\slopes@{\\161\\152\\143\\134\\255\\126\\357\\238\\349\\45{10}\\56{11}%
 \\11{12}\\65{13}\\54{14}\\43{15}\\32{16}\\53{17}\\21{18}\\52{19}\\31{20}%
 \\41{21}\\51{22}\\61{23}}
\newcount\tan@i
\newcount\tan@ip
\newcount\tan@ii
\newcount\tan@iip
\newdimen\slope@i
\newdimen\slope@ip
\newdimen\slope@ii
\newdimen\slope@iip
\newcount\angcount@
\newcount\extracount@
\def\slope@{{\slope@i\secondy@\advance\slope@i-\firsty@
 \ifN@\else\multiply\slope@i\m@ne\fi
 \slope@ii\secondx@\advance\slope@ii-\firstx@
 \ifE@\else\multiply\slope@ii\m@ne\fi
 \ifdim\slope@ii<\z@
  \global\tan@i6 \global\tan@ii\@ne\global\angcount@23
 \else
  \dimen@\slope@i\multiply\dimen@6
  \ifdim\dimen@<\slope@ii
   \global\tan@i\@ne\global\tan@ii6 \global\angcount@\@ne
  \else
   \dimen@\slope@ii\multiply\dimen@6
   \ifdim\dimen@<\slope@i
    \global\tan@i6 \global\tan@ii\@ne\global\angcount@23
   \else
    \global\tan@ip\z@\global\tan@iip\@ne
    \def\\##1##2##3{\global\angcount@##3\relax
     \slope@ip\slope@i\slope@iip\slope@ii
     \multiply\slope@iip##1\relax\multiply\slope@ip##2\relax
     \ifdim\slope@iip<\slope@ip
      \global\tan@ip##1\relax\global\tan@iip##2\relax
     \else
      \global\tan@i##1\relax\global\tan@ii##2\relax
      \def\\####1####2####3{}%
     \fi}%
    \slopes@
    \slope@i\secondy@\advance\slope@i-\firsty@
    \ifN@\else\multiply\slope@i\m@ne\fi
    \multiply\slope@i\tan@ii\multiply\slope@i\tan@iip\multiply\slope@i\tw@
    \count@\tan@i\multiply\count@\tan@iip
    \extracount@\tan@ip\multiply\extracount@\tan@ii
    \advance\count@\extracount@
    \slope@ii\secondx@\advance\slope@ii-\firstx@
    \ifE@\else\multiply\slope@ii\m@ne\fi
    \multiply\slope@ii\count@
    \ifdim\slope@i<\slope@ii
     \global\tan@i\tan@ip\global\tan@ii\tan@iip
     \global\advance\angcount@\m@ne
    \fi
   \fi
  \fi
 \fi}%
}
\def\slope@a#1{{\def\\##1##2##3{\ifnum##3=#1\global\tan@i##1\relax
 \global\tan@ii##2\relax\fi}\slopes@}}
\newcount\i@
\newcount\j@
\newcount\colcount@
\newcount\Colcount@
\newcount\tcolcount@
\newdimen\rowht@
\newdimen\rowdp@
\newcount\rowcount@
\newcount\Rowcount@
\newcount\maxcolrow@
\newtoks\colwidthtoks@
\newtoks\Rowheighttoks@
\newtoks\Rowdepthtoks@
\newtoks\widthtoks@
\newtoks\Widthtoks@
\newtoks\heighttoks@
\newtoks\Heighttoks@
\newtoks\depthtoks@
\newtoks\Depthtoks@
\newif\iffirstCDcr@
\def\dotoks@i{%
 \global\widthtoks@\expandafter{\the\widthtoks@\else\getdim@\z@\fi}%
 \global\heighttoks@\expandafter{\the\heighttoks@\else\getdim@\z@\fi}%
 \global\depthtoks@\expandafter{\the\depthtoks@\else\getdim@\z@\fi}}
\def\dotoks@ii{%
 \global\widthtoks@{\ifcase\j@}%
 \global\heighttoks@{\ifcase\j@}%
 \global\depthtoks@{\ifcase\j@}}
\def\preCD@#1\endCD{\setbox\ZER@
 \vbox{%
  \def\arrow@##1##2{{}}%
  \global\rowcount@\m@ne\global\colcount@\z@\global\Colcount@\z@
  \global\firstCDcr@true\toks@{}%
  \global\widthtoks@{\ifcase\j@}%
  \global\Widthtoks@{\ifcase\i@}%
  \global\heighttoks@{\ifcase\j@}%
  \global\Heighttoks@{\ifcase\i@}%
  \global\depthtoks@{\ifcase\j@}%
  \global\Depthtoks@{\ifcase\i@}%
  \global\Rowheighttoks@{\ifcase\i@}%
  \global\Rowdepthtoks@{\ifcase\i@}%
  \Let@
  \everycr{%
   \noalign{%
    \global\advance\rowcount@\@ne
    \ifnum\colcount@<\Colcount@
    \else
     \global\Colcount@\colcount@\global\maxcolrow@\rowcount@
    \fi
    \global\colcount@\z@
    \iffirstCDcr@
     \global\firstCDcr@false
    \else
     \edef\next@{\the\Rowheighttoks@\noexpand\or\noexpand\getdim@\the\rowht@}%
      \global\Rowheighttoks@\expandafter{\next@}%
     \edef\next@{\the\Rowdepthtoks@\noexpand\or\noexpand\getdim@\the\rowdp@}%
      \global\Rowdepthtoks@\expandafter{\next@}%
     \global\rowht@\z@\global\rowdp@\z@
     \dotoks@i
     \edef\next@{\the\Widthtoks@\noexpand\or\the\widthtoks@}%
      \global\Widthtoks@\expandafter{\next@}%
     \edef\next@{\the\Heighttoks@\noexpand\or\the\heighttoks@}%
      \global\Heighttoks@\expandafter{\next@}%
     \edef\next@{\the\Depthtoks@\noexpand\or\the\depthtoks@}%
      \global\Depthtoks@\expandafter{\next@}%
     \dotoks@ii
    \fi}}%
  \tabskip\z@
  \halign{&\setbox\ZER@\hbox{\vrule\height\ten@\p@\width\z@\depth\z@     
   $\m@th\displaystyle{##}$}\copy\ZER@
   \ifdim\ht\ZER@>\rowht@\global\rowht@\ht\ZER@\fi
   \ifdim\dp\ZER@>\rowdp@\global\rowdp@\dp\ZER@\fi
   \global\advance\colcount@\@ne
   \edef\next@{\the\widthtoks@\noexpand\or\noexpand\getdim@\the\wd\ZER@}%
    \global\widthtoks@\expandafter{\next@}%
   \edef\next@{\the\heighttoks@\noexpand\or\noexpand\getdim@\the\ht\ZER@}%
    \global\heighttoks@\expandafter{\next@}%
   \edef\next@{\the\depthtoks@\noexpand\or\noexpand\getdim@\the\dp\ZER@}%
    \global\depthtoks@\expandafter{\next@}%
   \cr#1\crcr}}%
 \Rowcount@\rowcount@
 \global\Widthtoks@\expandafter{\the\Widthtoks@\fi\relax}%
 \edef\Width@##1##2{\i@##1\relax\j@##2\relax\the\Widthtoks@}%
 \global\Heighttoks@\expandafter{\the\Heighttoks@\fi\relax}%
 \edef\Height@##1##2{\i@##1\relax\j@##2\relax\the\Heighttoks@}%
 \global\Depthtoks@\expandafter{\the\Depthtoks@\fi\relax}%
 \edef\Depth@##1##2{\i@##1\relax\j@##2\relax\the\Depthtoks@}%
 \edef\next@{\the\Rowheighttoks@\noexpand\fi\relax}%
 \global\Rowheighttoks@\expandafter{\next@}%
 \edef\Rowheight@##1{\i@##1\relax\the\Rowheighttoks@}%
 \edef\next@{\the\Rowdepthtoks@\noexpand\fi\relax}%
 \global\Rowdepthtoks@\expandafter{\next@}%
 \edef\Rowdepth@##1{\i@##1\relax\the\Rowdepthtoks@}%
 \global\colwidthtoks@{\fi}%
 \setbox\ZER@\vbox{%
  \unvbox\ZER@
  \count@\rowcount@
  \loop
   \unskip\unpenalty
   \setbox\ZER@\lastbox
   \ifnum\count@>\maxcolrow@\advance\count@\m@ne
   \repeat
  \hbox{%
   \unhbox\ZER@
   \count@\z@
   \loop
    \unskip
    \setbox\ZER@\lastbox
    \edef\next@{\noexpand\or\noexpand\getdim@\the\wd\ZER@\the\colwidthtoks@}%
     \global\colwidthtoks@\expandafter{\next@}%
    \advance\count@\@ne
    \ifnum\count@<\Colcount@
    \repeat}}%
 \edef\next@{\noexpand\ifcase\noexpand\i@\the\colwidthtoks@}%
  \global\colwidthtoks@\expandafter{\next@}%
 \edef\Colwidth@##1{\i@##1\relax\the\colwidthtoks@}%
 \global\colwidthtoks@{}\global\Rowheighttoks@{}\global\Rowdepthtoks@{}%
 \global\widthtoks@{}\global\Widthtoks@{}\global\heighttoks@{}%
 \global\Heighttoks@{}\global\depthtoks@{}\global\Depthtoks@{}%
}
\newcount\xoff@
\newcount\yoff@
\newcount\endcount@
\newcount\rcount@
\newdimen\firstx@
\newdimen\firsty@
\newdimen\secondx@
\newdimen\secondy@
\newdimen\tocenter@
\newdimen\charht@
\newdimen\charwd@
\def\outside@{\Err@{This arrow points outside the \string\CD}}
\newif\ifsvertex@
\newif\iftvertex@
\def\arrow@#1#2{\global\xoff@#1\relax\global\yoff@#2\relax
 \count@\rowcount@\advance\count@-\yoff@
 \ifnum\count@<\@ne\outside@\else\ifnum\count@>\Rowcount@\outside@\fi\fi
 \count@\colcount@\advance\count@\xoff@
 \ifnum\count@<\@ne\outside@\else\ifnum\count@>\Colcount@\outside@\fi\fi
 \tcolcount@\colcount@\advance\tcolcount@\xoff@
 \Width@\rowcount@\colcount@\divide\getdim@\tw@\tocenter@-\getdim@
 \ifdim\getdim@=\z@
  \firstx@\z@\firsty@\mathaxis@\svertex@true
 \else
  \svertex@false
  \ifHshort@
   \Colwidth@\colcount@\divide\getdim@\tw@
   \ifE@ \firstx@\getdim@ \else \firstx@-\getdim@ \fi
  \else
   \ifE@ \firstx@\getdim@ \else \firstx@-\getdim@ \fi
  \fi
  \ifE@
   \ifH@ \advance\firstx@\thr@@\p@ \else \advance\firstx@-\thr@@\p@ \fi  
  \else
   \ifH@ \advance\firstx@-\thr@@\p@ \else \advance\firstx@\thr@@\p@ \fi  
  \fi
  \ifN@
   \Height@\rowcount@\colcount@ \firsty@\getdim@                         
   \ifV@ \advance\firsty@\thr@@\p@ \fi                                   
  \else
   \ifV@
    \Depth@\rowcount@\colcount@ \firsty@-\getdim@                        
    \advance\firsty@-\thr@@\p@                                           
   \else
    \firsty@\z@                                                          
   \fi
  \fi
 \fi
 \ifV@
 \else
  \Colwidth@\colcount@\divide\getdim@\tw@
  \ifE@\secondx@\getdim@\else\secondx@-\getdim@\fi
  \ifE@\else\getcgap@\colcount@\advance\secondx@-\getdim@\fi
  \endcount@\colcount@\advance\endcount@\xoff@
  \count@\colcount@
  \ifE@
   \advance\count@\@ne
   \loop
    \ifnum\count@<\endcount@
    \Colwidth@\count@\advance\secondx@\getdim@
    \getcgap@\count@\advance\secondx@\getdim@
    \advance\count@\@ne
    \repeat
  \else
   \advance\count@\m@ne
   \loop
    \ifnum\count@>\endcount@
    \Colwidth@\count@\advance\secondx@-\getdim@
    \getcgap@\count@\advance\secondx@-\getdim@
    \advance\count@\m@ne
    \repeat
  \fi
  \Colwidth@\count@\divide\getdim@\tw@
  \ifHshort@
  \else
   \ifE@\advance\secondx@\getdim@\else\advance\secondx@-\getdim@\fi
  \fi
  \ifE@\getcgap@\count@\advance\secondx@\getdim@\fi
  \rcount@\rowcount@\advance\rcount@-\yoff@
  \Width@\rcount@\count@\divide\getdim@\tw@
  \tvertex@false
  \ifH@\ifdim\getdim@=\z@\tvertex@true\Hshort@false\fi\fi
  \ifHshort@
  \else
   \ifE@\advance\secondx@-\getdim@\else\advance\secondx@\getdim@\fi
  \fi
  \iftvertex@
   \advance\secondx@.4\p@
  \else
   \ifE@\advance\secondx@-\thr@@\p@\else\advance\secondx@\thr@@\p@\fi    
  \fi
 \fi
 \ifH@
 \else
  \ifN@
   \Rowheight@\rowcount@\secondy@\getdim@
  \else
   \Rowdepth@\rowcount@\secondy@-\getdim@
   \getrgap@\rowcount@\advance\secondy@-\getdim@
  \fi
  \endcount@\rowcount@\advance\endcount@-\yoff@
  \count@\rowcount@
  \ifN@
   \advance\count@\m@ne
   \loop
    \ifnum\count@>\endcount@
    \Rowheight@\count@\advance\secondy@\getdim@
    \Rowdepth@\count@\advance\secondy@\getdim@
    \getrgap@\count@\advance\secondy@\getdim@
    \advance\count@\m@ne
    \repeat
  \else
   \advance\count@\@ne
   \loop
    \ifnum\count@<\endcount@
    \Rowheight@\count@\advance\secondy@-\getdim@
    \Rowdepth@\count@\advance\secondy@-\getdim@
    \getrgap@\count@\advance\secondy@-\getdim@
    \advance\count@\@ne
    \repeat
  \fi
  \tvertex@false
  \ifV@\Width@\count@\colcount@\ifdim\getdim@=\z@\tvertex@true\fi\fi
  \ifN@
   \getrgap@\count@\advance\secondy@\getdim@
   \Rowdepth@\count@\advance\secondy@\getdim@
   \iftvertex@
    \advance\secondy@\mathaxis@
   \else
    \Depth@\count@\tcolcount@\advance\secondy@-\getdim@
    \advance\secondy@-\thr@@\p@                                          
   \fi
  \else
   \Rowheight@\count@\advance\secondy@-\getdim@
   \iftvertex@
    \advance\secondy@\mathaxis@
   \else
    \Height@\count@\tcolcount@\advance\secondy@\getdim@
    \advance\secondy@\thr@@\p@                                           
   \fi
  \fi
 \fi
 \ifV@\else\advance\firstx@\sxdimen@\fi
 \ifH@\else\advance\firsty@\sydimen@\fi
 \iftX@
  \advance\secondy@\tXdimen@ii
  \advance\secondx@\tXdimen@i
  \slope@
 \else
  \iftY@
   \advance\secondy@\tYdimen@ii
   \advance\secondx@\tYdimen@i
   \slope@
   \secondy@\secondx@\advance\secondy@-\firstx@
   \ifNESW@\else\multiply\secondy@\m@ne\fi
   \multiply\secondy@\tan@i\divide\secondy@\tan@ii\advance\secondy@\firsty@
  \else
   \ifa@
    \slope@
    \ifNESW@\global\advance\angcount@\exacount@\else
     \global\advance\angcount@-\exacount@\fi
    \ifnum\angcount@>23 \global\angcount@23 \fi
    \ifnum\angcount@<\@ne\global\angcount@\@ne\fi
    \slope@a\angcount@
    \ifY@
     \advance\secondy@\Ydimen@
    \else
     \ifX@
      \advance\secondx@\Xdimen@
      \dimen@\secondx@\advance\dimen@-\firstx@
      \ifNESW@\else\multiply\dimen@\m@ne\fi
      \multiply\dimen@\tan@i\divide\dimen@\tan@ii
      \advance\dimen@\firsty@\secondy@\dimen@
     \fi
    \fi
   \else
    \ifH@\else\ifV@\else\slope@\fi\fi
   \fi
  \fi
 \fi
 \ifH@\else\ifV@\else\ifsvertex@\else
  \dimen@6\p@\multiply\dimen@\tan@ii
  \count@\tan@i\advance\count@\tan@ii\divide\dimen@\count@
  \ifE@\advance\firstx@\dimen@\else\advance\firstx@-\dimen@\fi
  \multiply\dimen@\tan@i\divide\dimen@\tan@ii
  \ifN@\advance\firsty@\dimen@\else\advance\firsty@-\dimen@\fi
 \fi\fi\fi
 \ifp@
  \ifH@\else\ifV@\else
   \getcos@\pdimen@\advance\firsty@\dimen@\advance\secondy@\dimen@
   \ifNESW@\advance\firstx@-\dimen@ii\else\advance\firstx@\dimen@ii\fi
  \fi\fi
 \fi
 \ifH@\else\ifV@\else
  \ifnum\tan@i>\tan@ii
   \charht@\ten@\p@\charwd@\ten@\p@
   \multiply\charwd@\tan@ii\divide\charwd@\tan@i
  \else
   \charwd@\ten@\p@\charht@\ten@\p@
   \divide\charht@\tan@ii\multiply\charht@\tan@i
  \fi
  \ifnum\tcount@=\thr@@
   \ifN@\advance\secondy@-.3\charht@\else\advance\secondy@.3\charht@\fi
  \fi
  \ifnum\scount@=\tw@
   \ifE@\advance\firstx@.3\charht@\else\advance\firstx@-.3\charht@\fi
  \fi
  \ifnum\tcount@=12
   \ifN@\advance\secondy@-\charht@\else\advance\secondy@\charht@\fi
  \fi
  \iftY@
  \else
   \ifa@
    \ifX@
    \else
     \secondx@\secondy@\advance\secondx@-\firsty@
     \ifNESW@\else\multiply\secondx@\m@ne\fi
     \multiply\secondx@\tan@ii\divide\secondx@\tan@i
     \advance\secondx@\firstx@
    \fi
   \fi
  \fi
 \fi\fi
 \ifH@\harrow@\else\ifV@\varrow@\else\arrow@@\fi\fi}
\newdimen\mathaxis@
\mathaxis@90\p@\divide\mathaxis@36
\def\harrow@b{\ifE@\hskip\tocenter@\hskip\firstx@\fi}
\def\harrow@bb{\ifE@\hskip\xdimen@\else\hskip\Xdimen@\fi}
\def\harrow@e{\ifE@\else\hskip-\firstx@\hskip-\tocenter@\fi}
\def\harrow@ee{\ifE@\hskip-\Xdimen@\else\hskip-\xdimen@\fi}
\def\harrow@{\dimen@\secondx@\advance\dimen@-\firstx@
 \ifE@\let\next@\rlap\else\multiply\dimen@\m@ne\let\next@\llap\fi
 \next@{%
  \harrow@b
  \smash{\raise\pdimen@\hbox to\dimen@
   {\harrow@bb\arrow@ii
    \ifnum\arrcount@=\m@ne\else\ifnum\arrcount@=\thr@@\else
     \ifE@
      \ifnum\scount@=\m@ne
      \else
       \ifcase\scount@\or\or\char118 \or\char117 \or\or\or\char119 \or
       \char120 \or\char121 \or\char122 \or\or\or\arrow@i\char125 \or
       \char117 \hskip\thr@@\p@\char117 \hskip-\thr@@\p@\fi
      \fi
     \else
      \ifnum\tcount@=\m@ne
      \else
       \ifcase\tcount@\char117 \or\or\char117 \or\char118 \or\char119 \or
       \char120 \or\or\or\or\or\char121 \or\char122 \or\arrow@i\char125
       \or\char117 \hskip\thr@@\p@\char117 \hskip-\thr@@\p@\fi
      \fi
     \fi
    \fi\fi
    \dimen@\mathaxis@\advance\dimen@.2\p@
    \dimen@ii\mathaxis@\advance\dimen@ii-.2\p@
    \ifnum\arrcount@=\m@ne
     \let\leads@\null
    \else
     \ifcase\arrcount@
      \def\leads@{\hrule\height\dimen@\depth-\dimen@ii}\or
      \def\leads@{\hrule\height\dimen@\depth-\dimen@ii}\or
      \def\leads@{\hbox to\ten@\p@{%
       \leaders\hrule\height\dimen@\depth-\dimen@ii\hfil
       \hfil
      \leaders\hrule\height\dimen@\depth-\dimen@ii\hskip\z@ plus2fil\relax
       \hfil
       \leaders\hrule\height\dimen@\depth-\dimen@ii\hfil}}\or
     \def\leads@{\hbox{\hbox to\ten@\p@{\dimen@\mathaxis@\advance\dimen@1.2\p@
       \dimen@ii\dimen@\advance\dimen@ii-.4\p@
       \leaders\hrule\height\dimen@\depth-\dimen@ii\hfil}%
       \kern-\ten@\p@
       \hbox to\ten@\p@{\dimen@\mathaxis@\advance\dimen@-1.2\p@
       \dimen@ii\dimen@\advance\dimen@ii-.4\p@
       \leaders\hrule\height\dimen@\depth-\dimen@ii\hfil}}}\fi
    \fi
    \cleaders\leads@\hfil
    \ifnum\arrcount@=\m@ne\else\ifnum\arrcount@=\thr@@\else
     \arrow@i
     \ifE@
      \ifnum\tcount@=\m@ne
      \else
       \ifcase\tcount@\char119 \or\or\char119 \or\char120 \or\char121 \or
       \char122 \or \or\or\or\or\char123 \or\char124 \or
       \char125 \or\char119 \hskip-\thr@@\p@\char119 \hskip\thr@@\p@\fi
      \fi
     \else
      \ifcase\scount@\or\or\char120 \or\char119 \or\or\or\char121 \or\char122
      \or\char123 \or\char124 \or\or\or\char125 \or
      \char119 \hskip-\thr@@\p@\char119 \hskip\thr@@\p@\fi
     \fi
    \fi\fi
    \harrow@ee}}%
  \harrow@e}%
 \iflabel@i
  \dimen@ii\z@\setbox\ZER@\hbox{$\m@th\tsize@@\label@i$}%
  \ifnum\arrcount@=\m@ne
  \else
   \advance\dimen@ii\mathaxis@
   \advance\dimen@ii\dp\ZER@\advance\dimen@ii\tw@\p@
   \ifnum\arrcount@=\thr@@\advance\dimen@ii\tw@\p@\fi
  \fi
  \advance\dimen@ii\pdimen@
  \next@{\harrow@b\smash{\raise\dimen@ii\hbox to\dimen@
   {\harrow@bb\hskip\tw@\ldimen@i\hfil\box\ZER@\hfil\harrow@ee}}\harrow@e}%
 \fi
 \iflabel@ii
  \ifnum\arrcount@=\m@ne
  \else
   \setbox\ZER@\hbox{$\m@th\tsize@\label@ii$}%
   \dimen@ii-\ht\ZER@\advance\dimen@ii-\tw@\p@
   \ifnum\arrcount@=\thr@@\advance\dimen@ii-\tw@\p@\fi
   \advance\dimen@ii\mathaxis@\advance\dimen@ii\pdimen@
   \next@{\harrow@b\smash{\raise\dimen@ii\hbox to\dimen@
    {\harrow@bb\hskip\tw@\ldimen@ii\hfil\box\ZER@\hfil\harrow@ee}}\harrow@e}%
  \fi
 \fi}
\let\tsize@\tsize
\def\tsizeCDlabels{\let\tsize@\tsize}
\def\ssizeCDlabels{\let\tsize@\ssize}
\def\tsize@@{\ifnum\arrcount@=\m@ne\else\tsize@\fi}
\def\varrow@{\dimen@\secondy@\advance\dimen@-\firsty@
 \ifN@\else\multiply\dimen@\m@ne\fi
 \setbox\ZER@\vbox to\dimen@
  {\ifN@\vskip-\Ydimen@\else\vskip\ydimen@\fi
   \ifnum\arrcount@=\m@ne\else\ifnum\arrcount@=\thr@@\else
    \hbox{\arrow@iii
     \ifN@
      \ifnum\tcount@=\m@ne
      \else
       \ifcase\tcount@\char117 \or\or\char117 \or\char118 \or\char119 \or
       \char120 \or\or\or\or\or\char121 \or\char122 \or\char123 \or
       \vbox{\hbox{\char117}\nointerlineskip\vskip\thr@@\p@
       \hbox{\char117}\vskip-\thr@@\p@}\fi
      \fi
     \else
      \ifcase\scount@\or\or\char118 \or\char117 \or\or\or\char119 \or
      \char120 \or\char121 \or\char122 \or\or\or\char123 \or
      \vbox{\hbox{\char117}\nointerlineskip\vskip\thr@@\p@
      \hbox{\char117}\vskip-\thr@@\p@}\fi
     \fi}%
    \nointerlineskip
   \fi\fi
   \ifnum\arrcount@=\m@ne
    \let\leads@\null
   \else
    \ifcase\arrcount@\let\leads@\vrule\or\let\leads@\vrule\or
    \def\leads@{\vbox to\ten@\p@{%
     \hrule\height1.67\p@\depth\z@\width.4\p@
     \vfil
     \hrule\height3.33\p@\depth\z@\width.4\p@
     \vfil
     \hrule\height1.67\p@\depth\z@\width.4\p@}}\or
    \def\leads@{\hbox{\vrule\height\p@\hskip\tw@\p@\vrule}}\fi
   \fi
  \cleaders\leads@\vfill\nointerlineskip
   \ifnum\arrcount@=\m@ne\else\ifnum\arrcount@=\thr@@\else
    \hbox{\arrow@iv
     \ifN@
      \ifcase\scount@\or\or\char118 \or\char117 \or\or\or\char119 \or
      \char120 \or\char121 \or\char122 \or\or\or\arrow@iii\char123 \or
      \vbox{\hbox{\char117}\nointerlineskip\vskip-\thr@@\p@
      \hbox{\char117}\vskip\thr@@\p@}\fi
     \else
      \ifnum\tcount@=\m@ne
      \else
       \ifcase\tcount@\char117 \or\or\char117 \or\char118 \or\char119 \or
       \char120 \or\or\or\or\or\char121 \or\char122 \or\arrow@iii\char123 \or
       \vbox{\hbox{\char117}\nointerlineskip\vskip-\thr@@\p@
       \hbox{\char117}\vskip\thr@@\p@}\fi
      \fi
     \fi}%
   \fi\fi
   \ifN@\vskip\ydimen@\else\vskip-\Ydimen@\fi}%
 \ifN@
  \dimen@ii\firsty@
 \else
  \dimen@ii-\firsty@\advance\dimen@ii\ht\ZER@\multiply\dimen@ii\m@ne
 \fi
 \rlap{\smash{\hskip\tocenter@\hskip\pdimen@\raise\dimen@ii\box\ZER@}}%
 \iflabel@i
  \setbox\ZER@\vbox to\dimen@{\vfil
   \hbox{$\m@th\tsize@@\label@i$}\vskip\tw@\ldimen@i\vfil}%
  \rlap{\smash{\hskip\tocenter@\hskip\pdimen@
  \ifnum\arrcount@=\m@ne\let\next@\relax\else\let\next@\llap\fi
  \next@{\raise\dimen@ii\hbox{\ifnum\arrcount@=\m@ne\hskip-.5\wd\ZER@\fi
   \box\ZER@\ifnum\arrcount@=\m@ne\else\hskip\tw@\p@\fi}}}}%
 \fi
 \iflabel@ii
  \ifnum\arrcount@=\m@ne
  \else
   \setbox\ZER@\vbox to\dimen@{\vfil
    \hbox{$\m@th\tsize@\label@ii$}\vskip\tw@\ldimen@ii\vfil}%
   \rlap{\smash{\hskip\tocenter@\hskip\pdimen@
   \rlap{\raise\dimen@ii\hbox{\ifnum\arrcount@=\thr@@\hskip4.5\p@\else
    \hskip2.5\p@\fi\box\ZER@}}}}%
  \fi
 \fi
}
\newdimen\goal@
\newdimen\shifted@
\newcount\Tcount@
\newcount\Scount@
\newbox\shaft@
\newcount\slcount@
\def\getcos@#1{%
 \ifnum\tan@i<\tan@ii
  \dimen@#1%
  \ifnum\slcount@<8 \count@9 \else \ifnum\slcount@<12 \count@8 \else
   \count@7 \fi\fi
  \multiply\dimen@\count@\divide\dimen@\ten@
  \dimen@ii\dimen@\multiply\dimen@ii\tan@i\divide\dimen@ii\tan@ii
 \else
  \dimen@ii#1%
  \count@-\slcount@\advance\count@24
  \ifnum\count@<8 \count@9 \else \ifnum\count@<12 \count@8
   \else\count@7 \fi\fi
  \multiply\dimen@ii\count@\divide\dimen@ii\ten@
  \dimen@\dimen@ii\multiply\dimen@\tan@ii\divide\dimen@\tan@i
 \fi}
\newdimen\adjust@
\def\Nnext@{\ifN@\let\next@\raise\else\let\next@\lower\fi}
\def\arrow@@{\slcount@\angcount@
 \ifNESW@
  \ifnum\angcount@<\ten@
   \let\arrowfont@\arrow@i\global\advance\angcount@\m@ne
   \global\multiply\angcount@13
  \else
   \ifnum\angcount@<19
    \let\arrowfont@\arrow@ii\global\advance\angcount@-\ten@
    \global\multiply\angcount@13
   \else
    \let\arrowfont@\arrow@iii\global\advance\angcount@-19
    \global\multiply\angcount@13
  \fi\fi
  \Tcount@\angcount@
 \else
  \ifnum\angcount@<5
   \let\arrowfont@\arrow@iii\global\advance\angcount@\m@ne
   \global\multiply\angcount@13 \global\advance\angcount@65
  \else
   \ifnum\angcount@<14
    \let\arrowfont@\arrow@iv\global\advance\angcount@-5
    \global\multiply\angcount@13
   \else
    \ifnum\angcount@<23
     \let\arrowfont@\arrow@v\global\advance\angcount@-14
     \global\multiply\angcount@13
    \else
     \let\arrowfont@\arrow@i\global\angcount@117
  \fi\fi\fi
  \ifnum\angcount@=117 \Tcount@115 \else\Tcount@\angcount@\fi
 \fi
 \Scount@\Tcount@
 \ifE@
  \ifnum\tcount@=\z@\advance\Tcount@\tw@\else\ifnum\tcount@=13
   \advance\Tcount@\tw@\else\advance\Tcount@\tcount@\fi\fi
  \ifnum\scount@=\z@\else\ifnum\scount@=13 \advance\Scount@\thr@@\else
   \advance\Scount@\scount@\fi\fi
 \else
  \ifcase\tcount@\advance\Tcount@\thr@@\or\or\advance\Tcount@\thr@@\or
  \advance\Tcount@\tw@\or\advance\Tcount@6 \or\advance\Tcount@7
  \or\or\or\or\or\advance\Tcount@8 \or\advance\Tcount@9 \or
  \advance\Tcount@12 \or\advance\Tcount@\thr@@\fi
  \ifcase\scount@\or\or\advance\Scount@\thr@@\or\advance\Scount@\tw@\or
  \or\or\advance\Scount@4 \or\advance\Scount@5 \or\advance\Scount@\ten@
  \or\advance\Scount@11 \or\or\or\advance\Scount@12 \or\advance
  \Scount@\tw@\fi
 \fi
 \ifcase\arrcount@\or\or\global\advance\angcount@\@ne\else\fi
 \ifN@\shifted@\firsty@\else\shifted@-\firsty@\fi
 \ifE@\else\advance\shifted@\charht@\fi
 \goal@\secondy@\advance\goal@-\firsty@
 \ifN@\else\multiply\goal@\m@ne\fi
 \setbox\shaft@\hbox{\arrowfont@\char\angcount@}%
 \ifnum\arrcount@=\thr@@
  \getcos@{1.5\p@}%
  \setbox\shaft@\hbox to\wd\shaft@{\arrowfont@
   \rlap{\hskip\dimen@ii
    \smash{\ifNESW@\let\next@\lower\else\let\next@\raise\fi
     \next@\dimen@\hbox{\arrowfont@\char\angcount@}}}%
   \rlap{\hskip-\dimen@ii
    \smash{\ifNESW@\let\next@\raise\else\let\next@\lower\fi
      \next@\dimen@\hbox{\arrowfont@\char\angcount@}}}\hfil}%
 \fi
 \rlap{\smash{\hskip\tocenter@\hskip\firstx@
  \ifnum\arrcount@=\m@ne
  \else
   \ifnum\arrcount@=\thr@@
   \else
    \ifnum\scount@=\m@ne
    \else
     \ifnum\scount@=\z@
     \else
      \setbox\ZER@\hbox{\ifnum\angcount@=117 \arrow@v\else\arrowfont@\fi
       \char\Scount@}%
      \ifNESW@
       \ifnum\scount@=\tw@
        \dimen@\shifted@\advance\dimen@-\charht@
        \ifN@\hskip-\wd\ZER@\fi
        \Nnext@
        \next@\dimen@\copy\ZER@
        \ifN@\else\hskip-\wd\ZER@\fi
       \else
        \Nnext@
        \ifN@\else\hskip-\wd\ZER@\fi
        \next@\shifted@\copy\ZER@
        \ifN@\hskip-\wd\ZER@\fi
       \fi
       \ifnum\scount@=12
        \advance\shifted@\charht@\advance\goal@-\charht@
        \ifN@\hskip\wd\ZER@\else\hskip-\wd\ZER@\fi
       \fi
       \ifnum\scount@=13
        \getcos@{\thr@@\p@}%
        \ifN@\hskip\dimen@\else\hskip-\wd\ZER@\hskip-\dimen@\fi
        \adjust@\shifted@\advance\adjust@\dimen@ii
        \Nnext@
        \next@\adjust@\copy\ZER@
        \ifN@\hskip-\dimen@\hskip-\wd\ZER@\else\hskip\dimen@\fi
       \fi
      \else
       \ifN@\hskip-\wd\ZER@\fi
       \ifnum\scount@=\tw@
        \ifN@\hskip\wd\ZER@\else\hskip-\wd\ZER@\fi
        \dimen@\shifted@\advance\dimen@-\charht@
        \Nnext@
        \next@\dimen@\copy\ZER@
        \ifN@\hskip-\wd\ZER@\fi
       \else
        \Nnext@
        \next@\shifted@\copy\ZER@
        \ifN@\else\hskip-\wd\ZER@\fi
       \fi
       \ifnum\scount@=12
        \advance\shifted@\charht@\advance\goal@-\charht@
        \ifN@\hskip-\wd\ZER@\else\hskip\wd\ZER@\fi
       \fi
       \ifnum\scount@=13
        \getcos@{\thr@@\p@}%
        \ifN@\hskip-\wd\ZER@\hskip-\dimen@\else\hskip\dimen@\fi
        \adjust@\shifted@\advance\adjust@\dimen@ii
        \Nnext@
        \next@\adjust@\copy\ZER@
        \ifN@\hskip\dimen@\else\hskip-\dimen@\hskip-\wd\ZER@\fi
       \fi
      \fi
  \fi\fi\fi\fi
  \ifnum\arrcount@=\m@ne
  \else
   \loop
    \ifdim\goal@>\charht@
    \ifE@\else\hskip-\charwd@\fi
    \Nnext@
    \next@\shifted@\copy\shaft@
    \ifE@\else\hskip-\charwd@\fi
    \advance\shifted@\charht@\advance\goal@-\charht@
    \repeat
   \ifdim\goal@>\z@
    \dimen@\charht@\advance\dimen@-\goal@
    \divide\dimen@\tan@i\multiply\dimen@\tan@ii
    \ifE@\hskip-\dimen@\else\hskip-\charwd@\hskip\dimen@\fi
    \adjust@\shifted@\advance\adjust@-\charht@\advance\adjust@\goal@
    \Nnext@
    \next@\adjust@\copy\shaft@
    \ifE@\else\hskip-\charwd@\fi
   \else
    \adjust@\shifted@\advance\adjust@-\charht@
   \fi
  \fi
  \ifnum\arrcount@=\m@ne
  \else
   \ifnum\arrcount@=\thr@@
   \else
    \ifnum\tcount@=\m@ne
    \else
     \setbox\ZER@
      \hbox{\ifnum\angcount@=117 \arrow@v\else\arrowfont@\fi\char\Tcount@}%
     \ifnum\tcount@=\thr@@
      \advance\adjust@\charht@
      \ifE@\else\ifN@\hskip-\charwd@\else\hskip-\wd\ZER@\fi\fi
     \else
      \ifnum\tcount@=12
       \advance\adjust@\charht@
       \ifE@\else\ifN@\hskip-\charwd@\else\hskip-\wd\ZER@\fi\fi
      \else
       \ifE@\hskip-\wd\ZER@\fi
     \fi\fi
     \Nnext@
     \next@\adjust@\copy\ZER@
     \ifnum\tcount@=13
      \hskip-\wd\ZER@
      \getcos@{\thr@@\p@}%
      \ifE@\hskip-\dimen@\else\hskip\dimen@\fi
      \advance\adjust@-\dimen@ii
      \Nnext@
      \next@\adjust@\box\ZER@
     \fi
  \fi\fi\fi}}%
 \iflabel@i
  \rlap{\hskip\tocenter@
  \dimen@\firstx@\advance\dimen@\secondx@\divide\dimen@\tw@
  \advance\dimen@\ldimen@i
  \dimen@ii\firsty@\advance\dimen@ii\secondy@\divide\dimen@ii\tw@
  \global\multiply\ldimen@i\tan@i\global\divide\ldimen@i\tan@ii
  \ifNESW@\advance\dimen@ii\ldimen@i\else\advance\dimen@ii-\ldimen@i\fi
  \setbox\ZER@\hbox{\ifNESW@\else\ifnum\arrcount@=\thr@@\hskip4\p@\else
   \hskip\tw@\p@\fi\fi
   $\m@th\tsize@@\label@i$\ifNESW@\ifnum\arrcount@=\thr@@\hskip4\p@\else
   \hskip\tw@\p@\fi\fi}%
  \ifnum\arrcount@=\m@ne
   \ifNESW@\advance\dimen@.5\wd\ZER@\advance\dimen@\p@\else
    \advance\dimen@-.5\wd\ZER@\advance\dimen@-\p@\fi
   \advance\dimen@ii-.5\ht\ZER@
  \else
   \advance\dimen@ii\dp\ZER@
   \ifnum\slcount@<6 \advance\dimen@ii\tw@\p@\fi
  \fi
  \hskip\dimen@
  \ifNESW@\let\next@\llap\else\let\next@\rlap\fi
  \next@{\smash{\raise\dimen@ii\box\ZER@}}}%
 \fi
 \iflabel@ii
  \ifnum\arrcount@=\m@ne
  \else
   \rlap{\hskip\tocenter@
   \dimen@\firstx@\advance\dimen@\secondx@\divide\dimen@\tw@
   \ifNESW@\advance\dimen@\ldimen@ii\else\advance\dimen@-\ldimen@ii\fi
   \dimen@ii\firsty@\advance\dimen@ii\secondy@\divide\dimen@ii\tw@
   \global\multiply\ldimen@ii\tan@i\global\divide\ldimen@ii\tan@ii
   \advance\dimen@ii\ldimen@ii
   \setbox\ZER@\hbox{\ifNESW@\ifnum\arrcount@=\thr@@\hskip4\p@\else
    \hskip\tw@\p@\fi\fi
    $\m@th\tsize@\label@ii$\ifNESW@\else\ifnum\arrcount@=\thr@@\hskip4\p@
    \else\hskip\tw@\p@\fi\fi}%
   \advance\dimen@ii-\ht\ZER@
   \ifnum\slcount@<9 \advance\dimen@ii-\thr@@\p@\fi
   \ifNESW@\let\next@\rlap\else\let\next@\llap\fi
   \hskip\dimen@\next@{\smash{\raise\dimen@ii\box\ZER@}}}%
  \fi
 \fi
}
\def\outCD@#1{\def#1{\Err@{\noexpand#1must not be used within \string\CD}}}
\newskip\preCDskip@
\newskip\postCDskip@
\preCDskip@\z@
\postCDskip@\z@
\def\preCDspace#1{\RIfMIfI@
 \onlydmatherr@\preCDspace\else\advance\preCDskip@#1\relax\fi\else
 \onlydmatherr@\preCDspace\fi}
\def\postCDspace#1{\RIfMIfI@
 \onlydmatherr@\postCDspace\else\advance\postCDskip@#1\relax\fi\else
 \onlydmatherr@\postCDspace\fi}
\def\predisplayspace#1{\RIfMIfI@
 \onlydmatherr@\predisplayspace\else
 \advance\abovedisplayskip#1\relax
 \advance\abovedisplayshortskip#1\relax\fi
 \else\onlydmatherr@\preCDspace\fi}
\def\postdisplayspace#1{\RIfMIfI@
 \onlydmatherr@\postdisplayspace\else
 \advance\belowdisplayskip#1\relax
 \advance\belowdisplayshortskip#1\relax\fi
 \else\onlydmatherr@\postdisplayspace\fi}
\def\PreCDSpace#1{\global\preCDskip@#1\relax}
\def\PostCDSpace#1{\global\postCDskip@#1\relax}
\def\CD#1\endCD{%
 \outCD@\cgaps\outCD@\rgaps\outCD@\Cgaps\outCD@\Rgaps
 \preCD@#1\endCD
 \advance\abovedisplayskip\preCDskip@
 \advance\abovedisplayshortskip\preCDskip@
 \advance\belowdisplayskip\postCDskip@
 \advance\belowdisplayshortskip\postCDskip@
 \vcenter{\offinterlineskip
  \vskip\preCDskip@\Let@\global\colcount@\@ne\global\rowcount@\z@
  \everycr{%
   \noalign{%
    \ifnum\rowcount@=\Rowcount@
    \else
     \getrgap@\rowcount@\vskip\getdim@
     \global\advance\rowcount@\@ne\global\colcount@\@ne
    \fi}}%
  \tabskip\z@
  \halign{&\global\xoff@\z@\global\yoff@\z@
   \getcgap@\colcount@\hskip\getdim@
   \hfil\vrule\height\ten@\p@\width\z@\depth\z@
   $\m@th\displaystyle{##}$\hfil
   \global\advance\colcount@\@ne\cr
   #1\crcr}\vskip\postCDskip@}%
 \preCDskip@\z@\postCDskip@\z@
 \def\getcgap@##1{\ifcase##1\or\getdim@\z@\else\getdim@\standardcgap\fi}%
 \def\getrgap@##1{\ifcase##1\getdim@\z@\else\getdim@\standardrgap\fi}%
 \let\Width@\relax\let\Height@\relax\let\Depth@\relax\let\Rowheight@\relax
 \let\Rowdepth@\relax\let\Colwidth@\relax
}
\let\enddocument\bye
\def\alloc@#1#2#3#4#5{\global\advance\count1#1by\@ne
  \ch@ck#1#4#2%
  \allocationnumber=\count1#1%
  \global#3#5=\allocationnumber
  \wlog{\string#5=\string#2\the\allocationnumber}}
\catcode`\@=\active